\documentclass[11pt]{article}
\usepackage[T1]{fontenc}
\usepackage{lmodern}
\usepackage[margin=1in]{geometry}
\usepackage{amsmath,amssymb,amsthm,mathtools}
\usepackage{graphicx,booktabs,array}
\usepackage{microtype}
\usepackage{xcolor}
\usepackage{enumitem}
\usepackage[colorlinks=true,linkcolor=blue!45!black,citecolor=blue!45!black,filecolor=blue!45!black,urlcolor=blue!45!black]{hyperref}
\usepackage[nameinlink,noabbrev]{cleveref}
\crefname{theorem}{Theorem}{Theorems}
\crefname{lemma}{Lemma}{Lemmas}
\crefname{proposition}{Proposition}{Propositions}
\crefname{corollary}{Corollary}{Corollaries}
\crefname{remark}{Remark}{Remarks}
\crefname{conjecture}{Conjecture}{Conjectures}
\Crefname{conjecture}{Conjecture}{Conjectures}
\crefname{definition}{Definition}{Definitions}
\crefname{example}{Example}{Examples}
\crefname{section}{Section}{Sections}
\crefname{figure}{Figure}{Figures}
\crefname{table}{Table}{Tables}
\hypersetup{pdftitle={When Windowed Anderson Acceleration Converges Faster},pdfsubject={Convergence guarantees for Anderson acceleration with fixed memory},pdfkeywords={Anderson acceleration, polynomial certificate, moment equations, memory threshold, normaloid matrix, periodic Markov chain, deltoid}}
\setlist{itemsep=2pt,topsep=4pt}
\numberwithin{equation}{section}
\newtheorem{theorem}{Theorem}[section]
\newtheorem{lemma}[theorem]{Lemma}
\newtheorem{proposition}[theorem]{Proposition}
\newtheorem{corollary}[theorem]{Corollary}
\theoremstyle{definition}
\newtheorem{definition}[theorem]{Definition}
\newtheorem{example}[theorem]{Example}
\theoremstyle{remark}
\newtheorem{remark}[theorem]{Remark}

\DeclareMathOperator{\aff}{aff}
\DeclareMathOperator{\conv}{conv}
\DeclareMathOperator{\diag}{diag}
\DeclareMathOperator*{\argmin}{argmin}
\DeclareMathOperator{\spec}{\sigma}
\newcommand{\R}{\mathbb R}
\newcommand{\C}{\mathbb C}
\DeclarePairedDelimiter{\norm}{\lVert}{\rVert}
\DeclarePairedDelimiterX{\ip}[2]{\langle}{\rangle}{#1,#2}
\newcommand{\rhAA}{\varrho_{\mathrm{AA}}}

\title{When Windowed Anderson Acceleration Converges Faster}
\author{Gilad Lerman\thanks{\raggedright School of Mathematics, University of Minnesota, Minneapolis, MN 55455 USA (e-mail: lerman@umn.edu).}
\and Teng Zhang\thanks{\raggedright Department of Mathematics, University of Central Florida, Orlando, FL 32816 USA (e-mail: teng.zhang@ucf.edu).}}
\date{20 October 2026}

\begin{document}
\maketitle
\vspace{-1.5em}

\begin{abstract}
Windowed Anderson acceleration AA($m$) combines at most $m+1$ recent iterates
by least squares. For a real affine contraction $x\mapsto Wx+b$ with
$q=\norm W\in(0,1)$, we determine exactly when AA($m$) converges from every
starting point with an asymptotic factor below $q$. We give a polynomial
criterion for this improvement and an equivalent spectral condition: there
is no probability distribution on the phases of the eigenvalues of modulus
$q$ whose first $m$ cosine moments are $q,q^2,\ldots,q^m$. If such a
distribution exists, we construct a starting point from which AA($m$)
reproduces fixed-point iteration at every step. Otherwise, over a fixed
number of consecutive steps, every run contracts the residual norm by a
common factor strictly smaller than the corresponding power of $q$. The
polynomial criterion also gives explicit convergence-rate bounds with constants
independent of the dimension. When $q$ equals the spectral radius, this characterizes improvement over the asymptotic convergence factor of fixed-point iteration. 
For general affine contractions, the asymptotic factor never exceeds the numerical radius of $W$, and a related moment condition gives a uniform upper bound strictly below that radius.
\end{abstract}


\section{Introduction}
Anderson acceleration (AA)~\cite{anderson1965} combines previous
fixed-point iterates using coefficients chosen to minimize the norm
of the combined residual. Practical implementations usually retain
only a sliding window of at most $m+1$ recent iterates, keeping storage
requirements and least-squares problems small. Minimizing over this
window, however, need not improve the asymptotic convergence factor.
For an affine contraction $g(x)=Wx+b$, we characterize exactly when
this windowed method, denoted AA($m$), converges from every starting
point with an asymptotic factor strictly below the contraction factor
$\norm W$. This extends the sharp memory analysis for symmetric
contractions~\cite{garner2026} to general affine contractions.

Fix the inner product used in AA's least-squares problems (the
\emph{working inner product}), and let $q=\norm W\in(0,1)$ in its induced
norm. With residuals $r_k=g(x_k)-x_k$, both fixed-point iteration and
AA($m$) multiply the residual norm by at most $q$ per step. We ask when
every AA($m$) run has an asymptotic factor strictly below $q$, and then say
that AA($m$) \emph{beats} $q$. When $q$ equals the spectral radius
$\rho(W)$, beating $q$ means improving on the worst-case fixed-point
factor. Such matrices are called \emph{normaloid}~\cite{halmos1982} and
include all matrices normal in the working inner product. When
$\rho(W)<q$, fixed-point iteration itself beats $q$;
\cref{main:sec:numrad} studies comparison with the smaller benchmark
$\rho(W)$. Throughout, normality and self-adjointness refer to the working
inner product. 
Our main result shows that the eigenvalues of modulus $q$ alone
determine whether AA($m$) beats $q$ from every starting point.

\paragraph{An exact alternative.}
For real $W$, list its distinct eigenvalues of modulus $q$, up to
complex conjugacy, as $qe^{i\theta_1},\ldots,qe^{i\theta_{K_q}}$, with
$0\le\theta_i\le\pi$. Every AA($m$) run beats $q$ if and only if there
are no weights satisfying
\[
 w_i\ge0,\qquad \sum_{i=1}^{K_q}w_i=1,\qquad
 \sum_{i=1}^{K_q}w_i\cos(\ell\theta_i)=q^\ell
 \quad(1\le\ell\le m).
\]
These equations ask whether a probability distribution on the phases can
have first $m$ cosine moments $q,q^2,\ldots,q^m$. If they have a solution,
its weights determine an explicit initial residual for which AA
reproduces fixed-point iteration exactly, with
$\norm{r_k}=q^k\norm{r_0}$. Otherwise, a fixed number of consecutive steps
reduces the residual norm by a common factor strictly smaller than the
corresponding power of $q$, for every start
(\cref{main:thm:normalalternative}).

The equivalent polynomial criterion is the existence of a real
polynomial $P(z)=\sum_{j=s}^{d}a_jz^j$, with integers $d\ge1$ and
$0\le s\le d$, such that $P(1)=1$, $d-s\le m$, and
$\norm{P(W)}<q^d$. We call such a polynomial a \emph{certificate}.
The normalization makes its coefficients sum to one, as required for
AA's residual combinations, while the norm inequality improves on
the bound for $d$ fixed-point steps. The window restricts the width
$d-s$ of the polynomial's support rather than its degree. For normal
matrices one can take $s=0$; general contractions may require the
factor $z^s$, which delays the action of the remaining polynomial.

For self-adjoint contractions, the only possible phases are $0$ and
$\pi$. The moment equations recover the sharp memory threshold~\cite{garner2026}:
AA(1) beats $q$ unless both $q$ and $-q$ are eigenvalues, while AA(2)
always beats $q$.

\paragraph{Sharp memory thresholds.}
When $q=\rho(W)$, let $N_q$ be the number of distinct eigenvalues of
modulus $q$. Memory $N_q$ always beats $q$, and this bound is sharp for
every $N_q\ge2$ (\cref{main:prop:peripheralcount}). For an integer
$\nu\ge2$, suppose these eigenvalues are exactly the $\nu$ points
$qe^{2\pi ij/\nu}$, $0\le j<\nu$. Then the smallest memory that beats
$q$ is $m_{\min}=\lfloor\nu/2\rfloor+1$. Aliasing explains the upper
bound: moments of orders $\ell$ and $\nu-\ell$ agree on these phases,
whereas $q^\ell$ and $q^{\nu-\ell}$ differ when the orders differ
(\cref{main:thm:cusps}).

\begin{figure}[t]
\centering
\includegraphics[width=0.92\linewidth]{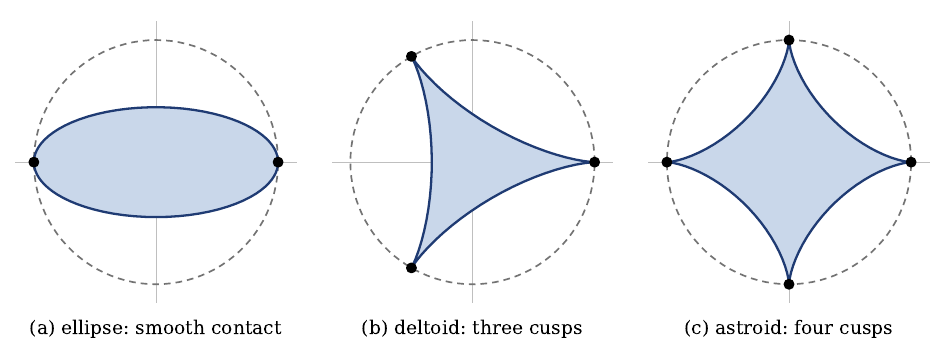}
\par\vspace{-3pt}{\footnotesize
\makebox[0.3067\linewidth]{memory~threshold~$2$}
\makebox[0.3067\linewidth]{memory~threshold~$2$}
\makebox[0.3067\linewidth]{memory~threshold~$3$}}
\caption{The shaded regions represent possible locations of the spectrum:
an ellipse, the deltoid, and the astroid, shown at unit scale. Black
dots mark contacts with the unit circle. For matrices with
$\norm W=q\in(0,1)$ and spectrum in the corresponding $q$-scaled
region, the displayed memory thresholds guarantee improvement below
$q$ from every starting point and are sharp over each matrix class
(\cref{main:thm:cusps}). The ellipse shares its contacts $\pm1$ with
$[-1,1]$. For normal matrices, rate bounds over each region also depend
on the boundary geometry near these contacts
(\cref{main:sec:sqrtstatus}).}
\label{main:fig:regions}
\end{figure}

The same memory bound remains sufficient when the eigenvalues of
modulus $q$ occupy only some of the $\nu$ equally spaced points.
Consequently, it also applies when the spectrum lies in a region
contained in $|z|\le q$ that meets $|z|=q$ only at those points.
\Cref{main:fig:regions} shows three such regions at unit scale:
an ellipse with contacts $\pm1$, and the deltoid and astroid with
three and four equally spaced contacts, respectively. After scaling
by $q$, these regions give sharp worst-case memory thresholds of
$2$, $2$, and $3$. The ellipse shares its contacts with $[-1,1]$.
Its smooth boundary provides a comparison with the cusped regions
when we study how the boundary geometry affects quantitative
convergence bounds for normal matrices
(\cref{main:sec:sqrtstatus}).

For an irreducible Markov chain of period $\nu$, multiplying the
transition matrix by $q$ gives exactly these eigenvalues of modulus $q$.
Its stationary distribution defines an inner product in which this
matrix is normaloid, whether the chain is reversible or not. The same
sharp threshold therefore applies to discounted policy evaluation
(\cref{main:ex:markov}).

\paragraph{Rates and their limits.}
A certificate bounds the contraction over a fixed number of consecutive
steps, with constants independent of the dimension and initialization.
AA does not apply $P(W)$ directly to a residual. Instead, if these steps
barely improve on $q$, the least-squares corrections must be small, so
the stored residuals nearly follow fixed-point iteration. The certificate
then supplies a feasible combination that forces a strict gain
(\cref{main:thm:certificate,main:prop:delayedcertificate}). For a
certificate $p$ of degree $d\le m$, a sharper estimate gives a
guaranteed gain linear in $1-\norm{p(W)}/q^d$, with constants depending
on $q$, $d$, and the coefficients, but not on the dimension,
initialization, or memory $m\ge d$ (\cref{main:thm:generalmargin}).

At every fixed memory $m\ge2$ there are matrices with a quadratic certificate
for which the asymptotic factor is only $O(1-q)$ below $q$ as
$q\uparrow1$. For normal matrices with spectra in an ellipse inside
$|z|\le q$ and tangent at $q$, polynomial approximation gives the same
limitation for dimension-independent finite-step bounds. For the cusped
regions considered here, stationary methods with a fixed window attain
gains of order $\sqrt{1-q}$; whether AA with fixed memory attains them
remains open (\cref{main:sec:limits}).

\paragraph{Beyond the contraction factor.}
When $\rho(W)<q$, AA(1) can be slower even than the worst-case
fixed-point factor $\rho(W)$ (\cref{main:prop:notbenchmark}). The
numerical range consists of the values $\ip{z}{Wz}$ for unit complex
vectors $z$, using the complex extension of the working inner product;
its largest modulus is the numerical radius $w(W)$. Every AA($m$) run
has asymptotic factor at most $w(W)$, with
$\rho(W)\le w(W)\le q$ and $w(W)<q$ when $\rho(W)<q$
(\cref{main:thm:numradceiling}). A related moment test uses the angles
where the numerical range touches the circle of radius $w(W)$.
Infeasibility of its moment equations guarantees a uniform factor below
$w(W)$ (\cref{main:thm:uniformceiling}). For normaloid matrices this
recovers the exact criterion above. In general it is a sufficient
condition; whether AA always improves uniformly on $w(W)$ when
$w(W)>\rho(W)$ remains open.

\paragraph{Inner products.}
The working inner product is part of the algorithm: changing it changes
the least-squares coefficients, $q$, and $w(W)$, but leaves $\rho(W)$
unchanged. If the eigenvalues of modulus $\rho(W)>0$ are semisimple,
some inner product makes $W$ normaloid, so the exact criterion compares
AA with the worst-case fixed-point factor (supplementary
\cref{si:cor:spectralbenchmark}). Constructing such an inner product
usually requires spectral information; Markov chains provide a natural
exception.

\paragraph{Nonlinear maps.}
If $g$ is differentiable at a fixed point $x_*$ and $g'(x_*)$ is
contractive and admits a certificate of degree at most $m$, we also
obtain local rate bounds. For this extension, we constrain the sum of
the absolute values of the mixing coefficients by a fixed $\Lambda>1$,
which can be chosen independently of the certificate coefficients.
The rate bounds depend on $\Lambda$. No Lipschitz assumption on the
derivative is required (\cref{main:thm:nonlinear}).

\paragraph{Relation to earlier work.}
Connections with multisecant methods and Shanks-type transformations
are developed in~\cite{fang2009,brezinski2022}. For affine problems,
AA with full memory, which retains all previous iterates, is closely
connected to the generalized minimal residual method
(GMRES)~\cite{walker2011,potra2013}. For affine contractions, the
nonstagnation hypothesis needed for this connection follows from
the argument of~\cite{eisenstat1983}. For finite windows, polynomial
representations $r_k=p_k(W)r_0$ and the iteration on residual windows
are analyzed in~\cite{desterck2024,desterck2022}. Our starting point is the symmetric
analysis of Garner, Lerman, and Zhang~\cite{garner2026}: we use its
observation that little residual improvement forces small least-squares
corrections to obtain a polynomial test, quantitative bounds, and an
exact moment obstruction for nonsymmetric contractions. Weighted
formulations appear in~\cite{yang2022weighted}; restarted and periodic
variants in~\cite{henderson2019,bertrand2021,banerjee2016,krzysik2025}.

Cowal, Marshall, and Pollock~\cite{cowal2026deltoid,cowal2026random}
construct rapidly growing polynomials from mean-zero random walks and
use them for momentum power methods. Their comparison of polynomial
growth near cusps and smooth boundaries underlies our discussion of
the possible size of AA's gain. These polynomials connect to Faber
polynomials~\cite{eiermann1989,eiermann1993} and classical semiiterative
methods~\cite{niethammervarga1983,eiermann1985}; our question concerns
AA's least-squares choice of coefficients. For nonlinear AA,
coefficient bounds and least-squares gains are central
in~\cite{toth2015,evans2020,pollock2021}; constrained extrapolation and
regularization are analyzed
in~\cite{barre2022,scieur2016,scieur2018,ouyang2023}.

\paragraph{Organization.}
After the projection identities in \cref{main:sec:geometry},
\cref{main:sec:certificate} proves the certificate bound,
\cref{main:sec:normal} the exact criterion and memory thresholds, and
\cref{main:sec:numrad} the numerical-radius bounds.
\Cref{main:sec:limits} treats quantitative gains and their limits,
and \cref{main:sec:nonlinear} the nonlinear extension.
\Cref{main:sec:conclusion} summarizes the results and open questions.
The supplement contains the remaining proofs and refinements.

\section{Windowed Anderson acceleration}\label{main:sec:geometry}

\paragraph{Notation and problem setup.}
Throughout, we fix an inner product $\ip\cdot\cdot$ on $\R^n$, the \emph{working inner product}, in which AA solves its least-squares problems; $\norm\cdot$ denotes the norm it induces on vectors and the corresponding operator norm on matrices. The working inner product is the Euclidean one unless
stated otherwise; the change of coordinates at the end of this section shows that this is no loss of generality. Further, $|\cdot|$ is the absolute value of a real
or complex number, $\spec(W)$ and $\rho(W)$ denote the spectrum and spectral
radius, and $\aff$ and $\conv$ denote affine and convex hulls. We write $A^T$ for the transpose and $A^*$ for the conjugate transpose of a matrix or vector, and $\bar\lambda$ for the complex conjugate of a number; for a real matrix $W$, $W^*=W^T$, which is also the adjoint of $W$ on $\C^n$ in the Euclidean inner product. Finally, $\norm\cdot_2$ denotes the Euclidean norm and the corresponding operator norm, and $R_\theta$ the rotation of the plane by the angle $\theta$.
Anderson acceleration accelerates a fixed-point iteration $x_{k+1}=g(x_k)$ for a map $g\colon\R^n\to\R^n$, using the residual $r(x)=g(x)-x$; the algorithm below needs only values of $g$. \Cref{main:sec:certificate,main:sec:normal,main:sec:numrad,main:sec:limits} treat the affine map
\[
 g(x)=Wx+b,\qquad W\in\R^{n\times n},\qquad q=\norm W\in(0,1),
\]
whose fixed point and residual are
\[
x_*=(I-W)^{-1}b, \qquad r(x)=(W-I)(x-x_*);
\]
the excluded case $W=0$ is solved by one fixed-point step. \Cref{main:sec:nonlinear} returns to nonlinear maps, with the derivative of $g$ at its fixed point in the role of $W$.

\paragraph{The sliding-window algorithm.}
For a memory $m\ge1$, also called the depth, set $m_k=\min\{m,k\}$. Starting
from a single point $x_0$, Anderson acceleration AA($m$) computes
\begin{equation}\label{main:eq:algorithm}
 (\alpha_j)_{j=k-m_k}^k\in
 \argmin_{\sum_{j=k-m_k}^k\alpha_j=1}
       \norm*{\sum_{j=k-m_k}^k\alpha_jr_j},\qquad r_j=r(x_j),\qquad
 x_{k+1}=\sum_{j=k-m_k}^k\alpha_jg(x_j).
\end{equation}
Thus AA($m$) retains at most $m+1$ residuals, and its first step is the
fixed-point step $x_1=g(x_0)$. 
\emph{Capped} AA($m$) with a cap $\Lambda\ge1$ adds the constraint $\sum_{j=k-m_k}^k|\alpha_j|\le\Lambda$ to the minimization in \eqref{main:eq:algorithm}. It is a variant of the modified AA($m$) of~\cite[(3.5)]{garner2026}, and is used for nonlinear maps in \cref{main:sec:nonlinear}. Where the distinction matters, we call the unconstrained method \emph{ordinary} AA($m$).

For ordinary AA, define $y_{k+1}$ as the orthogonal projection of the origin onto the affine hull $\aff\{r_{k-m_k},\ldots,r_k\}$ of the retained residuals; this is exactly the minimum-norm point within that affine space. 
Since $g$ is affine and
$\sum_j\alpha_j=1$, the new residual is
$r(\sum_j\alpha_jg(x_j))=\sum_j\alpha_jWr_j$, as
in~\cite[Section~3.2]{garner2026}. With the orthogonality of the
projection, and since $r_k$ lies in the affine hull, this gives
\begin{align}
 r_{k+1}&=Wy_{k+1},\label{main:eq:residualupdate}\\
 \ip{y_{k+1}}{r_j-y_{k+1}}&=0
       \qquad(k-m_k\le j\le k),\label{main:eq:orthogonality}\\
 \norm{r_k-y_{k+1}}^2
   &=\norm{r_k}^2-\norm{y_{k+1}}^2,\label{main:eq:pythagoras}\\
 \norm{r_{k+1}}&\le q\norm{y_{k+1}}\le q\norm{r_k}.
       \label{main:eq:basiccontraction}
\end{align}
The projection $y_{k+1}$ is unique by the geometry of Hilbert spaces, which makes the new residual $r_{k+1} = Wy_{k+1}$ unique. Because $\norm W<1$, the matrix $W-I$ is invertible, so the residual $r(x)=(W-I)(x-x_*)$ determines the iterate; hence $x_{k+1}$ is also unique, even if the minimizing coefficients $(\alpha_j)_{j=k-m_k}^k$ are not.

\paragraph{Asymptotic rates.}
We measure
the speed by the root-linear factor
\(\rhAA(x_0)=\limsup_{k\to\infty}\norm{x_k-x_*}^{1/k}\), which vanishes for
a terminating run. 
The inequalities
\begin{equation}\label{main:eq:errorequivalence}
 (1-q)\norm{x_k-x_*}\le\norm{r_k}\le(1+q)\norm{x_k-x_*}
\end{equation}
show that $\limsup_{k\to\infty}\norm{r_k}^{1/k} = \limsup_{k\to\infty}\norm{x_k-x_*}^{1/k}$, so the exact same asymptotic factor is obtained by measuring the residuals. 
The factor of
fixed-point (Picard) iteration from the same initial point is $\rho_{\mathrm{Picard}}(x_0)=\limsup_{k\to\infty} \norm{W^k(x_0-x_*)}^{1/k}\le\rho(W)$, 
with equality unless the initial error $x_0-x_*$ has no components in the directions of the dominant eigenvalues (i.e., it lies entirely in the invariant subspace of the eigenvalues of modulus strictly below $\rho(W)$).

\paragraph{The inner product is part of the algorithm.}
Every inner product on $\R^n$ takes the form $\ip xy_H=x^THy$ for a symmetric positive definite matrix $H$. By factoring $H=S^TS$ for some invertible $S$ (for instance, a Cholesky factor or the square root of $H$), the induced norm becomes $\norm x_H=\norm{Sx}_2$. Write AA$_H$($m$) for the iteration \eqref{main:eq:algorithm} with its least-squares problem posed in $\norm\cdot_H$; the update $x_{k+1}=\sum_j\alpha_jg(x_j)$ remains algebraically identical. However, for $H\ne I$, this operates as a different algorithm from Euclidean AA($m$), and they are related by a change of coordinates.

Precisely, $(x_k)$ is an AA$_H$($m$) sequence for $g$ if and only if
$(Sx_k)$ is a Euclidean AA($m$) sequence for the conjugated map
$\tilde g(\tilde x)=Sg(S^{-1}\tilde x)$, generating the same minimizing
coefficients, scaled residuals $\tilde r_k=Sr_k$, and identical root-linear factors. For $g(x)=Wx+b$, the conjugated matrix is $SWS^{-1}$, and $\norm{p(W)}_H=\norm{Sp(W)S^{-1}}_2$ for every polynomial $p$ (supplementary \cref{si:prop:conjugation}).

All results in this paper are statements about AA in the working inner product, and we write their proofs for the Euclidean one. By supplementary \cref{si:prop:conjugation}, each holds verbatim for AA$_H$ when the norm, the adjoint $W^{*_H}=H^{-1}W^TH$, orthogonality, and normality ($W^{*_H}W=WW^{*_H}$) are taken in the inner product of $H$; the adjoint has the same form on $\C^n$, with $\ip xy_H=x^*Hy$. A statement about $W$ in that inner product is simply the same statement applied to $SWS^{-1}$ in the Euclidean inner product (and likewise for the derivative $g'(x_*)$ of a nonlinear map). In particular, a normal matrix is one that is normal in the working inner product, ensuring $q=\rho(W)$. Normality is therefore a joint property of $W$ and the inner product; because every matrix diagonalizable over $\C$ is normal in \emph{some} inner product, supplementary \cref{si:sec:spectralradius} uses this freedom to compare AA with fixed-point iteration.

\section{A polynomial test for faster convergence}\label{main:sec:certificate}
A \emph{certificate} for AA($m$) is a real polynomial
$P(z)=\sum_{j=s}^d a_jz^j$ with $P(1)=1$, $d-s\le m$, and
$\norm{P(W)}<q^d$. It is \emph{undelayed} when $s=0$, so that
$d\le m$, and we then write $p=P$. Otherwise it is \emph{delayed}:
$P(z)=z^sp(z)$ with $\deg p\le d-s\le m$. The window can retain
$r_{k+s},\ldots,r_{k+d}$, so $\sum_j a_jr_{k+j}$ is a feasible
least-squares combination. If these residuals were exactly $W^jr_k$,
that combination would equal $P(W)r_k$. When a block gains little over
fixed-point iteration, its least-squares corrections are small, making
this comparison accurate enough to force a gain.

For an undelayed certificate, set $\eta=\norm{p(W)}/q^d$ and call
$1-\eta$ its \emph{margin}. The next theorem converts this margin into
an explicit block contraction, with constants depending only on $q$
and the certificate, independently of the dimension, the start, and
the memory $m\ge d$. \Cref{main:prop:delayedcertificate} extends the
bound to delayed certificates. These estimates give a gain quadratic
in the margin; \cref{main:sec:lineargain} proves a sharper bound that
is linear in the margin for undelayed certificates.

\begin{theorem}[Polynomial certificate]\label{main:thm:certificate}
Let $p(z)=\sum_{j=0}^d a_jz^j$ be an undelayed certificate for AA($m$): its coefficients are real, $1\le d\le m$, $p(1)=1$, and $\norm{p(W)}<q^d$. Then the number 
\begin{equation}\label{main:eq:constants}
 c=\frac{\eta+\beta\sqrt{1+\beta^2-\eta^2}}{1+\beta^2},\qquad\text{where}\quad
 \eta=\frac{\norm{p(W)}}{q^d},\quad
 \beta=\frac1{q^d}\Bigl\|\Bigl(\sum_{j=i+1}^d|a_j|q^j\Bigr)_{i=0}^{d-1}\Bigr\|_2,
\end{equation}
lies in $(0,1)$ and is the largest $u\in[0,1]$ with $u\le\eta+\beta\sqrt{1-u^2}$. Every AA($m$) trajectory satisfies
\begin{equation}\label{main:eq:blockbound}
 \norm{r_{k+d+1}}\le cq^{d+1}\norm{r_k}\qquad(k\ge0).
\end{equation}
Consequently,
\begin{equation}\label{main:eq:rootbound}
 \rhAA(x_0)\le q c^{1/(d+1)}<q.
\end{equation}
\end{theorem}
The factor $c$ always lies between $\eta$ and $1$: it is close to $\eta$ when the coefficients of $p$ are small, that is, when $\beta$ is small, and it tends to $1$ as they grow. No normality of $W$ and no linear independence of the stored residuals are assumed.

The proof extends the residual-space argument of Garner, Lerman, and Zhang~\cite{garner2026} for symmetric matrices. Let $t=\norm{r_{k+d+1}}/(q^{d+1}\norm{r_k})\in[0,1]$ be the factor that a block of $d+1$ steps achieves relative to the guarantee $q^{d+1}$. If the stored residuals were exactly $r_k,Wr_k,\ldots,W^dr_k$, the certificate's combination $\sum_ja_jr_{k+j}$ would equal $p(W)r_k$, and the least-squares property would give $t\le\eta$. The residuals differ from these powers only through the least-squares corrections $r_{k+j}-y_{k+j+1}$, whose squared norms, weighted by $q^{-2j}$, add up to at most $(1-t^2)\norm{r_k}^2$; this is the energy that the block releases. A block with little gain, $t$ near $1$, therefore has small corrections, its residuals are nearly the powers $W^jr_k$, and the certificate forces a gain after all. Balancing the two effects gives the inequality $t\le\eta+\beta\sqrt{1-t^2}$, in which $\beta$ weights the corrections by the coefficients of $p$, and hence $t\le c$.
\begin{proof}

If $r_k=0$, all later residuals vanish. Otherwise, by homogeneity
(scaling $x_0-x_*$ scales every residual), we may assume $\norm{r_k}=1$.
By \eqref{main:eq:pythagoras} and~\eqref{main:eq:basiccontraction}, each correction
satisfies $\norm{r_{k+j}-y_{k+j+1}}^2\le\norm{r_{k+j}}^2-q^{-2}\norm{r_{k+j+1}}^2$.
Weighting by $q^{-2j}$, summing over $j<d$, and using
$\norm{r_{k+d+1}}\le q\norm{r_{k+d}}$ and $t$ defined above give
\begin{equation}\label{main:eq:telescoping}
 \sum_{j=0}^{d-1}q^{-2j}\norm{r_{k+j}-y_{k+j+1}}^2\le1-t^2,
 \qquad t=q^{-d-1}\norm{r_{k+d+1}}\in[0,1].
\end{equation}
Since $r_{k+j+1}=Wy_{k+j+1}$, induction on $j$ gives $r_{k+j}-W^jr_k=-\sum_{i=0}^{j-1}W^{j-i}(r_{k+i}-y_{k+i+1})$.
Grouping the corrections by their time index gives $\sum_{j=0}^da_jr_{k+j}-p(W)r_k =-\sum_{i=0}^{d-1}\sum_{j=i+1}^d a_jW^{j-i}(r_{k+i}-y_{k+i+1})$.
The norm of this discrepancy is at most $\sum_{i=0}^{d-1}(\sum_{j=i+1}^d|a_j|q^j) q^{-i}\norm{r_{k+i}-y_{k+i+1}} \le q^d\beta\sqrt{1-t^2}$,
by Cauchy--Schwarz and \eqref{main:eq:telescoping}.
At iteration $k+d$, the window contains $r_k,\ldots,r_{k+d}$,
including during startup, because $m\ge d$. Since $\sum_ja_j=1$, the
combination $\sum_ja_jr_{k+j}$ is feasible, so $\norm{y_{k+d+1}} \le\norm{p(W)r_k} +\norm*{\sum_{j=0}^da_jr_{k+j}-p(W)r_k} \le\norm{p(W)}+q^d\beta\sqrt{1-t^2}$.
On the other hand, $\norm{y_{k+d+1}}\ge q^{-1}\norm{r_{k+d+1}}=q^dt$.
Thus, after division by $q^d$, $t\le\eta+\beta\sqrt{1-t^2}$.
Here $\beta>0$: otherwise $p\equiv1$, contrary to $\norm{p(W)}<q^d$. The function
$u\mapsto u-\eta-\beta\sqrt{1-u^2}$ is strictly increasing on $[0,1]$,
negative at $u=0$, and positive at $u=1$. Squaring the equation
$u-\eta=\beta\sqrt{1-u^2}$ shows that its zero is the number $c$ in
\eqref{main:eq:constants}. Hence $0<c<1$ and $t\le c$, proving the block
bound. Iterating blocks and using \eqref{main:eq:basiccontraction}
between block endpoints proves the root bound.
\end{proof}

\Cref{main:thm:certificate} uses $d\le m$ only so that the combination $\sum_ja_jr_{k+j}$ is available to the least-squares problem. Only the residuals with nonzero coefficients are needed, so the theorem extends to delayed certificates, with the same constant.

\begin{proposition}[Delayed polynomial certificates]\label{main:prop:delayedcertificate}
Let $P(z)=\sum_{j=s}^{d}a_jz^j=z^sp(z)$ be a certificate for AA($m$), possibly delayed, so that $d-s\le m$, $P(1)=1$, and $\norm{P(W)}<q^d$, which forces $d\ge1$; the norm condition is on $P$, not on $p$.
Then \eqref{main:eq:blockbound} and~\eqref{main:eq:rootbound} hold with the constant $c$ of \cref{main:thm:certificate}, computed from $P$ with
$a_j=0$ for $j<s$.
\end{proposition}
\begin{proof}
In the proof of \cref{main:thm:certificate}, the final competitor $\sum_{j=0}^d a_jr_{k+j}$ now involves only $r_{k+s}$, $\ldots$, $r_{k+d}$, which are retained because $d-s\le m$, including during startup. All preceding estimates are unchanged; in particular $\beta>0$, since $\beta=0$ would force $P\equiv a_0=1$.
\end{proof}

The factor $z^s$ is what nonnormal matrices may need. The eigenvalues of modulus $q$ split off a block, orthogonal to the rest, on which $W$ is $q$ times an isometry, as in \eqref{main:eq:peripheralsplit}; the remaining block has spectral radius below $q$, but its norm may equal $q$, and $z^s$ lets it decay before $p$ acts. For example, let $W=\diag(q,qJ)$ with $J=\bigl[\begin{smallmatrix}0&1\\0&0\end{smallmatrix}\bigr]$. No certificate of degree one exists. Indeed, a certificate $p(z)=a_0+a_1z$ would satisfy $|p(q)|\le\norm{p(W)}<q$, because $q$ is an eigenvalue of $W$, and together with $a_0+a_1=1$ this forces $a_1>1$; but then the off-diagonal entry $a_1q$ of $p(qJ)$ gives $\norm{p(W)}\ge\norm{p(qJ)}>q$. The delayed polynomial $z^2(z-q)/(1-q)$, on the other hand, annihilates $W$, because $J^2=0$; it is a delayed certificate with $s=2$ and $d=3$, usable already with memory one. A certificate of degree at most $m$ may therefore fail to exist even when every trajectory beats $q$; for normal $W$ this cannot happen (\cref{main:thm:normalalternative}).

Certificates are sufficient conditions with explicit rates. The next section shows that, once delays are allowed, they are also necessary: AA($m$) beats $q$ from every start exactly when a delayed certificate exists. \Cref{main:sec:lineargain} returns to undelayed certificates and sharpens the gain to one linear in the margin.

\section{An exact criterion and sharp memory thresholds}\label{main:sec:normal}
This section proves that certificates are necessary as well as
sufficient for AA($m$) to beat $q=\norm W$ from every start. The
criterion depends only on the eigenvalues of modulus $q$, which we
call the \emph{peripheral eigenvalues}. If there are none, the criterion
holds automatically. Otherwise, moment equations on their phases give
the exact obstruction. We then use these equations to determine sharp
memory thresholds for general configurations and for the equally
spaced phases forced by periodicity.

\subsection{When a uniform improvement is possible}
Let $W$ be real with $q=\norm W\in(0,1)$. List its distinct eigenvalues of modulus $q$, up to complex conjugacy, as $qe^{i\theta_1},\ldots,qe^{i\theta_{K_q}}$, with phases $0\le\theta_i\le\pi$; thus $K_q$ is the number of these eigenvalues with nonnegative imaginary part, and $K_q=0$ exactly when $\rho(W)<q$. The criterion rests on a splitting of $W$, the finite-dimensional case of the canonical decomposition of a contraction into a unitary and a completely nonunitary part~\cite[Theorem~I.3.2]{sznagy2010}. Let $M$ be the real span of the real and imaginary parts of the eigenvectors of $W$ for its eigenvalues of modulus $q$. If $Wx=\lambda x$ with $|\lambda|=q$ and $\norm x=1$, then $\norm{W^Tx-\bar\lambda x}^2\le q^2-2|\lambda|^2+|\lambda|^2=0$, so these eigenvectors are shared with $W^T$. Hence $M$ is invariant under $W$ and $W^T$, and so is $M^\perp$, and in an orthonormal basis adapted to $\R^n=M\oplus M^\perp$,
\begin{equation}\label{main:eq:peripheralsplit}
 W=Q\diag(qU,W_0)Q^T,\qquad Q,U\text{ orthogonal},\qquad
 \norm{W_0}\le q,\qquad \rho(W_0)<q,
\end{equation}
where either summand may be absent; $M=\{0\}$ exactly when $\rho(W)<q$. Supplementary \cref{si:lem:normattained} states this precisely, together with the further properties used in the proofs.

\begin{theorem}[Polynomial and moment alternatives]\label{main:thm:normalalternative}
For every real $W$ with $q=\norm W\in(0,1)$ and every $m\ge1$, the following
are equivalent:
\begin{enumerate}[label=(\roman*)]
\item There are integers $d\ge1$ and $0\le s\le d$, with $d-s\le m$,
      and a real polynomial $P(z)=\sum_{j=s}^{d}a_jz^j$ such that
      $P(1)=1$ and $\norm{P(W)}<q^d$.
\item There are no weights $w_i\ge0$ with $\sum_{i=1}^{K_q}w_i=1$ and
      \begin{equation}\label{main:eq:momentcondition}
       \sum_{i=1}^{K_q}w_i\cos(\ell\theta_i)=q^\ell,\qquad 1\le\ell\le m.
      \end{equation}
      In particular, (ii) holds when $\rho(W)<q$, that is, when $W$ is not normaloid.
\item There are an integer $k_0\ge1$ and a constant $c\in[0,1)$,
      independent of the initialization, such that every AA($m$)
      trajectory satisfies $\norm{r_{k+k_0}}\le c q^{k_0}\norm{r_k}$ for
      every $k\ge0$.
\item For every starting point $x_0$, the AA($m$) trajectory has root-linear factor $\rhAA(x_0)=\limsup_k\norm{x_k-x_*}^{1/k}<q$.
\end{enumerate}
If these conditions fail, there is a real initial residual $r_0\in M$ of
norm one from which AA($m$) follows fixed-point iteration exactly, with
$r_k=W^kr_0$ and $\norm{r_k}=q^k$ for all $k\ge0$.

If, moreover, the block $W_0$ in \eqref{main:eq:peripheralsplit} is absent or has norm below $q$, as for every normal $W$, then (i)--(iv) are also equivalent to
\begin{enumerate}[label=(\roman*),start=5]
\item there is a real polynomial $p$ of degree at most $m$ with $p(1)=1$
      and $\norm{p(W)}<q^m$,
\end{enumerate}
and (iii) holds with $k_0=m+1$.
\end{theorem}
\begin{proof}
We prove (ii)$\Rightarrow$(i)$\Rightarrow$(iii)$\Rightarrow$(iv)$\Rightarrow$(ii). The last implication constructs the stated fixed-point trajectory when the moment equations are feasible. We treat (v) at the end.

\emph{(ii)$\Rightarrow$(i).} Suppose (ii) holds. If $\rho(W)<q$, then $\norm{W^d}<q^d$ for some $d$,
and the monomial $P(z)=z^d$, with $s=d$, proves (i). Otherwise $\rho(W)=q$.
Since $Q$ in \eqref{main:eq:peripheralsplit} is orthogonal and $P(W)=Q\diag\bigl(P(qU),P(W_0)\bigr)Q^T$, we have $\norm{P(W)}=\max\{\norm{P(qU)},\norm{P(W_0)}\}$ for every polynomial $P$.
The origin lies outside the convex hull of the vectors
$\bigl(\cos(\ell\theta_i)-q^\ell\bigr)_{\ell=1}^m$, so strict separation
gives $h\in\R^m$ with
$\sum_{\ell=1}^mh_\ell(\cos(\ell\theta_i)-q^\ell)<0$ for every $i$. Set $p_\varepsilon(z)=z^m+\varepsilon\sum_{\ell=1}^{m}h_\ell q^\ell(z^{m-\ell}-z^m)$, so that $p_\varepsilon(1)=1$. At $\lambda_i=qe^{i\theta_i}$,
\[
 |p_\varepsilon(\lambda_i)|^2=q^{2m}+2\varepsilon q^{2m}\sum_{\ell=1}^{m}h_\ell\bigl(\cos(\ell\theta_i)-q^\ell\bigr)+O(\varepsilon^2),
\]
and the middle term is negative. As $qU$ is normal, $\norm{p_\varepsilon(qU)}$ is the largest value of $|p_\varepsilon|$ on its eigenvalues, which are the $\lambda_i$ and their conjugates, and $|p_\varepsilon(\overline{\lambda_i})|=|p_\varepsilon(\lambda_i)|$ because $p_\varepsilon$ has real coefficients. Hence $\norm{p_\varepsilon(qU)}<q^m$ for every sufficiently small $\varepsilon>0$; fix such an $\varepsilon$. If $W_0$ is absent, $P=p_\varepsilon$, with $s=0$
and $d=m$, proves (i). Otherwise, since $\rho(W_0)<q$, Gelfand's formula $\norm{W_0^s}^{1/s}\to\rho(W_0)$ gives $\norm{W_0^s}\le(\rho(W_0)+\delta)^s$ for every $\delta>0$ and all large $s$, so for sufficiently
large $s$ we have $\norm{W_0^s p_\varepsilon(W_0)}\le\norm{W_0^s}\,\norm{p_\varepsilon(W_0)}<q^{s+m}$, while $\norm{(qU)^sp_\varepsilon(qU)}=q^s\norm{p_\varepsilon(qU)}<q^{s+m}$; hence
$P(z)=z^sp_\varepsilon(z)$, with $d=s+m$, proves (i).

\emph{(i)$\Rightarrow$(iii)$\Rightarrow$(iv).} \Cref{main:prop:delayedcertificate} proves (i)$\Rightarrow$(iii) with
$k_0=d+1$, and block iteration together with
\eqref{main:eq:basiccontraction} and~\eqref{main:eq:errorequivalence} proves
(iii)$\Rightarrow$(iv).

\emph{(iv)$\Rightarrow$(ii), and a start that reproduces fixed-point iteration.} Suppose that weights solving \eqref{main:eq:momentcondition} exist, so that (ii) fails. We construct a start from which AA($m$) reproduces fixed-point iteration, so that (iv) fails as well. Since weights exist, $K_q\ge1$, so
$\rho(W)=q$ and the summand $qU$ is present. Its space $M$ is the
orthogonal sum of real invariant subspaces on which $U$ acts as $\pm1$ or
as a plane rotation $R_{\theta_i}$. Choose a unit vector $u_i$ in one such
subspace for each phase $\theta_i$, and set $r_0=\sum_i\sqrt{w_i}u_i$. Then
$\norm{r_0}=1$ and $\ip{U^\ell r_0}{r_0}=\sum_iw_i\cos(\ell\theta_i)=q^\ell$
for $1\le\ell\le m$. For the candidate sequence
$r_k=W^kr_0=q^kU^kr_0$ and each retained lag $1\le\ell\le m_k$, $\ip{r_k}{r_{k-\ell}-r_k} =q^{2k-\ell}\ip{U^\ell r_0}{r_0}-q^{2k}=0$.
Thus $r_k$ is orthogonal to every direction of its retained affine hull,
hence is the minimum-norm point of that hull, and
\eqref{main:eq:residualupdate} gives $r_{k+1}=Wr_k$. Induction from the first
step $r_1=Wr_0$ shows that AA($m$) started at $x_0=x_*+(W-I)^{-1}r_0$
generates this sequence, with $\norm{r_k}=q^k$. This proves
(iv)$\Rightarrow$(ii), and the theorem's assertion for the case that the conditions fail. The same $r_0$ also explains, without reference to AA, why no certificate exists in this case: $\ip{W^jr_0}{W^kr_0}=q^{j+k}\ip{U^{|j-k|}r_0}{r_0}=q^{2\max(j,k)}$ whenever $|j-k|\le m$, so every real $P(z)=\sum_{j=s}^da_jz^j$ with $P(1)=1$ and $d-s\le m$ satisfies $\norm{P(W)r_0}^2=\sum_{j,k}a_ja_kq^{2\max(j,k)}\ge q^{2d}$ by supplementary \eqref{si:eq:generalgrambaseline}.

\emph{Part (v).} Assume that $W_0$ is absent or $\norm{W_0}<q$. Since $\norm W=q$, the summand $qU$ is then present, so $M\ne\{0\}$ and $\rho(W)=q$. If (ii) holds, the polynomial
$p_\varepsilon$ above has degree at most $m$, and its definition gives
$\norm{p_\varepsilon(W_0)}\le\norm{W_0}^m+\varepsilon\sum_{\ell=1}^m|h_\ell|q^\ell\bigl(\norm{W_0}^{m-\ell}+\norm{W_0}^m\bigr)$,
in which the sum does not depend on $\varepsilon$; since $\norm{W_0}<q$, decreasing $\varepsilon$ makes the right side smaller than $q^m$;
hence $\norm{p_\varepsilon(W)}=\max\{\norm{p_\varepsilon(qU)},\norm{p_\varepsilon(W_0)}\}<q^m$, which is (v). Conversely, (v) gives
(iii) with $k_0=m+1$ by \cref{main:thm:certificate} with $d=m$. A normal $W$
satisfies the assumption, because then $W_0$ is normal and
$\norm{W_0}=\rho(W_0)<q$.
\end{proof}

For a fixed $q$, the phases of the peripheral eigenvalues decide
whether AA($m$) beats $q$. The remaining block $W_0$ affects the gain
and the possible delay, but not the answer. The theorem is qualitative:
its block length and factor may depend on $W$. Supplementary
\cref{si:prop:shortblock,si:prop:explicitdelay,si:prop:nonnormaloidexplicit}
give explicit choices, while \cref{main:sec:certificate} gives
dimension-independent rates for classes sharing a certificate.

For normal $W$, no delay is needed and
$\norm{p(W)}=\max_{\lambda\in\spec(W)}|p(\lambda)|$. Thus AA($m$)
beats $q$ from every start exactly when a real polynomial of degree
at most $m$, normalized by $p(1)=1$, is smaller than $q^m$ in modulus
on the spectrum. A bound on a containing region $q\Omega$ gives a
certificate for every normal matrix with spectrum there. More
generally, no delay is needed whenever $\norm{W_0}<q$. Normaloidity
alone does not ensure this: $\diag(q,qJ)$ from the preceding section
has $\norm{W_0}=q$ and no degree-one certificate, yet AA(1) terminates
within five steps from every start (supplementary
\cref{si:rem:normhypothesis}).

For a symmetric matrix the eigenvalues of modulus $q$ lie in $\{q,-q\}$, so the phases belong to $\{0,\pi\}$ and $\cos(\ell\theta_i)=\pm1$. Write $w_+$ and $w_-$ for the weights at $q$ and $-q$. The second moment equation would give $w_++w_-=q^2$, which contradicts $w_++w_-=1$, so every memory $m\ge2$ strictly improves the factor $q$. At memory one the only moment equation is $w_+-w_-=q$. If only one of $q$ and $-q$ is an eigenvalue, its weight is one, and the equation reads $1=q$ or $-1=q$, which is impossible; if both are, it has the solution $w_+=(1+q)/2$, $w_-=(1-q)/2$. Hence AA(1) fails to beat $q$ exactly when both $q$ and $-q$ are eigenvalues. This recovers the memory conditions of~\cite[Theorem~3.1 and Proposition~3.3]{garner2026}, and both ingredients of the present theory are already visible there. The weights $(1+q)/2$ and $(1-q)/2$ are the squared coordinates of the initial residual in~\cite[proof of Proposition~3.3]{garner2026} from which AA(1) reproduces fixed-point iteration. The degree-one certificate $p(z)=(z-c)/(1-c)$, with
$c=(\lambda_{\max}+\lambda_{\min})/2$ the midpoint of the spectrum, is the
polynomial form of the ball estimate~\cite[Lemma~5.1]{garner2026}: its norm
$\norm{p(W)}=(\lambda_{\max}-\lambda_{\min})/(2-\lambda_{\max}-\lambda_{\min})$
is less than $q$ exactly when $q$ and $-q$ are not both eigenvalues. When $-q$ is not an
eigenvalue, \cite[Theorem~3.1]{garner2026} also gives an explicit factor for
every memory through a variational formula; supplementary
\cref{si:cor:lineartwoeigen} gives the simpler bound
$\rhAA(x_0)\le q-[1-(5/8)^{1/4}](1-q)(q-q')$ when $\spec(W)\subseteq[-q',q]$
with $q'<q$, including memory one.

\begin{remark}[Complex matrices]\label{main:rem:complex}
\Cref{main:thm:certificate,main:prop:delayedcertificate} hold verbatim for complex matrices, residuals, and coefficients. \Cref{main:thm:normalalternative} holds for a complex contraction with complex polynomials and with \eqref{main:eq:momentcondition} replaced by $\sum_iw_ie^{i\ell\theta_i}=q^\ell$ for $1\le\ell\le m$, over all distinct eigenvalues $qe^{i\theta_i}$ of modulus $q$; these are $2m$ real conditions, for a real matrix the two systems are solvable together, and complex AA from a real start coincides with real AA.
\end{remark}

\subsection{How much memory is needed}\label{main:sec:memorycount}
The moment equations \eqref{main:eq:momentcondition} for memory $m$ contain those for every smaller memory, so by \cref{main:thm:normalalternative} the smallest memory that beats $q$ is the first $m$ for which they have no solution. Their right-hand sides $q^\ell$ are the cosine moments of the Poisson kernel $P_q(\theta)=(1-q^2)/(1-2q\cos\theta+q^2)$, since $\frac1{2\pi}\int_0^{2\pi}\cos(\ell\theta)P_q(\theta)\,d\theta=q^\ell$. Thus \eqref{main:eq:momentcondition} asks whether a probability distribution supported on the finitely many phases $\theta_i$ can have the same first $m$ cosine moments as the Poisson distribution $P_q(\theta)\,d\theta/(2\pi)$ on the circle. Two facts answer this, in terms of the number $N$ of distinct eigenvalues of modulus $q$. Some matrices with $N$ such eigenvalues match the moments through degree $N-1$: the moments are those of the unitary dilation of the scalar contraction $q$, and a finite-dimensional dilation through order $N-1$ realizes them (\cref{main:lem:momentmatrix}). No matrix with $N$ such eigenvalues matches them through degree $N$, because the Poisson distribution has a positive density while the phases are finitely many points (\cref{main:prop:peripheralcount}). Hence memory $N$ always suffices, and memory $N-1$ may not.

\begin{lemma}[Realizing Poisson moments]\label{main:lem:momentmatrix}
Let $m\ge1$ and $0<q<1$. On the standard basis $\mathbf e_0,\ldots,\mathbf e_m$ of $\R^{m+1}$ define
\begin{equation}\label{main:eq:momentmatrix}
 U\mathbf e_0=q\mathbf e_0+\sqrt{1-q^2}\,\mathbf e_1,\qquad
 U\mathbf e_j=\mathbf e_{j+1}\quad(1\le j<m),\qquad
 U\mathbf e_m=\sqrt{1-q^2}\,\mathbf e_0-q\mathbf e_1.
\end{equation}
Then $U$ is orthogonal, $v=\mathbf e_0$ satisfies
\begin{equation}\label{main:eq:moments}
 U^jv=q^j \mathbf e_0+\sqrt{1-q^2}\sum_{\ell=1}^jq^{j-\ell}\mathbf e_\ell,
 \qquad \ip v{U^jv}=q^j\qquad(0\le j\le m),
\end{equation}
the vectors $v,Uv,\ldots,U^mv$ are linearly independent, and $U$ has $m+1$ distinct eigenvalues. Consequently, the normal matrix $qU$ has $m+1$ distinct eigenvalues of modulus $q$, its moment equations \eqref{main:eq:momentcondition} have a solution through degree $m$, and AA($m$) does not beat $q$ on it.
\end{lemma}
\begin{proof}
The images in \eqref{main:eq:momentmatrix} are orthonormal, so $U$ is
orthogonal, and \eqref{main:eq:moments} follows by induction on $j$. The
coordinate matrix of $v,Uv,\ldots,U^mv$ is triangular with nonzero diagonal entries $1,\sqrt{1-q^2},\ldots,\sqrt{1-q^2}$. Thus $v$ is a cyclic vector, the minimal polynomial of $U$
has degree $m+1$, and the diagonalizable matrix $U$ has $m+1$ distinct eigenvalues. The squared norms $w_i$ of the components of $v$ in the real invariant subspaces of $U$, on which $U$ acts as $\pm1$ or as a rotation $R_{\theta_i}$, satisfy $\sum_iw_i\cos(\ell\theta_i)=\ip v{U^\ell v}=q^\ell$ for $0\le\ell\le m$, so they solve the moment equations of $qU$, and \cref{main:thm:normalalternative} gives the last assertion.
\end{proof}

The matrix $U$ is a finite-dimensional unitary dilation of the scalar contraction $q$ through order $m$~\cite[Chapter~I]{sznagy2010}: it keeps the fraction $q$ of $\mathbf e_0$, moves the remainder $\sqrt{1-q^2}$ into the fresh direction $\mathbf e_1$, shifts it along $\mathbf e_1,\ldots,\mathbf e_m$, and returns it to $\mathbf e_0$ only after $m+1$ steps. This is why $\ip v{U^jv}=q^j$ holds exactly for $j\le m$, the lags that the moment equations of memory $m$ involve, while $\ip v{U^{m+1}v}=q^{m+1}+1-q^2$.

Both bounds of the next proposition come from polynomials in $\cos\theta$ that vanish at every phase $\theta_i$: a nonnegative one of degree $N_q$, whose Poisson integral is always positive, and the node polynomial of degree $K_q$, whose Poisson integral vanishes only for exceptional phases.

\begin{proposition}[Counting the eigenvalues on the circle]
\label{main:prop:peripheralcount}
Let $W$ be real with $q=\norm W\in(0,1)$, let $N_q$ be the number of its distinct eigenvalues of modulus $q$, and let $K_q$ be the number of these with nonnegative imaginary part, written $qe^{i\theta_1},\ldots,qe^{i\theta_{K_q}}$ with $0\le\theta_i\le\pi$ as in \cref{main:thm:normalalternative}.
Then the equivalent conditions of \cref{main:thm:normalalternative} hold
for every $m\ge N_q$. They also hold for every $m\ge K_q$, unless the node polynomial $\omega(x)=\prod_{i=1}^{K_q}(x-\cos\theta_i)$ has Poisson mean zero:
\begin{equation}\label{main:eq:peripheralexception}
 \int_0^{2\pi}\omega(\cos\theta)P_q(\theta)\,d\theta=0.
\end{equation}
The bound $m\ge N_q$ cannot be improved in general: for every $N\ge2$, the matrix $qU$ of \cref{main:lem:momentmatrix} with $m=N-1$ has $N_q=N$ distinct eigenvalues of modulus $q$, and AA($N-1$) does not beat $q$ from every start.
\end{proposition}
\begin{proof}
By \cref{main:thm:normalalternative}, some trajectory fails to accelerate
only if there are weights $w_i\ge0$ with $\sum_iw_i=1$ and
$\sum_iw_i\cos(\ell\theta_i)=q^\ell$ for $1\le\ell\le m$. Since $\frac1{2\pi}\int_0^{2\pi}\cos(\ell\theta)P_q(\theta)\,d\theta=q^\ell$ for every $\ell\ge0$, and $\cos(\ell\theta)$ is a polynomial of degree $\ell$ in $\cos\theta$, such weights would give $\sum_iw_iQ(\cos\theta_i)=\frac1{2\pi}\int_0^{2\pi}Q(\cos\theta)P_q(\theta)\,d\theta$ for
every polynomial $Q$ of degree at most $m$. If $Q$ vanishes at every $\cos\theta_i$, the left side is zero, so its Poisson integral must vanish too. For the first claim, let $Q$ be the product of
$(x-\cos\theta_i)^2$ over the phases other than $0$ and $\pi$, multiplied by
$1-x$ if $q$ is an eigenvalue and by $1+x$ if $-q$ is. It has degree $N_q$, vanishes at every $\cos\theta_i$, and is nonnegative and not identically zero on $[-1,1]$, so its Poisson integral is positive, because $P_q>0$, and no weights exist when $m\ge N_q$. For the second claim, take $Q=\omega$, of degree $K_q$: weights can exist for $m\ge K_q$ only if \eqref{main:eq:peripheralexception} holds. The last claim is \cref{main:lem:momentmatrix} with $m=N-1$.
\end{proof}

For two or fewer eigenvalues of modulus $q$ the proposition gives the complete answer. A single eigenvalue, $q$ or $-q$, needs memory one. A conjugate pair $qe^{\pm i\theta}$ has $K_q=1$ and $\omega(x)=x-\cos\theta$, whose Poisson mean is $q-\cos\theta$; memory one suffices unless $\cos\theta=q$, that is, unless the pair has real part $q^2$, and then memory two is needed. The pair $\{q,-q\}$ has $K_q=N_q=2$ and needs memory two; this is the symmetric case, and the matrix $qU$ of \cref{main:lem:momentmatrix} with $N=2$.

With three or more eigenvalues of modulus $q$, at least one conjugate pair is present, so $K_q<N_q$ and the two bounds of \cref{main:prop:peripheralcount} differ. Memory $K_q$ suffices outside the exceptional set \eqref{main:eq:peripheralexception}, which has measure zero, and only configurations in that set need more. Both bounds are attained. For the equally spaced phases $2\pi j/\nu$, memory $K_q=\lfloor\nu/2\rfloor+1$ is needed (\cref{main:thm:cusps}). The matrix $qU$ of \cref{main:lem:momentmatrix} needs memory $N_q=N$, so it lies in the exceptional set; for $N=3$ its eigenvalues are $q$ and a pair with $\cos\theta=-(1-q)/2$, and $K_q=2$. In the variable $x=\cos\theta$ the moment equations ask for a quadrature rule: nonnegative weights at the nodes $\cos\theta_i$ that integrate every polynomial of degree at most $m$ exactly against the Poisson distribution, since $\cos(\ell\theta)$ is a polynomial of degree $\ell$ in $\cos\theta$. A configuration needs memory $N_q$ exactly when such a rule is exact through degree $N_q-1$, the largest degree that the polynomial $Q$ in the proof allows. Such nodes are those of the Gauss, Gauss--Radau, or Gauss--Lobatto rules for the Poisson measure, according to whether none, one, or both of $\pm q$ are eigenvalues (supplementary \cref{si:sec:peripheralconfigurations,si:tab:peripheral}). The examples above are of this kind: $\{q,-q\}$ gives the Lobatto nodes $\pm1$, and $q$ together with the pair $\cos\theta=-(1-q)/2$ gives a Radau rule.

\subsection{Equally spaced peripheral eigenvalues}\label{main:sec:cusps}\label{main:sec:deltoid}

Equally spaced phases are the configuration that periodicity forces. By the Perron--Frobenius theorem, the eigenvalues of maximal modulus of an irreducible nonnegative matrix of period $\nu$ are the $\nu$ points $\rho(W)e^{2\pi ik/\nu}$. When $\norm W=\rho(W)$, as for the periodic Markov chains of \cref{main:ex:markov} in the stationary inner product, these are the eigenvalues of modulus $q$. For these phases the count $K_q$ of \cref{main:prop:peripheralcount} is exactly the memory needed (\cref{main:thm:cusps}).

A useful family of spectral regions with this periphery is bounded by hypocycloids. For an integer $\nu\ge2$ let
\begin{equation}\label{main:eq:hypocycloid}
 \gamma_\nu(\theta)=\frac{\nu-1}{\nu}e^{i\theta}+\frac1\nu e^{-i(\nu-1)\theta},
 \qquad 0\le\theta\le2\pi,
\end{equation}
and let $\Gamma_\nu$ be the closed region it bounds, the hypocycloid with
$\nu$ cusps; $\Gamma_2=[-1,1]$, $\Gamma_3$ is the deltoid, and $\Gamma_4$ is the astroid. For $q>0$ we write $q\Gamma_\nu=\{qz:z\in\Gamma_\nu\}$ for the scaled region. The triangle inequality gives $|\gamma_\nu(\theta)|\le1$,
with equality exactly when $e^{i\nu\theta}=1$, so that
\begin{equation}\label{main:eq:hypocycloidperiphery}
 \Gamma_\nu\subseteq\{z:|z|\le1\},\qquad
 \Gamma_\nu\cap\{z:|z|=1\}=\{e^{2\pi ik/\nu}:0\le k<\nu\}.
\end{equation}
Hence every $W$ with $\norm W=q$ and spectrum in $q\Gamma_\nu$ has its eigenvalues of modulus $q$ among the points $qe^{2\pi ik/\nu}$. These regions arise from the mean-zero random walk with steps $1$ and $1-\nu$ in the work of Cowal, Marshall, and Pollock, whose polynomials, defined by the recurrence~\cite[Section~1.4.2, eq.~(16)]{cowal2026random}, are bounded on $\Gamma_\nu$.

\begin{theorem}[Memory for equally spaced peripheral eigenvalues]\label{main:thm:cusps}
Let $W$ be real with $q=\norm W\in(0,1)$, and suppose that every eigenvalue
of $W$ of modulus $q$ has the form $qe^{2\pi ik/\nu}$ for a fixed integer
$\nu\ge2$. If $m\ge\lfloor\nu/2\rfloor+1$, the equivalent conditions of
\cref{main:thm:normalalternative} hold; in particular $\rhAA(x_0)<q$ for every starting point $x_0$, uniformly over blocks of a fixed length, which is $m+1$ if $W$ is normal. The threshold is sharp:
for every $m\le\lfloor\nu/2\rfloor$ there are a real normal
$W\in\R^{\nu\times\nu}$ with $\spec(W)=\{qe^{2\pi ik/\nu}:0\le k<\nu\}$ and
an initialization from which AA($m$) follows fixed-point iteration exactly.
Thus $\lfloor\nu/2\rfloor+1$ is the smallest memory that guarantees strict acceleration for all such $W$, and also for all real $W$ with $\norm W=q$ and $\spec(W)\subseteq q\Gamma_\nu$.
\end{theorem}
The proof, in supplementary \cref{si:sec:cuspsproof}, rests on aliasing: at the phases $2\pi k/\nu$ one has $\cos(\ell\theta)=\cos((\nu-\ell)\theta)$, so for $m\ge\lfloor\nu/2\rfloor+1$ two of the conditions $\sum_iw_i\cos(\ell\theta_i)=q^\ell$, $0\le\ell\le m$, have the same left side but different right sides, while for $m\le\lfloor\nu/2\rfloor$ the discrete Fourier transform of $q^{\min\{j,\nu-j\}}$ gives nonnegative weights that solve them.

Periodic Markov chains give natural nonnormal examples.

\begin{corollary}[Periodic Markov chains]\label{main:ex:markov}
Let $P$ be an irreducible row-stochastic matrix with stationary distribution $\pi$ and period $\nu\ge2$, let $q\in(0,1)$, and let AA($m$) for the iteration $x\mapsto qPx+b$ pose its least-squares problems in the inner product $\ip xy_\pi=\sum_i\pi_ix_iy_i$, with induced norm $\norm\cdot_\pi$. Then $\norm{qP}_\pi=q=\rho(qP)$, and AA($m$) beats the factor $q$ of the iteration from every start if $m\ge\lfloor\nu/2\rfloor+1$, whereas for every $m\le\lfloor\nu/2\rfloor$ some start makes it reproduce the iteration exactly. For an aperiodic chain, $\nu=1$, memory one suffices.
\end{corollary}
\begin{proof}
Since every row of $P$ is a probability vector and $\pi^TP=\pi^T$, Jensen's inequality gives
$\norm{Px}_\pi^2=\sum_i\pi_i\bigl(\sum_jP_{ij}x_j\bigr)^2
\le\sum_i\pi_i\sum_jP_{ij}x_j^2=\norm x_\pi^2$, with equality for constant
$x$; hence $\norm{qP}_\pi=q=\rho(qP)$. By the Perron--Frobenius theorem the
eigenvalues of $qP$ of modulus $q$ are the $\nu$ points $qe^{2\pi ik/\nu}$, each simple. All results of this paper transfer to the inner product of $\pi$ by supplementary \cref{si:prop:conjugation}. For $m\ge\lfloor\nu/2\rfloor+1$, \cref{main:thm:cusps} gives the acceleration. For $m\le\lfloor\nu/2\rfloor$, the weights constructed in the proof of \cref{main:thm:cusps} solve the moment equations \eqref{main:eq:momentcondition} at the phases $2\pi k/\nu$, all of which occur, so \cref{main:thm:normalalternative} provides the start that reproduces the iteration. For $\nu=1$ the only eigenvalue of modulus $q$ is $q$ itself, and memory one suffices by \cref{main:thm:normalalternative}.
\end{proof}
The iteration $x\mapsto qPx+b$ is the policy-evaluation iteration for a fixed policy in a discounted Markov decision process with discount factor $q$~\cite{sutton2018}. The matrix $qP$ is self-adjoint in the inner product of $\pi$ when the chain is reversible, but in general it is not normal; supplementary \cref{si:sec:normalnormaloid} gives a four-state example. Posing the least-squares problem in this inner product needs $\pi$ only up to a positive factor, since rescaling an inner product does not change the minimizing coefficients. If $P$ is doubly stochastic, $\pi$ is uniform and the Euclidean inner product qualifies; otherwise $\norm P_2>1$, and, provided $q\norm P_2<1$, the Euclidean criterion compares AA with the larger factor $q\norm P_2$.

For $\nu=2$ the eigenvalues of modulus $q$ lie in $\{q,-q\}$, and \cref{main:thm:cusps} extends the symmetric case discussed after \cref{main:thm:normalalternative} to every real $W$ whose eigenvalues of modulus $q$ are real. For the astroid, $\nu=4$, memory three
is needed, and the obstruction weights at memory two are
$(1+q)^2/4,(1-q^2)/4,(1-q)^2/4,(1-q^2)/4$ at the phases
$0,\pi/2,\pi,3\pi/2$, written as a probability distribution on all four phases; in the form \eqref{main:eq:momentcondition}, with phases in $[0,\pi]$, the pair $\pm\pi/2$ carries the combined weight $(1-q^2)/2$. The memory needed grows like half the number
of cusps, and as $\nu\to\infty$ the regions $\Gamma_\nu$ fill the disk, on
which no memory accelerates uniformly. Only the points where the spectrum
touches the circle enter the threshold; the rest of the region, cusped or
smooth, affects the size of the margin (\cref{main:sec:sqrtstatus}).

The case $\nu=3$, the deltoid $\Gamma_3$, is our main example: it is the simplest of these regions whose contact points are not all real, and the square-root question of \cref{main:sec:limits} is posed on it. Its contact points with the unit circle are $1$, $\omega_3$, and $\omega_3^2$, where $\omega_3=e^{2\pi i/3}$, and it meets the real axis in $[-1/3,1]$; it is the spectral region studied by Cowal, Marshall, and Pollock~\cite{cowal2026deltoid}. By \cref{main:thm:cusps} with $\nu=3$, memory two is the smallest memory that guarantees acceleration, $\rhAA(x_0)<q$ for every starting point $x_0$, whenever the eigenvalues of $W$ of modulus $q$ lie among $q$, $q\omega_3$, and $q\omega_3^2$, in particular whenever $\norm W=q$ and $\spec(W)\subseteq q\Gamma_3$.

\section{Beyond normaloid matrices: the numerical radius}\label{main:sec:numrad}

For normaloid matrices, beating $q=\norm W$ means improving on
the worst-case fixed-point factor $\rho(W)$. When $\rho(W)<q$,
fixed-point iteration itself beats $q$, so the exact criterion of
\cref{main:thm:normalalternative} no longer settles that comparison.
We now keep the inner product fixed and bound AA's factor by the
numerical radius. This bound is sharper than $q$ and identifies a class
on which AA is never slower than worst-case fixed-point iteration;
an example shows that this comparison can fail outside that class.

The eigenvalues of modulus $\rho(W)$ then no longer decide the question. Every step contracts the residual by at most $\norm W$, so $\rhAA(x_0)\le\norm W$, but for a matrix that is not normaloid this bound can be improved. The numerical radius
\[
 w(W)=\max\{|z^*Wz|:\ z\in\C^n,\ \norm z=1\}
\]
satisfies $\rho(W)\le w(W)\le\norm W$. It equals $\norm W$ exactly for normaloid matrices, because a unit vector $z$ with $|z^*Wz|=\norm W$ must be an eigenvector with an eigenvalue of modulus $\norm W$. Matrices with $w(W)=\rho(W)$ are called \emph{spectraloid}~\cite{halmos1982}; every normal matrix is spectraloid. The next theorem shows that $w(W)$ is a ceiling for AA: an upper bound on $\rhAA(x_0)$ for every starting point and every memory, which is attained for some matrices (supplementary \cref{si:thm:numradequality}).

\begin{theorem}[The numerical radius is a ceiling]\label{main:thm:numradceiling}
Let $W\in\R^{n\times n}$ and $m\ge0$, where AA($0$) is fixed-point iteration. For every starting point, every $\delta>0$, and every $k\ge2$, the residuals of AA($m$) satisfy
\begin{equation}\label{main:eq:numradbound}
 \norm{r_k}\le\norm W\,\delta^{-1/2}\bigl((1+\delta)\,w(W)\bigr)^{k-1}\norm{r_0}.
\end{equation}
Whenever $I-W$ is invertible, the residual bound implies $\rhAA(x_0)\le w(W)$; in particular, this applies to every affine contraction. The same bounds hold for capped AA($m$) with any cap $\Lambda\ge1$.
\end{theorem}

\begin{proof}
If $w(W)=0$, then $W=0$ and the claim is immediate. Otherwise write $y_{k+1}$ for the minimum-norm point of the affine hull of the retained residuals, so that $r_{k+1}=Wy_{k+1}$ and $y_1=r_0$. For $k\ge1$ the residual $r_k=Wy_k$ belongs to the window, and $y_{k+1}$ is the orthogonal projection of the origin onto the affine hull, so $\ip{y_{k+1}}{r_k-y_{k+1}}=0$, that is,
\begin{equation}\label{main:eq:growthidentity}
 \norm{y_{k+1}}^2=\ip{y_{k+1}}{Wy_k}\qquad(k\ge1).
\end{equation}
Fix $s>w(W)$ and $L\ge2$, and consider the vector-valued trigonometric polynomial $\zeta(\theta)=\sum_{k=1}^Le^{ik\theta}s^{-k}y_k$. Since the $y_k$ are real, $\zeta(\theta)^*W\zeta(\theta)=\sum_{j,k}e^{i(k-j)\theta}s^{-j-k}\ip{y_j}{Wy_k}$, and averaging against $e^{i\theta}$ retains exactly the pairs $j=k+1$:
\[
 \frac1{2\pi}\int_0^{2\pi}e^{i\theta}\,\zeta(\theta)^*W\zeta(\theta)\,d\theta=\sum_{k=1}^{L-1}s^{-2k-1}\ip{y_{k+1}}{Wy_k}=s\sum_{k=2}^{L}s^{-2k}\norm{y_k}^2
\]
by \eqref{main:eq:growthidentity}. On the other hand $|\zeta^*W\zeta|\le w(W)\norm\zeta^2$ pointwise, and $\frac1{2\pi}\int_0^{2\pi}\norm{\zeta(\theta)}^2d\theta=\sum_{k=1}^Ls^{-2k}\norm{y_k}^2$ by Parseval's identity. Hence
\[
 \bigl(s-w(W)\bigr)\sum_{k=2}^Ls^{-2k}\norm{y_k}^2\le w(W)\,s^{-2}\norm{y_1}^2 .
\]
Every term on the left is bounded by the right side; with $s=(1+\delta)w(W)$ this gives $\norm{y_k}\le\delta^{-1/2}((1+\delta)w(W))^{k-1}\norm{y_1}$ for $k\ge2$, and \eqref{main:eq:numradbound} follows from $\norm{r_k}\le\norm W\norm{y_k}$ and $y_1=r_0$. For capped AA the mixed point minimizes the norm over a convex subset of the affine hull that still contains $r_k$, so $\ip{y_{k+1}}{r_k-y_{k+1}}\ge0$; the argument uses \eqref{main:eq:growthidentity} only in the form $\norm{y_{k+1}}^2\le\ip{y_{k+1}}{Wy_k}$.
\end{proof}

\begin{corollary}\label{main:cor:spectraloid}
(a) If $W$ is spectraloid, then $\rhAA(x_0)\le\rho(W)$ for every starting point: AA($m$) is never asymptotically slower than the worst case of fixed-point iteration. If $W$ is normal, the stronger comparison $\rhAA(x_0)\le\rho_{\mathrm{Picard}}(x_0)$ holds from the same start. (b) If $W$ is not normaloid, then $\rhAA(x_0)\le w(W)<\norm W$ for every starting point, with the constants in \eqref{main:eq:numradbound} independent of the start. (c) If $\rho(W)>0$, AA($m$) is at most the factor $w(W)/\rho(W)$ slower than the worst case of fixed-point iteration. When $\rho(W)<w(W)$ no run that reproduces fixed-point iteration, on $\R^n$ or on an invariant subspace, has factor $w(W)$, because such a run has factor at most $\rho(W)$.
\end{corollary}

\begin{proof}
Parts (a) and (b) follow from \cref{main:thm:numradceiling}, since $w(W)=\rho(W)$ for spectraloid $W$ and $w(W)<\norm W$ for $W$ that is not normaloid. For normal $W$ and $r_0\ne0$, restrict the run to $M_0=\operatorname{span}\{W^jr_0:j\ge0\}$. This restriction is normal, its spectral radius is $\rho_{\mathrm{Picard}}(x_0)$, and its AA run is unchanged (supplementary \cref{si:prop:invariantsubspace}); the case $r_0=0$ is immediate. For (c), \cref{main:thm:numradceiling} gives $\rhAA(x_0)\le w(W)=\bigl(w(W)/\rho(W)\bigr)\rho(W)$, and a run that reproduces fixed-point iteration, on $\R^n$ or on an invariant subspace, has $\norm{r_k}=\norm{W^kr_0}$ and hence factor at most $\rho(W)$.
\end{proof}

Part (c) is not vacuous: AA can be slower than fixed-point iteration even when the dominant eigenvalue is real and simple, a configuration in which memory one beats $\rho(W)$ whenever $W$ is normaloid (\cref{main:prop:peripheralcount}).

\begin{example}[AA(1) slower than fixed-point iteration]
\label{main:prop:notbenchmark}
Let
\begin{equation*}
 E=\frac13\begin{pmatrix}1&-1\\2&2\end{pmatrix},\qquad
 W_\lambda=\diag(\lambda,E)\in\R^{3\times3},\qquad 2/3<\lambda<(2/9)^{1/4}.
\end{equation*}
Then $\norm{W_\lambda}=2\sqrt2/3<1$ and $\rho(W_\lambda)=\lambda$ is a simple eigenvalue, yet AA(1) started with the residual $r_0=\mathbf e_2$ satisfies $r_{k+4}=-(2/9)\,r_k$ for all $k\ge0$, so $\rhAA(x_0)=(2/9)^{1/4}>\rho(W_\lambda)$, whereas fixed-point iteration from the same start converges at the rate $\rho(E)=2/3<\rho(W_\lambda)$. The verification is in supplementary \cref{si:sec:notbenchmark}.
\end{example}

The slowness comes from the block $E$, whose norm exceeds $\rho(W_\lambda)$. The run respects the ceiling, since $(2/9)^{1/4}<w(W_\lambda)=2/\sqrt7$. Supplementary \cref{si:sec:notbenchmark} verifies the orbit exactly.

\paragraph{Below the ceiling.}
A related moment condition guarantees a uniform factor strictly below the ceiling. Its equations are posed on the \emph{touching set} of $W$,
\begin{equation}\label{main:eq:touchingset}
 \Theta(W):=\{\theta\in(-\pi,\pi]:\ w(W)e^{-i\theta}\in\mathcal W(W)\},\qquad
 \mathcal W(W):=\{z^*Wz:\ z\in\C^n,\ \norm z=1\},
\end{equation}
the set of angles at which the numerical range $\mathcal W(W)$ touches the circle of radius $w(W)$. With angles viewed modulo $2\pi$, this set is compact and symmetric under $\theta\mapsto-\theta$, because $W$ is real. For normaloid $W$ it consists of the phases of the eigenvalues of modulus $\norm W$, and the moment equations below are those of \cref{main:thm:normalalternative}. For a probability measure $\nu$ on $\Theta(W)$ and $w=w(W)$ they read
\begin{equation}\label{main:eq:measuremoments}
 \int\cos(j\theta)\,d\nu(\theta)=w^j,\qquad j=1,\dots,m.
\end{equation}
The next theorem is the analogue, for the ceiling, of the uniform half of \cref{main:thm:normalalternative}; $y_k$ is the minimum-norm point of the window, as in \cref{main:sec:geometry}.

\begin{theorem}[Uniform improvement below the numerical radius]\label{main:thm:uniformceiling}
Let $m\ge1$ and $0<w=w(W)<1$, and suppose that no probability measure $\nu$ on $\Theta(W)$ satisfies \eqref{main:eq:measuremoments}.
This holds whenever $\Theta(W)$ is finite and $m\ge|\Theta(W)|$, in particular when $m\ge2$ and $\Theta(W)$ has at most two points; it also holds when $m=1$ and $\Theta(W)$ is a single real point or a pair $\{\theta_0,-\theta_0\}$ with $\cos\theta_0\ne w$.
Then there is a constant $C\ge1$, depending only on $W$ and $m$, such that every ordinary AA($m$) run satisfies
\begin{equation}\label{main:eq:uniformtail}
 \sum_{j\ge k}w^{-2j}\norm{y_j}^2\le C\,w^{-2k}\norm{y_k}^2\qquad(k\ge1),
\end{equation}
and therefore
\begin{equation}\label{main:eq:uniformrate}
 \norm{r_k}\le\norm W\sqrt C\,\bigl(w\sqrt{1-1/C}\bigr)^{k-1}\norm{r_0}\qquad(k\ge1).
\end{equation}
In particular $\rhAA(x_0)\le w\sqrt{1-1/C}<w$ for every starting point, and for every $k_0\ge2+C\ln(4C)$ there is $c<1$ with $\norm{r_{k+k_0}}\le c\,w^{k_0}\norm{r_k}$ for all $k\ge0$ and all runs.
\end{theorem}
The proof, in supplementary \cref{si:sec:uniformceilingproof},
applies the energy argument of \cref{main:thm:numradceiling} to every
tail of a run. Separation of the target moments from those supported
on $\Theta(W)$ forces a gain controlled by the energy at the start of
the tail, giving a common rate. The proof does not give an explicit
$C$. Under quadratic contact, supplementary
\cref{si:prop:touchingcertificate} gives explicit constants using a
polynomial vanishing on the touching set, such as $z-1$ for $\{0\}$ or
$z^2-2z\cos\theta_0+1$ for $\{\pm\theta_0\}$. Supplementary
\cref{si:sec:attaining} examines equality at the ceiling and the case
left open when the moment equations are feasible.

\section{The size of the acceleration}
\label{main:sec:limits}
For the spectral regions considered here, fixed-degree certificates provide gains of order $1-q$. This section first sharpens the certificate bounds of \cref{main:sec:certificate} to gains linear in the certificate margin (\cref{main:sec:lineargain}), the order that matters as $q\uparrow1$, and then asks whether that order reflects the behavior of windowed AA or only the strength of the bounds. The first answer is a worst-case limitation: at every fixed memory $m\ge2$ there are matrices with a fixed quadratic certificate on which AA converges at the rate $1-2(1-q)$ (\cref{main:sec:everymemoryobstruction}). The second depends on the spectral geometry: smooth contact imposes the same order on finite-step guarantees for every polynomial method, whereas cusps permit gains of order $\sqrt{1-q}$. Classical stationary methods, which combine a fixed number of past iterates with fixed coefficients, attain such gains (\cref{main:sec:sqrtstatus}), so a fixed window is not itself the obstruction; whether windowed AA attains them is open.

\subsection{Gains linear in the margin}\label{main:sec:lineargain}
\Cref{main:thm:certificate} guarantees an explicit gain per block, but, as the following corollary shows, only a gain quadratic in the margin $1-\eta$.

\begin{corollary}[Explicit certificate margin]\label{main:cor:explicitconstant}
Under the assumptions of \cref{main:thm:certificate}, with $\eta$ and $\beta$ as in \eqref{main:eq:constants},

\begin{equation}\label{main:eq:explicitrate}
 1-c\ge\frac{(1-\eta)^2}{2(1+\beta^2)},\qquad
 \rhAA(x_0)\le q\left[1-\frac{(1-\eta)^2}{2(d+1)(1+\beta^2)}\right]<q.
\end{equation}
\end{corollary}
\begin{proof}
Rationalizing the expression for $1-c$ gives $1-c=(1-\eta)^2/(1+\beta^2-\eta+\beta\sqrt{1+\beta^2-\eta^2}) \ge(1-\eta)^2/(2(1+\beta^2))$,
because $\beta\sqrt{1+\beta^2-\eta^2}\le\beta^2+1/2$.
Use $c^{1/(d+1)}\le1-(1-c)/(d+1)$ in \eqref{main:eq:rootbound}.
\end{proof}

The gain $1-c$ depends on two features of the certificate: its margin $1-\eta$ and the size of its coefficients, measured by $\beta$. By \eqref{main:eq:explicitrate} it is at least $(1-\eta)^2/(2(1+\beta^2))$. It is also less than $1/(2\beta^2)$, however small $p(W)$ is: $c$ increases with $\eta$, so every certificate has $c\ge\beta/\sqrt{1+\beta^2}$. A useful certificate therefore needs both a small $p(W)$ and moderate coefficients. For example, $p(z)=(z-q)/(1-q)$ annihilates $W=qI$, and AA(1) terminates after two steps, since the line through $r_0$ and $r_1=qr_0$ contains the origin; but $\beta=1/(1-q)$, so the theorem guarantees a gain of less than $(1-q)^2/2$ per block. For certificates with large coefficients the bound can thus be far from the actual behavior; for two-term certificates such as this one, supplementary \cref{si:thm:twoterm} gives a bound that degrades much more slowly with the coefficients, as discussed at the end of this subsection.

The gain of \cref{main:cor:explicitconstant} is quadratic in the margin
because the proof compares vectors: it bounds the distance between
$\sum_ja_jr_{k+j}$ and $p(W)r_k$ by the square root of the projection energy
that the block releases. Comparing squared norms instead, that is, comparing
$\norm{p(W)r_k}^2$ with the value it takes when the block gains nothing,
gives a gain linear in the margin. This matters when margins are small: on the deltoid, the best margins of certificates of any fixed degree $d\ge2$ are of order $1-q$ as $q\uparrow1$ (supplementary \cref{si:prop:fixeddegreecertificate}), so the quadratic bound would give gains of order $(1-q)^2$ where the linear one gives order $1-q$. The improvement concerns the order of the gain in the margin, not the size of the factor: for large margins \cref{main:thm:certificate} can give the smaller factor, and neither bound dominates the other (supplementary \cref{si:rem:compare}).

For a real polynomial $p(z)=\sum_{i=0}^da_iz^i$ with $p(1)=1$, define its \emph{excess energy}
\begin{equation}\label{main:eq:generalexcess}
 \Delta_p=\sum_{i,j=0}^da_ia_jq^{2\max(i,j)}-\norm{p(W)}^2.
\end{equation}
The double sum has a simple meaning. Suppose that AA gains nothing over $q$ on the block of $d+1$ steps that starts at $r_k$, that is, $y_{k+j+1}=r_{k+j}$ and $\norm{r_{k+j}}=q^j\norm{r_k}$ for $0\le j\le d$. Then $r_{k+j}=W^jr_k$ by \eqref{main:eq:residualupdate}, and \eqref{main:eq:orthogonality} gives $\ip{W^ir_k}{W^jr_k}=q^{2\max(i,j)}\norm{r_k}^2$, so the double sum equals $\norm{p(W)r_k}^2/\norm{r_k}^2$. It is at least $q^{2d}$, with equality only for $p(z)=z^d$ (supplementary \eqref{si:eq:generalgrambaseline}). Hence every certificate has $\Delta_p\ge q^{2d}-\norm{p(W)}^2>0$, and in general $\Delta_p>0$ says that $p(W)$ is too small to be compatible with a block without gain.

\begin{theorem}[Linear certificate margin]\label{main:thm:generalmargin}
Let $p(z)=\sum_{i=0}^da_iz^i$ be an undelayed certificate for AA($m$) of degree $d\ge1$. Then
there is a number $\tau_p\in(0,1/(2(d+1))]$, given explicitly in supplementary \eqref{si:eq:taudef} and depending only on $q$, $\norm{p(W)}$, and the coefficients of $p$, such that every AA($m$) trajectory satisfies
\begin{align}
 \norm{r_{k+d+1}}^2
   &\le(1-\tau_p)q^{2(d+1)}\norm{r_k}^2
       \qquad(k\ge0),\label{main:eq:generalblocklinear}\\
 \rhAA(x_0)&\le q(1-\tau_p)^{1/[2(d+1)]}<q.\label{main:eq:generalrootlinear}
\end{align}
Moreover, with $\eta$ as in \eqref{main:eq:constants}, 
\begin{equation}\label{main:eq:generalnormalizedmargin}
 \tau_p\ge\frac{q^{2d}(1-\eta)}{4(d+1)^4\bigl(\sum_{i=0}^d|a_i|\bigr)^2}.
\end{equation}
\end{theorem}
The gain is therefore linear in the margin $1-\eta$, with a constant that depends only on $q$, the degree, and the size of the coefficients. The proof, in supplementary \cref{si:sec:generalmarginproof}, bounds $\Delta_p$ by the energy $L$ that a block releases, through an inequality $\Delta_p\le K_pL+G_p\sqrt L$ whose constants $K_p$ and $G_p$ depend only on $q$ and the coefficients of $p$ (supplementary \eqref{si:eq:generalK} and~\eqref{si:eq:generalG}): a block that releases little energy behaves almost like a block without gain, which $\Delta_p>0$ rules out. The argument uses the first residual $r_k$ in every least-squares problem of the block, which is why it needs an undelayed certificate. The norm test $\norm{p(W)}<q^d$ enters only at the end, to turn this inequality into the simple constant $\tau_p$. Without it, the argument needs only $\Delta_p>0$, which is a weaker requirement: the norm test asks that $\norm{p(W)}^2$ be smaller than $q^{2d}$, whereas $\Delta_p>0$ asks only that it be smaller than the double sum in \eqref{main:eq:generalexcess}, which is at least $q^{2d}$. For instance, $p(z)=2-z$ fails the norm test for $W=qI$, since $\norm{p(W)}=2-q>q$, but $\Delta_p=4q(1-q)>0$.

\begin{corollary}[Excess-energy certificate]\label{main:cor:excesscertificate}
Let $m\ge d\ge1$, and let $p(z)=\sum_{i=0}^da_iz^i$ be real with $p(1)=1$ and $\Delta_p>0$. Then there is a number $\widehat\tau_p\in(0,1)$, given explicitly in supplementary \eqref{si:eq:generalsharptau} and depending only on $q$, $\Delta_p$, and the coefficients of $p$, such that \eqref{main:eq:generalblocklinear} and~\eqref{main:eq:generalrootlinear} hold with $\tau_p$ replaced by $\widehat\tau_p$. If $\norm{p(W)}<q^d$, then $\widehat\tau_p\ge\tau_p$.
\end{corollary}

The number $\widehat\tau_p$ is the positive root of
$K_p\tau+G_p\sqrt\tau=\Delta_p$. It increases with $\Delta_p$ and satisfies
$\widehat\tau_p\ge\min\{\Delta_p/(2K_p),\,\Delta_p^2/(4G_p^2)\}$,
with the second term interpreted as $+\infty$ when $G_p=0$:
at this value each term on the left is at most $\Delta_p/2$.
For a quadratic certificate on the deltoid it gives
$\widehat\tau_p\approx0.0197\,(1-q)$ as $q\uparrow1$
(supplementary \cref{si:cor:quadraticlinearmargin}). It also covers
$\Delta_p>0$ when the norm test fails, although the gain may then be
quadratic in $\Delta_p$. \Cref{main:sec:nonlinear} allows either
hypothesis for nonlinear maps.

All three constants $c$, $\tau_p$, and $\widehat\tau_p$ depend on $W$ only through $q$ and $\norm{p(W)}$, and each bound weakens as $\norm{p(W)}$ grows. One may therefore replace $\norm{p(W)}$ by any upper bound $B$, provided that $B<q^d$ for $c$ and $\tau_p$, and that $\Delta_p$ stays positive with $B$ in place of $\norm{p(W)}$ for $\widehat\tau_p$. This is how the bounds are used in practice: $\norm{p(W)}$ is rarely known exactly, but for a normal matrix with spectrum in a region $\Omega$ it is at most $\max_\Omega|p|$, which gives one rate for all such matrices.

The coefficient dependence can be improved for two-term
certificates $p(z)=a_0+a_dz^d$. Supplementary \cref{si:thm:twoterm}
proves \eqref{main:eq:generalblocklinear} with a constant at least
$(1-\eta)/(1+da_d)$, which decreases only linearly in $a_d$ and $d$.
For $(z-q)/(1-q)$ on $W=qI$, this constant is $(1-q)/(2-q)$,
giving a gain of order $1-q$ instead of the quadratic-order gain
from \cref{main:thm:certificate}. Supplementary
\cref{si:cor:cubiclinearmargin} applies the same refinement to the
best two-term cubic certificate on the deltoid,
$(27z^3-13q^3)/(27-13q^3)$.

\subsection{Counterexamples at every fixed memory}
\label{main:sec:everymemoryobstruction}

The examples of this subsection show that the gain guaranteed by a fixed certificate can be of linear order only. Each matrix has two blocks. An auxiliary block raises the norm to $q$ without raising $\norm{W^2}$, which makes $z^2$ a certificate; the residuals of the run have no component in this block, so it does not affect the iteration. On the other, active block, AA($m$) produces exactly the iterates of fixed-point iteration, which converges there with factor $\rho(W)=1-2(1-q)<q$. Thus AA beats $q$, as the certificate guarantees, but no more than fixed-point iteration does: only by $1-q$, not by an amount of order $\sqrt{1-q}$.

\begin{theorem}[A fixed quadratic certificate without square-root acceleration]
\label{main:thm:everymemoryquadratic}
Fix $m\ge2$. For every $0<\varepsilon<1/4$, there are an
invertible real nonsymmetric matrix
$W_\varepsilon\in\R^{(m+3)\times(m+3)}$, diagonalizable over
$\C$, and a nonzero initial point $x_{0,\varepsilon}$ for the iteration
$g(x)=W_\varepsilon x$, such that
\begin{equation}
 \norm{W_\varepsilon}=1-\varepsilon, 
 \qquad
 \norm{W_\varepsilon^2}=(1-2\varepsilon)^2,
 \qquad
 \rhAA(x_{0,\varepsilon})
 =\rho_{\mathrm{Picard}}(x_{0,\varepsilon})
 =\rho(W_\varepsilon)=1-2\varepsilon.
 \label{main:eq:everymemoryfactor}    
\end{equation}
Here AA is ordinary AA($m$), every minimizing coefficient vector is unique,
and it selects the newest residual. Thus the fixed polynomial $p(z)=z^2$ has the margin $1-\norm{W_\varepsilon^2}/\norm{W_\varepsilon}^2=(2\varepsilon-3\varepsilon^2)/(1-\varepsilon)^2$, but no constant $c_m>0$ gives
$\rhAA(x_0)\le1-c_m\sqrt{1-\norm W}$ uniformly over this class.
\end{theorem}

The proof, in supplementary \cref{si:sec:everymemoryfurther}, assembles the two blocks described above. That section also shows that the obstruction is robust: it applies equally to capped AA($m$) with any cap $\Lambda\ge1$ (\cref{main:sec:geometry}), it occurs near a fixed point of an analytic map, and it is not an artifact of the Euclidean inner product. Two further facts are proved there: the factor cannot be bounded away from $q$ uniformly over matrices with a strict certificate, and normality does not restore square-root acceleration from individual starts.

\subsection{Smooth contact, cusps, and the square-root question}
\label{main:sec:sqrtstatus}
The matrices of \cref{main:thm:everymemoryquadratic} are not normal, and their spectra are not contained in the deltoid $q\Gamma_3$ (supplementary \cref{si:sec:contactproof}). For normal matrices, the attainable gain depends on how the spectral region meets the circle $|z|=q$, through a contrast established by Cowal, Marshall, and Pollock~\cite[Theorems~1.2 and~1.4]{cowal2026random}: polynomials bounded on a region with a cusp at $1$ can grow like $(1+c\sqrt\varepsilon)^n$ at $1+\varepsilon$, whereas polynomials bounded on a region containing an ellipse tangent at $1$ grow at most like $(1+c\varepsilon)^n$. The next theorem states what this contrast means for convergence. We call a method \emph{polynomial} if its residuals have the form $r_n=p_n(W)r_0$ with $\deg p_n\le n$ and $p_n(1)=1$, where $p_n$ may depend on $W$ and $r_0$; AA($m$) for every $m$, and GMRES, are polynomial methods. For $0<\delta<1$ let $E_\delta$ be the closed region bounded by the ellipse with vertices $\pm1$ and co-vertices $\pm i\delta$, and let $\Gamma_\nu$ be the hypocycloid region of \eqref{main:eq:hypocycloid}.

\begin{theorem}[Smooth contact and cusps]\label{main:thm:contact}
Let $1/2<q<1$.
\begin{enumerate}[label=(\alph*)]
\item Let $0<\delta<1$. For every $n\ge1$ there are a real normal matrix $W$ of dimension at most $2n+3$ with $\norm W=q$ and $\spec(W)\subseteq qE_\delta$ and a unit vector $r_0$ on which every polynomial method satisfies
\[
 \norm{r_n}\ \ge\ \Bigl(1+\frac{3(1-q)}{2\delta q}\Bigr)^{-n}.
\]
\item Let $\nu\ge2$, and let $R_\nu>1$ solve $(\nu-1)R_\nu+R_\nu^{-(\nu-1)}=\nu/q$, so that
\[
 R_\nu^{-1}=1-\sqrt{\frac{2(1-q)}{\nu-1}}+O(1-q).
\]
For every $n\ge1$ there are a real normal matrix $W$ of dimension at most $2n+3$ with $\norm W=q$ and $\spec(W)\subseteq q\Gamma_\nu$ and a unit vector $r_0$ on which every polynomial method satisfies $\norm{r_n}\ge R_\nu^{-n}$. The stationary $\nu$-step method
\[
 x_{k+1}=\mu_0g(x_k)+(1-\mu_0)x_{k+1-\nu},\qquad \mu_0=\frac{\nu}{(\nu-1)qR_\nu},
\]
initialized by $x_j=x_0$ for $1-\nu\le j\le0$, has worst-case root-linear factor $R_\nu^{-1}$ on the real normal matrices with spectrum in $q\Gamma_\nu$.
\end{enumerate}
\end{theorem}

The proof is in supplementary \cref{si:sec:contactproof}. With smooth contact, part (a) shows that no polynomial method, AA included, has a dimension-independent guarantee better than a per-step factor $1-O(1-q)$, so the linear-order bounds of the certificates have the correct order. Because part (a) holds for every polynomial method, including GMRES, which terminates on any fixed matrix, its worst-case matrix must depend on $n$: it bounds the first $n$ steps, not the asymptotic factor on a given matrix. For AA, \cref{main:thm:everymemoryquadratic} supplies the asymptotic statement: a single matrix of fixed size on which the asymptotic factor of AA($m$) is only linearly below $q$. At a cusp, part (b) shows that gains of order $\sqrt{1-q}$ are possible and that a fixed window is not itself the obstruction. The question is whether AA attains them: for each fixed $m\ge2$, are there constants $c_m>0$ and $q_m<1$ such that AA($m$) satisfies
\begin{equation}\label{main:eq:sqrttarget}
 \rhAA(x_0)\le1-c_m\sqrt{1-q}
 \qquad(q_m<q<1),
\end{equation}
for every dimension, every normal matrix with spectrum in $q\Gamma_3$, and every starting point? The class contains the symmetric matrices with spectrum in $[-q/3,q]$, so the question includes symmetric matrices.

The available results do not decide \eqref{main:eq:sqrttarget}. Memory two gives a rate $1-C(1-q)$ on this class, with explicit constants in supplementary \cref{si:sec:deltoidrates}, and certificates of fixed degree cannot prove more, since their best margin on $q\Gamma_3$ is of order $1-q$ (supplementary \cref{si:prop:fixeddegreecertificate}). On the other hand, AA($m$) with $m\ge2$ improves strictly on the fixed-point factor of every start whose factor is close enough to one, by a margin of linear order in one minus that factor (supplementary \cref{si:thm:deltoidactivelinear}). A square-root gain, if it exists, must come from the dynamics of AA itself, whose residuals, like those of the stationary methods of \cref{main:thm:contact}, are polynomials of growing degree in $W$.

\section{Local convergence for nonlinear maps}
\label{main:sec:nonlinear}
Near a fixed point $x_*$ of a nonlinear map $g$, the derivative
$J=g'(x_*)$ plays the role of $W$, and the certificate and its margin
refer to $J$. The residuals satisfy the linear recurrence of the previous
sections only up to a remainder, and the remainder is small only if the
mixing coefficients are not too large. 
Earlier analyses of nonlinear AA bound the coefficients~\cite{toth2015,evans2020}. Garner, Lerman, and Zhang~\cite[(3.5) and Theorem~3.6]{garner2026} build the bound into the algorithm: their modified AA($m$) constrains each coefficient by a constant $C_0$ in the least-squares problem, and they prove local convergence when $C_0$ exceeds a threshold, the Jacobian is symmetric and not a multiple of the identity, and the remainders are quadratic. Capped AA($m$) below is their modification with the bound placed on the $\ell^1$ norm of the coefficients; the cap may be any number greater than one, and only differentiability at $x_*$ is assumed.

Capped AA($m$), defined in \cref{main:sec:geometry}, starts from a single point $x_0$ and, at every step
with $r(x_k)\ne0$, chooses any minimizer of
\begin{equation}\label{main:eq:nonlinearcap}
 \min_{\substack{\sum_{j=k-m_k}^k\alpha_j=1\\
                         \sum_{j=k-m_k}^k|\alpha_j|\le\Lambda}}
       \norm*{\sum_{j=k-m_k}^k\alpha_jr(x_j)},
 \qquad
 x_{k+1}=\sum_{j=k-m_k}^k\alpha_jg(x_j);
\end{equation}
it stops when the current residual vanishes. The feasible set is compact,
convex, and nonempty, so the minimizing residual combination $y_{k+1}$ is
unique, although the coefficients, and hence $x_{k+1}$, need not be; all
statements below cover every minimizing selection. Projection onto a convex
set turns the identities \eqref{main:eq:orthogonality} and~\eqref{main:eq:pythagoras}
into inequalities, and the feasible points $(1+\delta)r_k-\delta r_i$, with
$\delta=(\Lambda-1)/2$ and $r_i$ retained, recover the moments needed by the
energy argument. Supplementary \cref{si:sec:cappedbounds} proves two
estimates on this basis. First, if $p$ is a real polynomial of degree
$d\le m$ with $p(1)=1$ and $\norm{p(J)}<q^d$, capped affine AA($m$) satisfies
\begin{equation}\label{main:eq:cappedfactor}
 \norm{r_{k+d+1}}^2\le(1-\tau_{p,\Lambda})\,q^{2(d+1)}\norm{r_k}^2,
\end{equation}
with an explicit $\tau_{p,\Lambda}\in(0,1)$ that depends only on $q$, $\norm{p(J)}$,
$\Lambda$, and the coefficients of $p$, and equals $\tau_p$ of
\cref{main:thm:generalmargin} when $\Lambda\ge7/5$
(\cref{si:thm:cappedgeneral}); the same holds under the excess-energy
condition $\Delta_p>0$ with the analogue of $\widehat\tau_p$
(\cref{si:cor:cappedexcesscertificate}). Second, this block estimate is
robust under small perturbations of the residual recurrence
(\cref{si:lem:compactblock}).

\begin{theorem}[Local nonlinear convergence with any fixed cap]
\label{main:thm:nonlinear}
Let $g$ be defined in a neighborhood of $x_*$, with $g(x_*)=x_*$, and
assume that $g$ is Fr\'echet differentiable at $x_*$. Let
$J=g'(x_*)$ and $q=\norm J\in(0,1)$. Suppose a real polynomial $p$ of
degree $d\ge1$ with $p(1)=1$ satisfies $\norm{p(J)}<q^d$ or $\Delta_p>0$, fix
$m\ge d$ and $\Lambda>1$, and let $\tau_{p,\Lambda}\in(0,1)$ be the corresponding constant in \eqref{main:eq:cappedfactor}. Every sufficiently close
initialization generates a sequence under \eqref{main:eq:nonlinearcap} that
remains near $x_*$ and converges to it, and every choice of minimizing
coefficients satisfies
\begin{equation}\label{main:eq:nonlinearroot}
 \limsup_{k\to\infty}\norm{x_k-x_*}^{1/k}
 \le q\,(1-\tau_{p,\Lambda})^{1/(2(d+1))}<q,
\end{equation}
where a terminating run is extended by $x_*$. For $\Lambda\ge7/5$, $\tau_{p,\Lambda}$ is the constant $\tau_p$ of \cref{main:thm:generalmargin}, or $\widehat\tau_p$ of \cref{main:cor:excesscertificate}, for $J$.
\end{theorem}
The proof, given in supplementary \cref{si:sec:nonlinearproof}, follows the
standard pattern of local analyses of nonlinear AA (see, e.g., \cite[Section~5.5]{garner2026}). Differentiability gives
$g(x_*+w)=x_*+Jw+o(\norm w)$. The cap bounds the new error by a multiple of the
envelope $T_k$ of the retained residuals, so that $r_{k+1}=Jy_{k+1}+o(T_k)$.
An induction keeps the iterates in a ball around $x_*$ and contracts $T_k$ over
every $m+1$ steps, which proves convergence. Once these remainders are small
relative to $T_k$, the perturbed block estimate of supplementary
\cref{si:lem:compactblock} gives $\norm{r_{k+d+1}}\le(1-\tau_{p,\Lambda})^{1/2}q^{d+1}\norm{r_k}+o(T_k)$,
and a delayed comparison gives \eqref{main:eq:nonlinearroot}.

Supplementary \cref{si:prop:nonlinearinexact} allows inexact coefficient solves, and \cref{si:rem:uncapped} explains why the cap cannot simply be dropped.

\section{Concluding remarks}\label{main:sec:conclusion}
The main theorem gives a necessary and sufficient condition for AA($m$) to beat $q=\norm W$ from every start: the phases of the eigenvalues of modulus $q$ must support no probability distribution with the prescribed cosine moments. A probability distribution with these moments gives an explicit start from which AA follows fixed-point iteration; if none exists, a polynomial certificate gives uniform contraction over a fixed number of steps. The window constrains the width of the certificate's support, and aliasing explains the sharp memory threshold for periodic Markov chains.

With the inner product fixed, every run has factor at most $w(W)$, so AA is never slower than worst-case fixed-point iteration when $w(W)=\rho(W)$. The touching-set moment test can improve this bound uniformly, although a factor below $w(W)$ need not be below $\rho(W)$. When $q=\rho(W)$, the main criterion gives the exact comparison. If $\rho(W)>0$ and the dominant eigenvalues are semisimple, an adapted inner product can achieve this equality, but constructing it generally requires spectral information; stationary distributions supply a natural choice for Markov chains.

Two questions remain. Can ordinary AA with fixed memory attain the dimension-independent square-root gain permitted by cusped spectral regions? And, for a nonspectraloid matrix, is the factor uniformly below $w(W)$ even when the touching set carries the prescribed moments?

\section*{Acknowledgments}
GL was supported in part by the U.S. National Science Foundation under award DMS-2427955 and TZ was supported in part by the U.S. National Science Foundation under award DMS-2318926. GL thanks Nicholas Marshall for a very clear talk on his work with Peter Cowal and Sara Pollock at the IMA Data Science Seminar. GL is grateful to Anthropic for providing a free team account, which was also shared with TZ. 

\section*{Statement on the use of AI tools}
The authors conceived the extension of the analysis of Garner, Lerman, and
Zhang~\cite{garner2026} to polynomial certificates and designed the overall
approach that led to \cref{main:thm:certificate,main:thm:normalalternative}; they stated and partially
proved a first version of these theorems for a normal matrix $W$, and recognized that the
numerical radius is needed for general $W$. They then used AI tools, mainly
Claude Opus~5.5 and Claude Fable~5.1 (Anthropic), to formulate mathematical
statements, develop and check proofs, run numerical tests (most of them not
reported here, though they guided the choice of next steps), answer questions
raised by the authors, and draft and revise text. The authors gave
feedback throughout, but the use of AI tools was heavy. ChatGPT (OpenAI) was also used for a subsequent mathematical and editorial review. The authors are
responsible for the results, the verification of the mathematical arguments,
and the final text.

\begingroup
\raggedright
\bibliographystyle{plain}
\bibliography{references}
\endgroup

\clearpage
\pdfbookmark[0]{Supplementary material}{supplementary-material}

\section*{Supplementary material}
\setcounter{section}{0}
\setcounter{subsection}{0}
\setcounter{equation}{0}
\setcounter{theorem}{0}
\setcounter{table}{0}
\setcounter{figure}{0}
\renewcommand{\thesection}{S\arabic{section}}
\renewcommand{\thetable}{S\arabic{table}}
\renewcommand{\thefigure}{S\arabic{figure}}
\renewcommand{\theHsection}{S\arabic{section}}
\renewcommand{\theHsubsection}{\theHsection.\arabic{subsection}}
\renewcommand{\theHsubsubsection}{\theHsubsection.\arabic{subsubsection}}
\renewcommand{\theHequation}{\theHsection.\arabic{equation}}
\renewcommand{\theHtheorem}{\theHsection.\arabic{theorem}}
\renewcommand{\theHtable}{S\arabic{table}}
\renewcommand{\theHfigure}{S\arabic{figure}}

The supplement uses the algorithm and notation of the main article. Ordinary AA($m$) starts from one point and retains at most $m+1$ residuals; a zero residual terminates the iteration. Vector norms are Euclidean unless a fixed weighted inner product is specified. For affine maps $q=\norm W$; for nonlinear maps $q=\norm{g'(x_*)}$.

\paragraph{Organization.} \Cref{si:tab:roadmap} lists the contents of the supplement and the parts of the article that each section supports. Its sections follow the order of the article and contain the proofs of the results that are not proved there, together with the explicit constants, refinements, and examples to which the article refers.

\begin{center}
\footnotesize\setlength{\tabcolsep}{4pt}
\renewcommand{\arraystretch}{1.05}
\begin{tabular}{@{}p{0.09\linewidth}p{0.56\linewidth}p{0.31\linewidth}@{}}
\toprule
\raggedright Section & \raggedright Contents & \raggedright Supports in the article\tabularnewline
\midrule
\raggedright \Cref{si:sec:certificateextensions} & \raggedright Proofs of \cref{main:thm:generalmargin,main:cor:excesscertificate}; the bound $\norm{W_0^t}<q^t$ on the complementary block of \eqref{main:eq:peripheralsplit}, with explicit block lengths and delays for \cref{main:thm:normalalternative}; what simplifies for normal matrices and what the normaloid case adds; an exact test for degree-one certificates; the proof of \cref{main:thm:cusps}; a catalogue of peripheral configurations and the smallest memory they need & \raggedright \Cref{main:sec:certificate,main:sec:lineargain}; \cref{main:thm:normalalternative,main:prop:peripheralcount,main:thm:cusps} and \cref{main:ex:markov} in \cref{main:sec:normal}\tabularnewline[2pt]
\raggedright \Cref{si:sec:deltoidrates} & \raggedright Explicit quadratic certificates on the deltoid, the memory-two rates they give, and the exact quadratic minimax value & \raggedright The deltoid rates quoted in \cref{main:sec:cusps,main:sec:sqrtstatus}\tabularnewline[2pt]
\raggedright \Cref{si:sec:metricdynamics} & \raggedright Change of inner product (\cref{si:prop:conjugation}); the peripheral splitting \eqref{main:eq:peripheralsplit} (\cref{si:lem:normattained}); adapted inner products and the spectral radius as the benchmark; invariant subspaces; verification of \cref{main:prop:notbenchmark}; runs at the ceiling; the proof of \cref{main:thm:uniformceiling} and explicit constants & \raggedright \Cref{main:sec:geometry}; \cref{main:thm:normalalternative}; \cref{main:sec:numrad}\tabularnewline[2pt]
\raggedright \Cref{si:sec:sispecialized} & \raggedright Comparison of \cref{main:thm:certificate,main:thm:generalmargin}; refinements of the energy constants; the best two-term cubic certificate on the deltoid; the linear-order limit of fixed-degree certificates; two-term certificates at every depth, with the symmetric and the cubic deltoid rates & \raggedright Constants for \cref{main:sec:certificate,main:sec:lineargain}; the bounds quoted after \cref{main:thm:normalalternative} and in \cref{main:sec:sqrtstatus}\tabularnewline[2pt]
\raggedright \Cref{si:sec:sinormal} & \raggedright A quantitative improvement over the active fixed-point factor on the deltoid; the proof of \cref{main:thm:everymemoryquadratic} and the robustness of its examples; a normal counterexample at every memory; polynomial approximation on the scaled hypocycloids; the proof of \cref{main:thm:contact} & \raggedright \Cref{main:sec:limits}\tabularnewline[2pt]
\raggedright \Cref{si:sec:cappedbounds} & \raggedright The energy bound \eqref{main:eq:cappedfactor} for capped AA with any cap $\Lambda>1$, with its excess-energy form, and the robust block estimate under small errors in the residual recurrence & \raggedright The two estimates behind \cref{main:thm:nonlinear}\tabularnewline[2pt]
\raggedright \Cref{si:sec:sinonlinear} & \raggedright An explicit perturbation modulus; failure without control of the minimizing coefficients; the proof of \cref{main:thm:nonlinear} and its version with inexact coefficient solves & \raggedright \Cref{main:sec:nonlinear}\tabularnewline
\bottomrule
\end{tabular}
\par\vspace{\abovecaptionskip}
\refstepcounter{table}\label{si:tab:roadmap}
{\normalsize Table~\thetable: Contents of the supplement and the parts of the article that each section supports.}
\end{center}

\section{Certificate bounds and peripheral configurations}
\label{si:sec:certificateextensions}

\subsection{Proofs of \texorpdfstring{\cref{main:thm:generalmargin}}{the linear-margin theorem} and \texorpdfstring{\cref{main:cor:excesscertificate}}{the excess-energy corollary}}\label{si:sec:generalmargin}\label{si:sec:generalmarginproof}
The two statements are in \cref{main:sec:lineargain}, together with the excess energy $\Delta_p$ of \eqref{main:eq:generalexcess}. The proof bounds $\Delta_p$ by the energy $L$ that a block of $d+1$ steps releases, through the inequality $\Delta_p\le K_pL+G_p\sqrt L$ of \eqref{si:eq:generalenergyinequality}, whose constants $K_p$ and $G_p$ are defined in \eqref{si:eq:generalK} and~\eqref{si:eq:generalG}; the numbers $\tau_p$ and $\widehat\tau_p$ quoted in the statements are defined in \eqref{si:eq:taudef} and~\eqref{si:eq:generalsharptau}. Refinements of these constants are collected in \cref{si:sec:sispecialized}.

\begin{proof}[Proof of \cref{main:thm:generalmargin}]
As in the proof of \cref{main:thm:certificate}, assume $\norm{r_k}=1$.
Rescaled by powers of $q$, the vectors $u_j=q^{-j}r_{k+j}$ and
$v_j=q^{-j}y_{k+j+1}$ satisfy $u_{j+1}=Tv_j$ with the contraction $T=W/q$,
and the corrections $f_j=u_j-v_j$ are orthogonal to $v_j$. Let
$D=(I-T^*T)^{1/2}$ be the defect operator of $T$, so that
$\norm x^2-\norm{Tx}^2=\norm{Dx}^2$ for every $x$. Each step
loses $\norm{u_j}^2-\norm{u_{j+1}}^2=\norm{f_j}^2+\norm{Dv_j}^2$.
Consequently the energy lost over the block, $L=1-\norm{u_{d+1}}^2$,
satisfies
\begin{equation}\label{si:eq:generalenergy}
 \sum_{j=0}^d\bigl(\norm{f_j}^2+\norm{Dv_j}^2\bigr)=L.
\end{equation}
All $\norm{u_j},\norm{v_j}$ are at most one, and
$\norm{v_j}^2\ge1-L$.

\emph{Raw powers.}
Put $t_j=T^ju_0$ and $\psi_j=t_j-v_j$.
The recurrence gives
$t_\ell=v_\ell+\sum_{i=0}^\ell T^{\ell-i}f_i$.
For each $i<j$, contraction-defect telescoping gives
\(\norm{T^{j-i}f_i}^2+\sum_{\ell=i}^{j-1}\norm{DT^{\ell-i}f_i}^2=\norm{f_i}^2\).
Apply the triangle inequality in the product space consisting of the
$j$ defect coordinates and the terminal coordinate $\psi_j$, taking
$(Dv_0,\ldots,Dv_{j-1},f_j)$ as the base vector:
\[
 \left(\sum_{\ell<j}\norm{Dt_\ell}^2+\norm{\psi_j}^2\right)^{1/2}
 \le\left(\sum_{\ell<j}\norm{Dv_\ell}^2+\norm{f_j}^2\right)^{1/2}
       +\sum_{i<j}\norm{f_i}
 \le\sqrt{(j+1)(1-\norm{v_j}^2)}.
\]
The last step uses Cauchy--Schwarz and the energy identity through
the projection at index $j$. Since $\sum_{\ell<j}\norm{Dt_\ell}^2=1-\norm{t_j}^2$,
\begin{equation}\label{si:eq:generalrawenergy}
 1-\norm{t_j}^2+\norm{\psi_j}^2
 \le(j+1)(1-\norm{v_j}^2)\le(j+1)L.
\end{equation}

\emph{Moments and the Gram matrix.}
The projection producing $v_j$ retains $r_k=u_0$ for $j\le d$.
Thus $\ip{v_j}{u_0}=q^j\norm{v_j}^2$, so $h_j:=u_0-q^jv_j$ is orthogonal to $v_j$ and $\norm{h_j}^2=1-q^{2j}\norm{v_j}^2\le1-q^{2j}+q^{2j}L$.
For $\mu_j=\ip{u_0}{t_j}$, write $\mu_j=q^j\norm{v_j}^2+q^j\ip{v_j}{\psi_j}+\ip{h_j}{\psi_j}$.
The first two error terms combine by polarization:
\(q^j(\norm{v_j}^2-1)+q^j\ip{v_j}{\psi_j}
 =q^j(\norm{t_j}^2+\norm{v_j}^2-\norm{\psi_j}^2-2)/2\).
This quantity is nonpositive and its magnitude is at most
$q^j(j+2)L/2$. Bounding $\ip{h_j}{\psi_j}$ by Cauchy--Schwarz and using
$\sqrt{x+y}\le\sqrt x+\sqrt y$ therefore gives
\begin{equation}\label{si:eq:generalmomenterror}
 |\mu_j-q^j|
 \le M_jL+\sqrt{(j+1)(1-q^{2j})L},\qquad
 M_j=q^j\Bigl(\frac{j+2}{2}+\sqrt{j+1}\Bigr).
\end{equation}
For $0\le i\le j\le d$, telescoping also gives $\ip{t_i}{t_j} =\mu_{j-i}-\sum_{\ell=0}^{i-1} \ip{Dt_\ell}{Dt_{\ell+j-i}}$,
and both index sets in the last sum lie in $\{0,\ldots,d-1\}$.
By Cauchy--Schwarz and \eqref{si:eq:generalrawenergy},
\begin{equation}\label{si:eq:generalgramerror}
 \bigl|\ip{t_i}{t_j}-\mu_{|j-i|}\bigr|\le(d+1)L.
\end{equation}
Put $\hat a_i=a_iq^i$ and $C_h=2\sum_{i=0}^{d-h}\hat a_i\hat a_{i+h}$ for $1\le h\le d$. Since $p(W)u_0=\sum_i \hat a_it_i$ and $\mu_0=1$, 
\(\norm{p(W)u_0}^2
 \ge\sum_i \hat a_i^2+\sum_{h=1}^dC_h\mu_h-(d+1)(\sum_i|\hat a_i|)^2L\).
The ideal Gram quadratic form is
\begin{equation}\label{si:eq:generalgrambaseline}
 \sum_i \hat a_i^2+\sum_{h=1}^dC_hq^h
 =\sum_{i,j=0}^d a_ia_jq^{2\max(i,j)}
 =q^{2d}+E_p,\qquad
 E_p=(1-q^2)\sum_{\ell=0}^{d-1}q^{2\ell}(a_0+\cdots+a_\ell)^2,
\end{equation}
where the last equality follows from 
$q^{2\max(i,j)}=q^{2d}+(1-q^2)
\sum_{\ell=\max(i,j)}^{d-1}q^{2\ell}$
and $\sum_i a_i=1$. Thus $E_p\ge0$, with equality only for $p(z)=z^d$.

\emph{Absorbing the moment error.}
The identity
\(\hat a_i=q^i\sum_{\ell=0}^i a_\ell-q(q^{i-1}\sum_{\ell=0}^{i-1}a_\ell)\)
for \(0\le i<d\),
with the empty sum equal to zero, gives
$(\sum_{i=0}^{d-1}\hat a_i^2)^{1/2}\le(1+q)\sqrt{E_p/(1-q^2)}$.
Thus Cauchy--Schwarz in each lag sum gives
\(|C_h|\le2(1+q)(\sum_{i=0}^d\hat a_i^2)^{1/2}\sqrt{E_p/(1-q^2)}\),
and since $1-q^{2h}=(1-q^2)\sum_{\ell<h}q^{2\ell}$, the numbers
\begin{equation}\label{si:eq:generalG}
 G_p=\sum_{h=1}^d|C_h|\sqrt{(h+1)(1-q^{2h})},\qquad
 H_p=2(1+q)\Bigl(\sum_{i=0}^d\hat a_i^2\Bigr)^{1/2}
       \sum_{h=1}^d\Bigl((h+1)\sum_{\ell=0}^{h-1}q^{2\ell}\Bigr)^{1/2}
\end{equation}
satisfy
\begin{equation}\label{si:eq:generalabsorb}
 G_p\le H_p\sqrt{E_p}.
\end{equation}
Combining \eqref{si:eq:generalmomenterror}, \eqref{si:eq:generalgramerror}, \eqref{si:eq:generalgrambaseline}, and~\eqref{si:eq:generalabsorb}
with $\norm{p(W)u_0}\le\norm{p(W)}$ proves the central inequality
\begin{equation}\label{si:eq:generalenergyinequality}
 \Delta_p\le K_pL+G_p\sqrt L\le K_pL+H_p\sqrt{E_pL},
\end{equation}
where
\begin{equation}\label{si:eq:generalK}
 K_p=(d+1)\Bigl(\sum_{i=0}^d|\hat a_i|\Bigr)^2+\sum_{h=1}^d|C_h|M_h.
\end{equation}
Since $\norm{p(W)}<q^d$, one has $E_p\le\Delta_p$, and Young's inequality gives $\Delta_p\le K_pL+\Delta_p/2+H_p^2L/2$,
hence $L\ge\tau_p$, where
\begin{equation}\label{si:eq:taudef}
 \tau_p=\frac{\Delta_p}{2K_p+H_p^2}.
\end{equation}
 Moreover $\sum_i|\hat a_i|>0$, $K_p\ge(d+1)(\sum_i|\hat a_i|)^2$,
and $\Delta_p\le q^{2d}+E_p\le(\sum_i|\hat a_i|)^2$
by \eqref{si:eq:generalgrambaseline}.
Therefore $0<\tau_p\le1/[2(d+1)]<1$.
Restoring the residual scale proves \eqref{main:eq:generalblocklinear},
including the zero-residual case, and block iteration gives the root
bound. Finally $\Delta_p\ge q^{2d}-\norm{p(W)}^2\ge q^{2d}(1-\eta)$. With $\norm a_1=\sum_i|a_i|$ one has $|\hat a_i|\le|a_i|$, $|C_h|\le2\norm a_1^2$, $M_h\le h+2$ because $\sqrt{h+1}\le(h+2)/2$, $(\sum_i\hat a_i^2)^{1/2}\le \norm a_1$, and $\sum_{\ell<h}q^{2\ell}\le h$. Hence $K_p\le(d^2+6d+1)\norm a_1^2$ and $H_p\le4\norm a_1\sum_{h=1}^d\sqrt{h(h+1)}\le2d(d+2)\norm a_1$, and since $2(d^2+6d+1)+4d^2(d+2)^2\le4(d+1)^4$, we get $2K_p+H_p^2\le4(d+1)^4\norm a_1^2$, which proves \eqref{main:eq:generalnormalizedmargin}.
\end{proof}
\begin{proof}[Proof of \cref{main:cor:excesscertificate}]
The proof of \eqref{si:eq:generalenergyinequality} did not use
$\norm{p(W)}<q^d$. Solving it as a quadratic inequality in $\sqrt L$ gives $L\ge\widehat\tau_p$, where
\begin{equation}\label{si:eq:generalsharptau}
 \widehat\tau_p=\frac{4\Delta_p^2}{\left(\sqrt{G_p^2+4K_p\Delta_p}+G_p\right)^2},
\end{equation}
with $K_p$ from \eqref{si:eq:generalK} and $G_p$ from \eqref{si:eq:generalG}. Its defining equality
$K_p\widehat\tau_p+G_p\sqrt{\widehat\tau_p}=\Delta_p$ and
$\Delta_p\le(\sum_i|\hat a_i|)^2$ imply
$0<\widehat\tau_p\le\Delta_p/K_p\le1/(d+1)<1$.
If $\norm{p(W)}<q^d$, then $E_p\le\Delta_p$ and
\(K_p\tau_p+G_p\sqrt{\tau_p}\le K_p\tau_p+\Delta_p/2+H_p^2\tau_p/2=\Delta_p\),
so monotonicity of $L\mapsto K_pL+G_p\sqrt L$ gives
$\widehat\tau_p\ge\tau_p$.
\end{proof}
In $K_p$, only the lag sums $C_h$ with positive sign need to enter (\cref{si:rem:onesidedmoments}).

\subsection{The peripheral splitting and shorter blocks}\label{si:sec:shorterblocks}
The results here supply the contraction decomposition used in
\cref{main:thm:normalalternative} and quantitative refinements of its
uniform block conclusion. Norms and adjoints refer to the working inner product, which, as in the main text, is written as the Euclidean one, so that adjoints are conjugate transposes; every assertion concerns ordinary single-point startup.

We work with the contraction $T=W/q$, whose summands in
\eqref{main:eq:peripheralsplit} are $U$ and $T_s=W_0/q$, and we use the index ${j_0}:=\min\{j\ge1:\|T_s^j\|<1\}$,
with ${j_0}=0$ when the stable summand is absent. The bound $\norm{W_0^t}<q^t$ of \cref{si:lem:normattained}, proved next, says that ${j_0}\le t$
when $T_s$ acts on a space of dimension $t\ge1$.\label{si:lem:peripheralsplitting}

\begin{proof}[Proof of the bound $\norm{W_0^t}<q^t$ in \cref{si:lem:normattained}]
The decomposition itself is proved in \cref{si:proof:normattained}; here we prove the bound on $W_0$.

Suppose $\|T_s^t x\|=\|x\|=1$. For
$D_s=(I-T_s^*T_s)^{1/2}$, contraction-defect telescoping gives $0=\|x\|^2-\|T_s^t x\|^2 =\sum_{j=0}^{t-1}\|D_sT_s^jx\|^2$.
By Cayley--Hamilton, $D_sT_s^jx=0$ for every $j\ge0$. Hence every power of
$T_s$ preserves the norm of $x$, contradicting $\rho(T_s)<1$. Since the
unit sphere is compact, this proves $\|T_s^t\|<1$.
\end{proof}

\begin{proposition}[A shorter uniform block]
\label{si:prop:shortblock}
Let $q=\|W\|\in(0,1)$ and $m\ge1$. Use the decomposition of $T=W/q$ in
\cref{si:lem:normattained}, and the index ${j_0}$ of \cref{si:lem:peripheralsplitting}. If $\rho(W)<q$, then the uniform block estimate in
\cref{main:thm:normalalternative} holds with ${k_0}={j_0}\le n$, at every memory.
If $\rho(W)=q$ and the moment equations have no solution, it holds with ${k_0}=\max\{{j_0},m+1\}$.
\end{proposition}

\begin{proof}
Suppose that, for the stated ${k_0}$, no constant $c<1$ works. There is then
a sequence of blocks with nonzero initial residual and normalized ratios
$\|r_{k+{k_0}}\|/(q^{k_0}\|r_k\|)$ tending to one. Normalize each block to have
$\|r_k\|=1$, and write $u_j=q^{-j}r_{k+j}$, $v_j=q^{-j}y_{k+j+1}$, $f_j=u_j-v_j$, $D=(I-T^*T)^{1/2}$.
Projection orthogonality and $u_{j+1}=Tv_j$ give the exact identity $\sum_{j=0}^{{k_0}-1}\bigl(\|f_j\|^2+\|Dv_j\|^2\bigr) =1-\|u_{k_0}\|^2\to0$.
Passing to a subsequence gives $u_0\to v$, $\|v\|=1$. The recurrence then
implies $u_j\to T^jv$ for $0\le j\le {k_0}$ and $v_j\to T^jv$ for $j<{k_0}$.
Monotonicity of the normalized norms shows $\|T^jv\|=1$ for every
$0\le j\le {k_0}$.

If $\rho(W)<q$, then ${k_0}={j_0}$ and $\|T^{k_0}\|<1$, a contradiction. Otherwise,
${k_0}\ge{j_0}$ forces the stable component of $v$ to vanish: the splitting is
orthogonal, the peripheral restriction is an isometry, and the stable
restriction of $T^{j_0}$ has norm strictly below one. If the stable summand
is absent, this conclusion is automatic.

Since ${k_0}\ge m+1$, projection orthogonality at the step indexed by $j=m$
can be applied to each retained residual $r_{k+m-\ell}$, $1\le\ell\le m$.
After rescaling, it reads $\langle v_m,u_{m-\ell}\rangle=q^\ell\|v_m\|^2$.
Taking limits, and using that $v$ is peripheral, gives $\langle U^m v,U^{m-\ell}v\rangle =\langle U^\ell v,v\rangle=q^\ell$ for $1\le\ell\le m$.
The squared components of $v$ on the orthogonal peripheral blocks are
nonnegative weights satisfying the prohibited moment equations. This
contradiction proves a uniform $c<1$.

The argument uses only the projection identities and compactness of
normalized initial residuals. It does not assume that the AA update is
continuous when the rank of a residual history changes.
\end{proof}

\begin{proposition}[An explicit delay and gain at the peripheral level]
\label{si:prop:explicitdelay}
Let $q=\|W\|\in(0,1)$ with $\rho(W)=q$, write $W=qU\oplus W_0$ as in
\cref{si:lem:normattained}, and, when $W_0$ is present, let ${j_0}$ be the
index above and $\sigma=\|(W_0/q)^{j_0}\|<1$. Let $p(z)=\sum_{j=0}^ma_jz^j$ be
real with $p(1)=1$ and $\eta_0:=\max\{|p(\lambda)|:\lambda\in\spec(W),\ |\lambda|=q\}<q^m$. Fix $\eta'$ with $\eta_0\le\eta'<q^m$ and $\eta'>0$, and put
\[
 A_p=\sum_{j=0}^m|a_j|q^j,\qquad
 B_p=\Bigl\|\Bigl(\sum_{j=i+1}^m|a_j|q^j\Bigr)_{i=0}^{m-1}\Bigr\|,\qquad
 s={j_0}\Bigl\lceil\frac{\log^+(A_p/\eta')}{\log(1/\sigma)}\Bigr\rceil,
\]
where $\log^+t=\max\{0,\log t\}$ and the displayed formula for $s$ applies when $0<\sigma<1$. Set $s=j_0$ when $W_0$ is present and $\sigma=0$, and $s=0$ when $W_0$ is absent. Then $P(z)=z^sp(z)$ satisfies
$\|P(W)\|\le q^s\eta'$, and every ordinary AA($m$) trajectory has $\rhAA(x_0)\le q\bigl[1-(q^m-\eta')^2/(2(s+m+1)\,(q^{2m}+sA_p^2+B_p^2))\bigr]<q$.
\end{proposition}
\begin{proof}
On the peripheral summand, $qU$ is normal and its eigenvalues are those of $W$ of modulus $q$, so $\|(qU)^sp(qU)\|=\max_\lambda|\lambda^sp(\lambda)|=q^s\eta_0\le q^s\eta'$, the maximum taken over these eigenvalues. On the stable summand,
$\|W_0^j\|\le q^j$ for every $j$ and $\|W_0^{k{j_0}}\|\le\|W_0^{j_0}\|^k=q^{k{j_0}}\sigma^k$,
while $\|p(W_0)\|\le A_p$. If $\sigma=0$, then $W_0^{j_0}=0$ and the required bound is immediate. Otherwise, with $s=k{j_0}$ and the stated $k$, this gives
$\|W_0^sp(W_0)\|\le q^s\sigma^kA_p\le q^s\eta'$. The splitting is orthogonal,
so $\|P(W)\|\le q^s\eta'<q^{s+m}$, and
\cref{main:prop:delayedcertificate} applies with $d=s+m$. The coefficients
of $P$ are those of $p$ shifted by $s$. Hence each of the first $s$ inner sums in the definition of $\beta$ in \eqref{main:eq:constants} equals $q^sA_p$,
the remaining ones are $q^s$ times those of $p$, and
$q^{2(s+m)}\beta^2=q^{2s}(sA_p^2+B_p^2)$, while the number $\eta$ of \eqref{main:eq:constants} for $P$ satisfies $\eta\le\eta'/q^m$. Inserting this into \eqref{main:eq:explicitrate} with $d=s+m$, the factors
$q^{2s}$ cancel.
\end{proof}

The delay grows only like $\log(A_p/\eta')/\log(1/\sigma)$, and the gain
in the displayed estimate decreases like $s^{-2}$ as $s$ grows, with $q,m,p,\eta'$ fixed. For a normal matrix, \cref{main:thm:normalalternative}
instead chooses $p$ small on the whole spectrum and needs no delay.

\begin{proposition}[An explicit bound when the operator norm exceeds the spectral radius]
\label{si:prop:nonnormaloidexplicit}
Suppose $\rho(W)<q=\|W\|<1$. Set $T=W/q$ and choose an integer $d\ge1$
with $a=\|T^d\|<1$. Such a $d$ exists with $d\le n$. Define $\beta=(\sum_{h=1}^{d}\|T^h\|^2)^{1/2}$, $c_d=(a+\beta\sqrt{1+\beta^2-a^2})/(1+\beta^2)$.
Then $0<c_d<1$, and every ordinary AA trajectory, at every memory
$m\ge1$, satisfies $\|r_{k+d}\|\le c_d q^d\|r_k\|$ for $k\ge0$, and $\rhAA(x_0)\le q c_d^{1/d}<q$.
The simpler choice $\beta=\sqrt d$ gives a valid, possibly weaker bound.
\end{proposition}

\begin{proof}
The existence of $d\le n$ follows from
\cref{si:lem:normattained}, since the peripheral summand is absent.
Normalize a nonzero block by $\|r_k\|=1$, and use $u_j,v_j,f_j,D$ from
the preceding proof. Put $t=\|u_d\|$ and $L=1-t^2$. The energy identity
gives $\sum_{j=0}^{d-1}\|f_j\|^2\le L$, while the recurrence gives $u_d=T^d u_0-\sum_{j=0}^{d-1}T^{d-j}f_j$.
Therefore $t\le a+\sum_{j=0}^{d-1}\|T^{d-j}\|\,\|f_j\| \le a+\beta\sqrt{1-t^2}$.
The unique root in $(0,1)$ of $t=a+\beta\sqrt{1-t^2}$ is the displayed
$c_d$, so $t\le c_d$. Restore scale and iterate blocks. Finally,
$\|T^h\|\le1$ gives $\beta\le\sqrt d$, and increasing $\beta$ weakens
the scalar bound.
\end{proof}

\subsection{Normal and normaloid matrices}\label{si:sec:normalnormaloid}
For a normal matrix the proof of \cref{main:thm:normalalternative} is shorter, while the periodic chains of \cref{main:ex:markov} need the theorem in its normaloid generality. This subsection records both points.

\emph{What simplifies.} If $W$ is normal, the splitting \eqref{main:eq:peripheralsplit} is the spectral theorem, the block $W_0$ is normal with $\norm{W_0}=\rho(W_0)<q$, and $\norm{P(W)}=\max_{\lambda\in\spec(W)}|P(\lambda)|$ for every polynomial $P$. In the step (ii)$\Rightarrow$(i) the polynomial $p_\varepsilon$ is then itself a certificate: on the eigenvalues of modulus $q$ it satisfies $|p_\varepsilon|<q^m$ as shown there, and on the others $|p_\varepsilon(\lambda)|\le|\lambda|^m+O(\varepsilon)<q^m$ for small $\varepsilon$. No delay is needed, the step (i)$\Rightarrow$(iii) is \cref{main:thm:certificate} with blocks of length $m+1$, part (v) applies, and neither \cref{main:prop:delayedcertificate} nor the bound $\norm{W_0^t}<q^t$ of \cref{si:lem:normattained} is used. The same shortening holds whenever $W_0$ is absent or $\norm{W_0}<q$, with \cref{si:lem:normattained} in place of the spectral theorem and $\norm{p_\varepsilon(W_0)}\le\norm{W_0}^m+O(\varepsilon)$ in place of the eigenvalue bound. The separation step, the projection-energy argument behind \cref{main:thm:certificate}, and the explicit start in (iv)$\Rightarrow$(ii) are the same for every normaloid $W$.

\emph{Diagonalizable matrices may still need a delay.} When $\norm{W_0}=q$, a delay may be unavoidable. \Cref{si:rem:normhypothesis} shows this for the defective matrix $\diag(q,qJ)$; the following example is diagonalizable. For $0<c\le q^2$ let
\[
 W=\diag(q,W_0),\qquad W_0=\begin{pmatrix}0&\beta\\0&c\end{pmatrix},\qquad \beta=\sqrt{q^2-c^2}.
\]
Then $W$ has the distinct eigenvalues $q$, $c$ and $0$, so it is diagonalizable, and $W_0=(\beta,c)^T\mathbf e_2^T$ has norm $q$, so $W$ is normaloid. Condition (ii) holds at $m=1$, since the only phase is $0$ and $\cos0=1\ne q$. A polynomial $P(z)=1-a+az$ with $|P(q)|<q$ has $a>1$, and with $v=\mathbf e_2$ and $u=W_0v/q$,
\[
 \ip{u}{P(W_0)v}=(1-a)\frac cq+aq=q+(a-1)\Bigl(q-\frac cq\Bigr)\ge q,
\]
so $\norm{P(W_0)}\ge q$ and no certificate of degree one exists. The delayed certificate $P(z)=z(z-c)/(1-c)$ annihilates $W_0$, since $W_0^2=cW_0$, and $|P(q)|=q(q-c)/(1-c)<q^2$.

\emph{A normaloid chain that is not normal.} The discounted chains $qP$ of \cref{main:ex:markov} are normaloid in the inner product of the stationary distribution, but they need not be normal there. Four states suffice. Let
\[
 P=\begin{pmatrix}0&I\\ \tfrac12J&0\end{pmatrix},\qquad J=\begin{pmatrix}1&1\\1&1\end{pmatrix},
\]
so that states $1$ and $2$ move to states $3$ and $4$, and each of states $3$ and $4$ moves to state $1$ or $2$ with probability $\tfrac12$. The chain is irreducible with period $2$, and $P$ is doubly stochastic, so $\pi$ is uniform and $\ip\cdot\cdot_\pi$ is a multiple of the Euclidean inner product. In the orthonormal basis $(1,1,1,1)^T/2$, $(1,1,-1,-1)^T/2$, $(1,-1,0,0)^T/\sqrt2$, $(0,0,1,-1)^T/\sqrt2$,
\[
 qP=\diag\bigl(q,\,-q,\,qN\bigr),\qquad N=\begin{pmatrix}0&1\\0&0\end{pmatrix}.
\]
This is the splitting \eqref{main:eq:peripheralsplit} with $U=\diag(1,-1)$ and $W_0=qN$. Hence $\norm{qP}_\pi=q=\rho(qP)$, so $qP$ is normaloid; it is not normal, since $N$ is not; and, because $\norm{W_0}=q$, its norm is attained also off the peripheral subspace. The moment equations on the phases $0$ and $\pi$ have the solution $w_0=\frac{1+q}2$, $w_\pi=\frac{1-q}2$ for $m=1$ and none for $m=2$, since $\cos0=\cos2\pi=1\ne q^2$. So AA(2) beats $q$ uniformly, as \cref{main:thm:cusps} guarantees for $\nu=2$, while by \cref{main:thm:normalalternative} some start makes AA(1) reproduce fixed-point iteration. Three states do not suffice: an irreducible three-state chain of period $3$ is a cyclic permutation, hence orthogonal, and one of period $2$ has a star-shaped transition graph, so it is reversible and $qP$ is self-adjoint in $\ip\cdot\cdot_\pi$.

\subsection{How often does a certificate exist?}\label{si:sec:howoften}
Delayed certificates give a complete qualitative test in
\cref{main:sec:normal}. For explicit rates, it is useful to know when a
certificate of small degree already works. At degree one there is an exact
test, and it succeeds for almost every contraction.

\begin{proposition}[An exact test for degree-one certificates]\label{si:prop:generic}
Let $W\in\R^{n\times n}$ with $0<q=\norm W<1$. A real polynomial $p$ of
degree at most one with $p(1)=1$ and $\norm{p(W)}<q$ exists if and only if
there is no unit vector $v$ with $\norm{Wv}=q$ and $\ip v{Wv}=q^2$.
The matrices with $\norm W<1$ that admit such a certificate form an open
and dense set of full Lebesgue measure. If the largest singular value of
$W$ is simple, with $Wv=qu$ for unit vectors $u,v$, the condition reads
$\ip uv\ne q$, and then $p_a(z)=a+(1-a)z$ satisfies
$\norm{p_a(W)}=q+a(\ip uv-q)+O(a^2)$. For a real normal matrix the
condition fails exactly when the eigenvalues $qe^{\pm i\theta_i}$ of modulus
$q$ carry weights $w_i\ge0$ summing to one with
$\sum_iw_i\cos\theta_i=q$, in agreement with
\cref{main:thm:normalalternative} for $m=1$.
\end{proposition}
\begin{proof}
Every normalized polynomial of degree at most one has the form
$p_a(z)=z+a(1-z)$. If a unit vector $v$ as in the statement exists, then $\norm{p_a(W)v}^2=\norm{Wv}^2+2a\bigl(\ip v{Wv}-\norm{Wv}^2\bigr)+a^2\norm{(I-W)v}^2 =q^2+a^2\norm{(I-W)v}^2\ge q^2$,
so no $p_a$ is a certificate.

Conversely, let $E=\ker(W^TW-q^2I)$ be the span of the right singular
vectors for $q$, so that $\norm{Wv}=q$ for unit $v\in E$. The values
$\ip v{Wv}-q^2$ over the unit sphere of $E$ are the Rayleigh quotients of
the compression of $(W+W^T)/2-q^2I$ to $E$ and therefore fill an
interval. If this interval does not contain $0$, choose $s\in\{-1,1\}$
with $s(\ip v{Wv}-q^2)<0$ for every unit $v\in E$. For
$F(t)=\norm{W+st(I-W)}^2=\max_{\norm v=1}\norm{(W+st(I-W))v}^2$, the
maximizers at $t=0$ are the unit vectors of $E$, and by Danskin's theorem
the right derivative at $0$ is $F'(0+)=2\max_{v\in E,\ \norm v=1}s\ip{Wv}{(I-W)v} =2\max_{v\in E,\ \norm v=1}s\bigl(\ip v{Wv}-q^2\bigr)<0$.
Hence $p_{st}$ is a certificate for all sufficiently small $t>0$.

The set of matrices with a certificate is open, since
$\norm{W+a(I-W)}<\norm W$ is an open condition. Its complement lies in the
union of the matrices with a repeated largest singular value, which are
contained in the zero set of the discriminant of the characteristic
polynomial of $W^TW$, a nonzero polynomial in the entries of $W$, and of
the matrices with a simple largest singular value and $\ip uv=\norm W$.
Along a ray $\{sW:s>0\}$ the vectors $u,v$ do not change while
$\norm{sW}$ does, so the latter set meets each ray in at most one point,
and integration in polar coordinates shows that it has measure zero. A
set of full measure is dense. When the largest singular value is simple,
$E=\operatorname{span}\{v\}$ and $\ip v{Wv}=q\ip uv$, which gives the stated
form of the condition, and the expansion of $\norm{p_a(W)}$ is the
first-order perturbation of a simple singular value, $\sigma_{\max}(W+aA)=\sigma_{\max}(W)+a\,u^TAv+O(a^2)$ for $A=I-W$, which follows from the analytic perturbation theory of the simple eigenvalue $q$ of the symmetric matrix $\begin{psmallmatrix}0&W\\W^T&0\end{psmallmatrix}$, with
$\ip u{(I-W)v}=\ip uv-q$. For a real normal matrix, $E$ is spanned by the
real and imaginary parts of the eigenvectors of the eigenvalues of modulus
$q$; for a unit vector $v\in E$ with squared components $w_i$ on the
corresponding real eigenblocks, $\ip v{Wv}=q\sum_iw_i\cos\theta_i$, and
every choice of weights is realized by some $v$.
\end{proof}

The margin controls the guaranteed gain, and two effects make it small or
misleading. For a normal matrix, $\norm{p(W)}$ is the maximum of $|p|$ over
the eigenvalues, so the margin is small when the spectrum comes close to
the circle $|z|=q$ in many directions. If the spectrum fills the disk
$|z|\le q$, as can happen for operators or in the limit of large matrices,
no polynomial does better than $z^d$: if $|p|<q^d$ on $|z|=q$, the maximum
principle applied to $p(z)/z^d$ on $|z|\ge q$ gives $|p(1)|<1$. On the
interval and on the deltoid, fixed-degree certificates have margins of order
$1-q$ (\cref{main:sec:sqrtstatus}). For a nonnormal matrix the margin can be
large for a different reason: $\norm{W^d}$ may be much smaller than
$\norm W^d$, so that $z^d$ itself is a certificate, and the guaranteed
improvement over $q$ then reflects the gap between $\norm W$ and the
asymptotic factor $\rho(W)$ of fixed-point iteration
(\cref{si:sec:spectralradius}).

\subsection{Proof of \texorpdfstring{\cref{main:thm:cusps}}{the memory theorem for equally spaced eigenvalues}}\label{si:sec:cuspsproof}
\begin{proof}[Proof of \cref{main:thm:cusps}]
For $\theta\in(2\pi/\nu)\mathbb Z$ and every integer $\ell$ we have
$\cos(\ell\theta)=\cos((\nu-\ell)\theta)$. Put $\ell=\lfloor\nu/2\rfloor+1$
and $\ell'=\nu-\ell=\lceil\nu/2\rceil-1$, so that $0\le\ell'<\ell\le m$.
The moment equations \eqref{main:eq:momentcondition} for $\ell$ and $\ell'$, the latter reading $\sum_iw_i=1$ when $\ell'=0$, have the same
left-hand side and the different right-hand sides $q^\ell\ne q^{\ell'}$.
No weights exist, and \cref{main:thm:normalalternative} gives the first
claim.

For sharpness let $\delta(j)=\min\{j,\nu-j\}$ be the distance from $j$ to $0$
on the cycle $\mathbb Z_\nu$, and set
\begin{equation*}
 \widetilde w_k=\frac1\nu\sum_{j=0}^{\nu-1}q^{\delta(j)}e^{2\pi ijk/\nu},
 \qquad k\in\mathbb Z_\nu.
\end{equation*}
Since $\delta(j)=\delta(\nu-j)$, the numbers $\widetilde w_k$ are real with
$\widetilde w_k=\widetilde w_{\nu-k}$, and Fourier inversion on $\mathbb Z_\nu$
gives $\sum_{k=0}^{\nu-1}\widetilde w_k\cos\bigl(2\pi\ell k/\nu\bigr) =q^{\delta(\ell)}=q^\ell$ for $0\le\ell\le\lfloor\nu/2\rfloor$,
which is \eqref{main:eq:momentcondition} for every
$m\le\lfloor\nu/2\rfloor$, together with $\sum_k\widetilde w_k=1$. The
$\widetilde w_k$ are nonnegative because $\nu\widetilde w_k$ are the
eigenvalues of the circulant matrix $C=(q^{\delta(j-k)})_{j,k\in\mathbb Z_\nu}$,
and $C$ is positive semidefinite: with $u_j\in\{0,1\}^\nu$ the indicator of
the cyclic window $\{j,j+1,\ldots,j+\lceil\nu/2\rceil-1\}$, one has
$\norm{u_j-u_k}^2=2\delta(j-k)$, so that $C_{jk}=q^{\norm{u_j-u_k}^2/2}$ is a
Gaussian kernel matrix, since $0<q<1$. Group the phases $2\pi k/\nu$ into
conjugate pairs, with weights $w_0=\widetilde w_0$, $w_k=2\widetilde w_k$
for $0<k<\nu/2$, and $w_{\nu/2}=\widetilde w_{\nu/2}$ when $\nu$ is even, and
take $W=q\diag(1,R_{2\pi/\nu},\ldots,R_{2\pi\lfloor(\nu-1)/2\rfloor/\nu})$,
extended by the diagonal entry $-q$ when $\nu$ is even. The construction in the proof of
\cref{main:thm:normalalternative} gives the initialization. For the last
claim, $\spec(W)\subseteq q\Gamma_\nu$ places the eigenvalues of modulus $q$
among the points $qe^{2\pi ik/\nu}$ by \eqref{main:eq:hypocycloidperiphery},
and the sharpness matrix has its spectrum in $q\Gamma_\nu$.
\end{proof}

\subsection{Peripheral configurations and the smallest memory}
\label{si:sec:peripheralconfigurations}
\Cref{si:tab:peripheral} collects what is known about the smallest memory for which AA($m$) beats $q=\norm W$ from every starting point, in terms of the eigenvalues of modulus $q$. Its rules are explained below.

\begin{table}[!tbp]
\centering\small
\renewcommand{\arraystretch}{1.12}
\begin{tabular}{@{}>{\raggedright\arraybackslash}p{0.40\linewidth}>{\raggedright\arraybackslash}p{0.38\linewidth}>{\raggedright\arraybackslash}p{0.16\linewidth}@{}}
\toprule
Eigenvalues of modulus $q$ & Smallest memory & Source\\
\midrule
\multicolumn{3}{@{}l}{\emph{(a) General rules}}\\
Any nonempty configuration & at most $N$ & \cref{main:prop:peripheralcount}\\
Any nonempty configuration, unless the phases satisfy \eqref{main:eq:peripheralexception} & at most $K$ & \cref{main:prop:peripheralcount}\\
All with real part above $q^2$, or all with real part below $q^2$ & $1$; otherwise at least $2$ & \cref{si:sec:peripheralconfigurations}\\
A subset of another configuration & at most what that configuration needs & \cref{si:sec:peripheralconfigurations}\\
Any, with given phases & decided exactly, for each $m$, by a linear feasibility problem in the weights & \cref{main:thm:normalalternative}\\
\midrule
\multicolumn{3}{@{}l}{\emph{(b) Explicit values}}\\
All $\nu$ points $qe^{2\pi ik/\nu}$ & exactly $\lfloor\nu/2\rfloor+1$ & \cref{main:thm:cusps}\\
Some of these $\nu$ points & at most $\lfloor\nu/2\rfloor+1$ & \cref{main:thm:cusps}\\
$q$ only, or $-q$ only & $1$ & \cref{main:prop:peripheralcount}; symmetric case \cite[Thm.~3.1]{garner2026}\\
$q$ and $-q$ & $2$ & \cite[Prop.~3.3]{garner2026}, \cref{main:thm:normalalternative}\\
One pair, with phase $\theta$ & $1$; but $2$ if $\cos\theta=q$ & \cref{si:sec:peripheralconfigurations}\\
$q$ and a pair & $1$ if $\cos\theta>q$; else $2$, but $3$ if $\cos\theta=-(1-q)/2$ & \cref{si:sec:peripheralconfigurations}\\
$-q$ and a pair & $1$ if $\cos\theta<q$; else $2$, but $3$ if $\cos\theta=(1+q)/2$ & \cref{si:sec:peripheralconfigurations}\\
Two pairs, $\cos\theta_1>\cos\theta_2$ & $1$ if $q\notin[\cos\theta_2,\cos\theta_1]$; else $2$, but $3$ if $(\cos\theta_1-q)(\cos\theta_2-q)=-(1-q^2)/2$, and $4$ at one point of that curve & \cref{si:sec:peripheralconfigurations}\\
\bottomrule
\end{tabular}
\caption{The smallest memory $m$ for which AA($m$) beats $q=\norm W$ from
every starting point, for a real matrix $W$. By
\cref{main:thm:normalalternative} only the eigenvalues of modulus $q$ enter,
through their phases $\theta_i\in[0,\pi]$: one for each of $\pm q$ that
occurs and one for each complex pair $qe^{\pm i\theta_i}$. $N$ is the number
of these eigenvalues and $K$ the number of phases. For a normal $W$, $m$ is
also the smallest degree of a certificate.
\Cref{si:sec:peripheralconfigurations} explains how the two parts fit together.}
\label{si:tab:peripheral}
\end{table}

By \cref{main:thm:normalalternative}, for every real $W$ with $\norm W=q$ only the eigenvalues of modulus $q$ matter. They are
$q$, $-q$, and complex pairs $qe^{\pm i\theta}$ with $0<\theta<\pi$, and we
describe them by phases $\theta_i\in[0,\pi]$, one for each pair and one for
each of $\pm q$ that occurs. Thus $N$ is $2K$ minus the number of the eigenvalues $\pm q$ present.

Part (a) of the table applies to every real matrix. Its first two rows are
\cref{main:prop:peripheralcount}. The third row is the case $m=1$ of
\cref{main:thm:normalalternative}: AA(1) fails to beat $q$ from every start exactly when some weighted average of the numbers $\cos\theta_i$
equals $q$, that is, when $q$ lies between the smallest and the largest of
them. An eigenvalue $qe^{i\theta}$ has real part $q\cos\theta$, so this
happens exactly when some of these eigenvalues have real part at least
$q^2$ and some have real part at most $q^2$; geometrically, when some of
them lie in the closed sector $|\arg z|\le\arccos q$ and some in the closed
sector $|\arg z|\ge\arccos q$. The fourth row holds because a polynomial that is small at all the
eigenvalues of modulus $q$ of a matrix is also small at any part of them;
for example, a matrix whose eigenvalues of modulus $q$ are only some of the
$\nu$ points of \cref{main:thm:cusps} needs at most $\lfloor\nu/2\rfloor+1$.
The last row says that for given phases the exact value is always
computable, since the moment equations of \cref{main:thm:normalalternative} form a linear program with $K$ unknowns and $m+1$ equations.

Part (b) lists the configurations whose value we know in closed form. For
$K\le2$ these are the rules of part (a) worked out. Take $q$ together with a
pair with phase $\theta$, so that $K=2$ and $N=3$. The third rule gives memory
one exactly when $\cos\theta>q$. Otherwise the second rule gives memory two,
unless the phases satisfy \eqref{main:eq:peripheralexception}, which for this
configuration reads $\cos\theta=-(1-q)/2$; there memory two fails, and the
first rule gives three. The other rows follow in the same way. The two rows
with real eigenvalues only are the symmetric case, treated
in~\cite[Theorem~3.1 and Proposition~3.3]{garner2026}.

The largest values in part (b) are instances of a general fact: for each
$N$, and each admissible choice of which of $\pm q$ occur, there is exactly one choice
of phases for which memory $N$ is needed. To describe it, say that points
$x_1,\ldots,x_K\in[-1,1]$ with weights $w_i>0$ form a quadrature rule of
degree $D$ for a probability measure $\mu$ if $\sum_iw_if(x_i)=\int f\,d\mu$
for every polynomial $f$ of degree at most $D$. Since $\cos(\ell\theta)$ is a
polynomial of degree $\ell$ in $\cos\theta$, the moment equations of \cref{main:thm:normalalternative} say precisely that the numbers
$\cos\theta_i$ carry a rule of degree $m$ for the measure $\mu_q$ obtained
from the Poisson density $(1-q^2)/\bigl(2\pi(1-2q\cos\theta+q^2)\bigr)$ by
$x=\cos\theta$. With $K$ points the largest possible degree is $2K-1$, or
$2K-2$ and $2K-3$ when one or both endpoints $\pm1$ are prescribed, and it is
reached only by the Gauss, Gauss--Radau, and Gauss--Lobatto rules. Hence
memory $N$ is needed exactly when the numbers $\cos\theta_i$ are the nodes
of the corresponding rule for $\mu_q$. For $K=2$ these nodes are the
exceptional values in part (b). For $q$, $-q$ and one pair, the Lobatto
node is $\cos\theta=q/2$, which needs memory four, whereas the equally
spaced choice $\cos\theta=0$ needs three.

Beyond part (b), the smallest memory is determined by a finite linear
program rather than a closed formula. For any normaloid matrix, the same
peripheral rules determine the smallest memory that accelerates every
initialization, by \cref{main:thm:normalalternative}; the existence of an
undelayed certificate of degree at most $m$ additionally depends on the
interior block. For a nonnormaloid contraction every memory improves on
$q=\norm W$, although this need not improve on $\rho(W)$.
\Cref{si:sec:spectralradius} treats that distinction and the choice of
inner product. Classes whose peripheral spectrum contains arbitrarily
many points need separate uniform estimates. Complex matrices use the
exponential moments of \cref{main:rem:complex}.

\section{Explicit rates on the deltoid}\label{si:sec:deltoidrates}
By \cref{main:thm:cusps} with $\nu=3$ and part (v) of \cref{main:thm:normalalternative}, the degree-two margin is positive for every normal matrix with spectrum in $q\Gamma_3$. This subsection computes the margin. It gives explicit quadratics, among them the one with the largest margin, and rates in which one quadratic certificate and one three-step factor serve every matrix in the class.

\begin{proposition}[Cusp-balanced quadratics]\label{si:prop:cuspbalanced}
For $t\ge0$, define $P_t(\zeta)=\zeta^2+t\zeta-t/(1+t)$, $M_t=(1+t+t^2)/(1+t)$.
Then $\max_{\zeta\in\Gamma_3}|P_t(\zeta)|=M_t$, and equality on
the boundary occurs exactly at the three cusps.
\end{proposition}
\begin{proof}
On the boundary, $|\gamma_3(\theta)|^2=(5+16\cos^3\theta-12\cos\theta)/9$,
$\operatorname{Re}\gamma_3(\theta)=(2\cos^2\theta+2\cos\theta-1)/3$, and
$\operatorname{Re}\gamma_3(\theta)^2=(8\cos^4\theta+4\cos\theta-3)/9$.
Expanding $|P_t(\gamma_3(\theta))|^2$ with these formulas gives
\begin{equation}\label{si:eq:quadraticfactorization}
 M_t^2-|P_t(\gamma_3(\theta))|^2
 =\frac{4(1-\cos\theta)(1+2\cos\theta)^2}{81(1+t)}\,\Sigma(\cos\theta,t),
\end{equation}
where $\Sigma$ is the polynomial $\Sigma(x,t)=9t^3+3(1+2x)^2t^2+(16x^3+12x^2-9x+8)t+2(8x^3-6x+7)$.
For $-1\le x\le1$ and $t\ge0$, every term of $\Sigma(x,t)$ is
nonnegative: the coefficient of $t$ equals
$27/4+(4x-1)^2(4x+5)/4\ge27/4$, and the constant term
equals $4(4x^3-3x)+14\ge10$ because $|4x^3-3x|\le1$ on $[-1,1]$. Hence
$\Sigma(\cos\theta,t)\ge10$, and \eqref{si:eq:quadraticfactorization}
gives the bound. Equality on the boundary requires $\cos\theta=1$ or $\cos\theta=-1/2$, which correspond to the three
cusps. The maximum-modulus
principle extends the bound to the filled region.
\end{proof}

\begin{proposition}[An explicit uniform quadratic certificate]\label{si:prop:quadratic}
For $q\in(0,1)$, the polynomial
\begin{equation}\label{si:eq:quadratic}
 p_2(z)=\frac{(4+2q)z^2+(2q^2+q^3)z-2q^3}
                 {4+2q+2q^2-q^3}
\end{equation}
satisfies $p_2(1)=1$ and
\begin{equation}\label{si:eq:quadraticbound}
 \max_{z\in q\Gamma_3}|p_2(z)|=q^2a_2,\qquad
 a_2:=\frac{4+2q+q^2}{4+2q+2q^2-q^3}<1.
\end{equation}
More generally, for $0<t<q$, the normalized polynomial
\begin{equation}\label{si:eq:quadraticfamily}
 p_{2,t}(z)=
 \frac{(1+t)z^2+t(1+t)qz-tq^2}{(1+t)(1+tq)-tq^2}
\end{equation}
satisfies $\max_{q\Gamma_3}|p_{2,t}|=q^2a_{2,t}$, where $1-a_{2,t}=t(1-q)(q-t)/(1+t+tq(1+t-q))>0$.
\end{proposition}
\begin{proof}
The denominator in \eqref{si:eq:quadraticfamily} is positive and equals
its numerator at $z=1$. For $z=q\zeta$, $p_{2,t}(q\zeta) =q^2(1+t)P_t(\zeta)/((1+t)(1+tq)-tq^2)$.
By \cref{si:prop:cuspbalanced}, its maximum modulus is
$q^2(1+t+t^2)/[(1+t)(1+tq)-tq^2]$.
Subtracting the numerator from the denominator gives
$t(1-q)(q-t)$, proving the margin formula.
Taking $t=q/2$ gives \eqref{si:eq:quadratic} and~\eqref{si:eq:quadraticbound}.
\end{proof}

Normality converts the scalar bound into the operator certificate
$\norm{p_2(W)}\le q^2a_2$. For $p_2$, the numbers $\eta$ and $\beta$ of \eqref{main:eq:constants} are $a_2$ (the ratio in \eqref{si:eq:quadraticbound}, not the leading coefficient of $p_2$) and $\beta_2$, where $\beta_2=\sqrt{(4+4q+q^2)^2+(4+2q)^2}/(4+2q+2q^2-q^3)$, $c_2=(a_2+\beta_2\sqrt{1+\beta_2^2-a_2^2})/(1+\beta_2^2)<1$.
Consequently, for every $m\ge2$, $\norm{r_{k+3}}\le c_2q^3\norm{r_k}$ for $k\ge0$, and $\rhAA(x_0)\le q c_2^{1/3}<q$.
The factor depends only on $q$, not on the dimension, the spectrum, the
initialization, or the memory $m\ge2$.

The energy estimate turns the same certificate into a linear-order
improvement as $q$ tends to one.
\begin{corollary}[A sharper quadratic deltoid rate]\label{si:cor:quadraticlinearmargin}
Let $W$ be real and normal with $\spec(W)\subseteq q\Gamma_3$. For every
$m\ge2$, ordinary exact AA($m$) satisfies $\norm{r_{k+3}}^2\le(1-\tau_2(q))q^6\norm{r_k}^2$, $\rhAA(x_0)\le q(1-\tau_2(q))^{1/6}<q$,
where $\tau_2(q)\in(0,1)$ is the explicit function
\eqref{si:eq:quadraticlineartau}. As $q\uparrow1$, $\tau_2(q)/(1-q)$
tends to $0.0197\ldots$, so the certified improvement
$q-q(1-\tau_2(q))^{1/6}$ is of order $1-q$.
\end{corollary}
\begin{proof}
We apply \cref{main:cor:excesscertificate} to the polynomial $p_2$ of
\eqref{si:eq:quadratic}, with $\norm{p_2(W)}$ replaced by its bound
$q^2a_2$ and with the sharper constant of supplementary
\cref{si:rem:sharpgram}. Put $\mathcal D(q)=4+2q+2q^2-q^3$. Direct
substitution gives $E_{p_2}=(1-q^2)q^6(q^2-4q+8)/(\mathcal D(q)^2)$, $\Delta_{p_2}=q^4+E_{p_2}-q^4a_2^2=8q^6(1-q)(2+q)/(\mathcal D(q)^2)$.
The lag coefficients satisfy $C_1=-C_2=8q^5(2+q)/\mathcal D(q)^2>0$, and
\eqref{si:eq:generalG} and~\eqref{si:eq:generalKsharp} give
\[
 \begin{split}
 K_{p_2}&=\frac{q^4}{\mathcal D(q)^2}
 \left[(2+q)^2(\sqrt2q+2\sqrt3)^2+8q^2(2+q)(3/2+\sqrt2)\right],\\
 G_{p_2}&=\frac{8q^5(2+q)\sqrt{1-q^2}}{\mathcal D(q)^2}
              \left[\sqrt2+\sqrt{3(1+q^2)}\right].
 \end{split}
\]
The maximum in \eqref{si:eq:generalKsharp} selects this expression for every
$0<q<1$. Indeed, $\mathcal D(q)\le(7/3)(2+q)$ and
$2(1-q^2)q^2(q^2-4q+8)\le4$. After removing the common factor
$q^4/\mathcal D(q)^2$, the Gram term alone exceeds $2(q^4+E_{p_2})$ by at
least $(12-98/9)(2+q)^2-4\ge4/9$. The factor of
\cref{main:cor:excesscertificate} is therefore
\begin{equation}\label{si:eq:quadraticlineartau}
 \tau_2(q)=\frac{4\Delta_{p_2}^2}
 {\left(\sqrt{G_{p_2}^2+4K_{p_2}\Delta_{p_2}}+G_{p_2}\right)^2}.
\end{equation}
As $q\uparrow1$,
\[
 K_{p_2}\longrightarrow\frac{162+36\sqrt6+24\sqrt2}{49},\qquad
 \frac{G_{p_2}}{\sqrt{1-q}}\longrightarrow\frac{48(1+\sqrt3)}{49},\qquad
 \frac{\Delta_{p_2}}{1-q}\longrightarrow\frac{24}{49}.
\]
Substituting these limits gives the exact positive constant
\[
 \lim_{q\uparrow1}\frac{\tau_2(q)}{1-q}
 =\left(1+\sqrt3+
 \sqrt{\frac{43}{4}+2\sqrt3+\frac32\sqrt6+\sqrt2}\right)^{-2}.
\]
Expanding $(1-\tau_2)^{1/6}$ proves the last assertion.
\end{proof}

The certified gain is about $0.0033(1-q)$; this is a guaranteed constant, not an estimate of the actual asymptotic gain. For memory at least three, the two-term cubic certificate $(27z^3-13q^3)/(27-13q^3)$ does much better and certifies a gain of about $0.103(1-q)$ (supplementary \cref{si:cor:cubiclinearmargin}). It is the unique best polynomial of the form $a_0+a_3z^3$ on $q\Gamma_3$, because $|27z^3-13|\le14$ on $\Gamma_3$ with equality only at the cusps and the arc midpoints, and it gives an explicit factor for every $m\ge3$ (supplementary \cref{si:lem:cubic,si:thm:cubic}).

The rational choice $t=q/2$ is convenient but does not minimize the
scalar polynomial norm. The minimization problem has an exact solution, even when complex
coefficients are allowed.

\begin{proposition}[The exact quadratic minimax value]\label{si:prop:quadraticminimax}
For every $q\in(0,1)$,
\[
 \inf_{\substack{p\in\C[z],\ \deg p\le2\\p(1)=1}}
       \max_{z\in q\Gamma_3}|p(z)|
 =\frac{3q^2}{1+q+q^2+2(1-q)\sqrt{1+q+q^2}}.
\]
The real polynomial $p_{2,t_{\mathrm{op}}}$ is the unique minimizer, where
$t_{\mathrm{op}}=q/(1+\sqrt{1+q+q^2})$.
\end{proposition}
\begin{proof}[Proof of \cref{si:prop:quadraticminimax}]
Let $\mathcal L_j$ be the Lagrange cardinal polynomials at $q\omega_3^j$,
$j=0,1,2$. Their values at one are $\mathcal L_0(1)=(1+q+q^2)/(3q^2)$, $\mathcal L_1(1)=(q^2+q\omega_3^2+\omega_3)/(3q^2)$, $\mathcal L_2(1)=\overline{\mathcal L_1(1)}$,
and $|\mathcal L_1(1)|=(1-q)\sqrt{1+q+q^2}/(3q^2)$.
Interpolation at these three points gives, for every admissible $p$, $1=\bigl|\sum_{j=0}^2\mathcal L_j(1)p(q\omega_3^j)\bigr|\le\bigl(1+q+q^2+2(1-q)\sqrt{1+q+q^2}\bigr)\max_{z\in q\Gamma_3}|p(z)|/(3q^2)$.
This proves the lower bound. Conversely, \cref{si:prop:cuspbalanced} gives $\max_{q\Gamma_3}|p_{2,t}|= q^2(1+t+t^2)/(1+t+tq(1+t-q))$.
Substituting $t=t_{\mathrm{op}}$ attains the lower bound. Equality in the interpolation inequalities fixes the values at all three
cusps, because each $\mathcal L_j(1)$ is nonzero. A quadratic is determined
by three values, which proves uniqueness even among complex coefficients.
\end{proof}

Writing the minimax value as $q^2a_{2,\mathrm{op}}$ yields $1-a_{2,\mathrm{op}}=(2/\sqrt3-1)(1-q)+O((1-q)^2)$.
Thus, uniformly over the class, no quadratic has a margin
$1-\max_{q\Gamma_3}|p|/q^2$ of larger than linear order in $1-q$, whatever its
coefficients. This is an approximation
result, not an optimality statement for the AA rate: the block estimate
also depends on the coefficient term $\beta$ of \eqref{main:eq:constants}, and the actual trajectory may
converge faster than its certified factor.

\section{Peripheral splitting and inner products}\label{si:sec:metricdynamics}

\subsection{Change of inner product}\label{si:sec:conjugation}
This is the change of coordinates stated in \cref{main:sec:geometry}.

\begin{proposition}[Change of inner product]\label{si:prop:conjugation}
Let $H=S^TS$ with $S$ real and invertible, and let $g$ be a map with fixed
point $x_*$. A sequence $(x_k)$ is an AA$_H$($m$) sequence for $g$ if and
only if $(Sx_k)$ is a Euclidean AA($m$) sequence for
$\tilde g(\tilde x)=Sg(S^{-1}\tilde x)$, with the same minimizing
coefficients at every step; the residuals correspond by $\tilde r_k=Sr_k$.
For $g(x)=Wx+b$, the conjugated map is $\tilde g(\tilde x)=SWS^{-1}\tilde
x+Sb$, and $\norm{p(W)}_H=\norm{Sp(W)S^{-1}}_2$ for every polynomial $p$.
The error sequences $x_k-x_*$ and $S(x_k-x_*)$ have the same root-linear factor, as do the residual sequences $r_k$ and $Sr_k$. For an affine map with $I-W$ invertible, the error and residual factors also agree.
\end{proposition}
\begin{proof}
Since $\tilde g(Sx)-Sx=S(g(x)-x)$, the residual of $\tilde g$ at $Sx$ is
$Sr(x)$, and $\norm{\sum_j\alpha_jr_j}_H=\norm{\sum_j\alpha_jSr_j}_2$, so the
two coefficient problems coincide, together with any constraint on the
coefficients, and $S\sum_j\alpha_jg(x_j)=\sum_j\alpha_j\tilde g(Sx_j)$. Norms
on $\R^n$ are equivalent, so root-linear factors agree.
\end{proof}

This section proves the peripheral splitting and the existence of adapted inner products, and shows that the norm-attainment hypothesis in the last part of \cref{main:thm:normalalternative} cannot be dropped.

\subsection{The normaloid decomposition}\label{si:proof:normattained}
\begin{lemma}[Peripheral splitting]\label{si:lem:normattained}
Let $W\in\R^{n\times n}$ with $q=\norm W>0$, and let $M$ be the real span of the real and imaginary parts of the eigenvectors of $W$ for its eigenvalues of modulus $q$. Then $M$ and $M^\perp$ are invariant under $W$ and $W^T$, and in an orthonormal basis adapted to $\R^n=M\oplus M^\perp$ the matrix $W$ has the form \eqref{main:eq:peripheralsplit}, where either summand may be absent. Every eigenvalue of modulus $q$ is semisimple, and $M=\{0\}$ exactly when $\rho(W)<q$. If $W_0$ acts on a space of dimension $t\ge1$, then $\norm{W_0^t}<q^t$. Finally, $\norm{Wx}=q\norm x$ holds only for $x\in M$ if and only if $W_0$ is absent or $\norm{W_0}<q$.
\end{lemma}
We prove the lemma when the peripheral summand is
present, including its last assertion, writing $V=W$ and $\rho=q$; the bound $\norm{W_0^t}<q^t$ is proved in \cref{si:lem:peripheralsplitting}.

\begin{proof}
Put $T=V/\rho$, a real matrix of norm one, and use the inner product
$\ip xy=y^*x$ on $\C^n$. Let $Tx=\lambda x$ with $|\lambda|=1$ and
$\norm x=1$. Then $\ip{T^*x}{\bar\lambda x}=\lambda\ip x{Tx}=\lambda\bar\lambda=1$,
so $\norm{T^*x-\bar\lambda x}^2 =\norm{T^*x}^2-2\operatorname{Re}\ip{T^*x}{\bar\lambda x}+1 \le1-2+1=0$,
and $x$ is a common eigenvector of $T$ and $T^*$. If also $Ty=\mu y$ with
$|\mu|=1$ and $\mu\ne\lambda$, then
$\lambda\ip xy=\ip{Tx}y=\ip x{T^*y}=\ip x{\bar\mu y}=\mu\ip xy$, so
$x\perp y$. Let $M_\C\subseteq\C^n$ be the span of all eigenvectors of $T$
for eigenvalues of modulus one. It is invariant under $T$ and $T^*$, hence
so is $M_\C^\perp$, and orthonormal bases of the eigenspaces form an
orthonormal eigenbasis of $M_\C$ in which $T|_{M_\C}$ is diagonal with
unimodular entries; thus $T|_{M_\C}$ is unitary. Since $V$ is real, $M_\C$
is closed under complex conjugation, so $M_\C=M+iM$ with
$M=M_\C\cap\R^n$, and $M_\C^\perp=M^\perp+iM^\perp$. Hence $M$ and
$M^\perp$ are invariant under $V$ and $V^T$, the restriction $V|_M=\rho T|_M$
is $\rho$ times an isometry of $M$, and $V|_{M^\perp}$ has no eigenvalue
of modulus $\rho$, because every eigenvector of $V$ for such an eigenvalue
lies in $M_\C$. For semisimplicity, let $|\lambda|=\rho$, $(V-\lambda)^2y=0$,
and write $y=y_1+y_2$ with $y_1\in M_\C$, $y_2\in M_\C^\perp$. Then
$(V-\lambda)y_2=(V-\lambda)y-(V-\lambda)y_1\in M_\C\cap M_\C^\perp=\{0\}$,
so $y_2=0$ because $\lambda$ is not an eigenvalue of $V|_{M_\C^\perp}$, and
$(V-\lambda)y_1=0$ because $V|_{M_\C}$ is diagonalizable. This proves
\eqref{main:eq:peripheralsplit} with $\norm{W_0}\le\norm V=\rho$. For
$x=x_1+x_2\in M\oplus M^\perp$, invariance gives
$\norm{Vx}^2=\rho^2\norm{x_1}^2+\norm{W_0x_2}^2$, so $\norm{Vx}=\rho\norm x$
for some $x\notin M$ exactly when $\norm{W_0}=\rho$. The converse is clear.
\end{proof}

\subsection{Adapted inner products}\label{si:sec:spectralradius}
\label{si:sec:weightedAA}
The exact criterion of \cref{main:thm:normalalternative} compares AA($m$) with fixed-point iteration whenever $W$ is normaloid, and whether $W$ is normaloid depends on the inner product, which is part of the algorithm (\cref{si:prop:conjugation}). This subsection characterizes when an inner product can make $W$ normaloid and thereby give an exact comparison with the worst-case fixed-point factor (\cref{si:cor:spectralbenchmark}). Since the inner product is now chosen, we assume only $\rho(W)<1$ rather than $\norm W<1$, unless a stronger hypothesis is stated.

We write $\ip xy_H=x^THy$ for the inner product of a symmetric positive
definite $H$, $\norm\cdot_H$ for its norm and the induced operator norm, and
AA$_H$($m$) for AA($m$) with its least-squares problem posed in
$\norm\cdot_H$, as in \cref{si:prop:conjugation}. Throughout, $S$ denotes
an invertible matrix with $S^TS=H$, so that $\norm W_H=\norm{SWS^{-1}}_2$.
For a matrix $W$ with $\rho(W)>0$, the \emph{peripheral subspace} $M_W$ is the real span of the real and imaginary parts of the eigenvectors of $W$ for its eigenvalues of modulus $\rho(W)$. When $\norm W=\rho(W)$ it is the subspace $M$ of \eqref{main:eq:peripheralsplit}; unlike $M$, it does not depend on the inner product, and the peripheral subspace of $SWS^{-1}$ is $SM_W$.

\begin{definition}[Adapted inner products]\label{si:def:adapted}
Let $\rho(W)>0$. The inner product of $H$ is \emph{adapted} to $W$ if
$\norm W_H=\rho(W)$ and $\norm{Wx}_H=\rho(W)\norm x_H$ if and only if $x\in M_W$.
\end{definition}

By \cref{si:lem:normattained}, $H$ is adapted to $W$ exactly when $SWS^{-1}$ has the form \eqref{main:eq:peripheralsplit}, with $\rho(W)$ in place of $q$ and with $W_0$ absent or $\norm{W_0}_2<\rho(W)$. Every inner product in which $W$ is normal
is adapted, but adaptedness does not require normality off the peripheral
subspace: the block $W_0$ may be nonnormal, even defective.

\begin{lemma}[Existence of adapted inner products]
\label{si:lem:adaptedexists}
Let $W\in\R^{n\times n}$ with $\rho=\rho(W)>0$.
\begin{enumerate}[label=(\alph*)]
\item For every $\varepsilon>0$ there is an inner product in which
      $\norm W\le\rho+\varepsilon$.
\item An inner product in which $\norm W=\rho$ exists if and only if every
      eigenvalue of $W$ of modulus $\rho$ is semisimple. In that case an
      adapted inner product exists. If $W$ is diagonalizable over $\C$, then
      $W$ is normal in some inner product.
\end{enumerate}
\end{lemma}
The proof, in \cref{si:proof:adaptedexists}, rescales the basis vectors of a real Jordan form of $W$~\cite[Section~3.4]{horn2013}.

Combining \cref{si:prop:conjugation} with \cref{main:thm:normalalternative}
gives the comparison with the spectral radius.

\begin{corollary}[The spectral radius as the benchmark]\label{si:cor:spectralbenchmark}
Let $W\in\R^{n\times n}$ with $\rho=\rho(W)\in(0,1)$, and let
$\rho e^{\pm i\theta_1},\ldots,\rho e^{\pm i\theta_{K_\rho}}$, $0\le\theta_i\le\pi$,
be its eigenvalues of modulus $\rho$.
\begin{enumerate}[label=(\alph*)]
\item If $\norm W_H=\rho$, every trajectory of AA$_H$($m$) has root-linear factor $\rhAA(x_0)\le\rho$. Every trajectory has $\rhAA(x_0)<\rho$, and then uniformly over blocks of a fixed length, if and only
      if the moment equations \eqref{main:eq:momentcondition} with $q=\rho$
      have no solution; otherwise some start makes AA$_H$($m$) reproduce
      fixed-point iteration exactly. Such an $H$ exists if and only if the
      eigenvalues of modulus $\rho$ are semisimple.
\item If $H$ is adapted to $W$, the conditions in (a) are also equivalent to
      the existence of a real polynomial $p$ of degree at most $m$ with
      $p(1)=1$ and $\norm{p(W)}_H<\rho^m$, and the uniform bound then holds
      over blocks of $m+1$ steps.
\item If $W$ is normal in the inner product of $H$, then $H$ is adapted to $W$,
      $\norm{p(W)}_H=\max_{\lambda\in\spec(W)}|p(\lambda)|$ for every
      polynomial $p$, and every result of
      \cref{main:sec:certificate,main:sec:normal} and of the supplement about
      AA for real normal matrices holds for AA$_H$($m$) applied to $W$, with
      $q=\rho(W)$ and the spectrum of $W$.
\end{enumerate}
\end{corollary}
\begin{proof}
By \cref{si:prop:conjugation}, AA$_H$($m$) for $W$ is Euclidean AA($m$) for
$V=SWS^{-1}$, with $\norm V_2=\norm W_H$, $\spec(V)=\spec(W)$, and
$\norm{p(V)}_2=\norm{p(W)}_H$. Part (a) is \cref{main:thm:normalalternative}
applied to $V$, together with \eqref{main:eq:basiccontraction} and
\cref{si:lem:adaptedexists}(b). Adaptedness means that the norm of $V$ is
attained only on its peripheral subspace, so the last part of
\cref{main:thm:normalalternative} gives (b). In (c), $V$ is normal, hence its
norm is attained only on its peripheral subspace and $\norm{p(V)}_2$ is the
maximum of $|p|$ over $\spec(V)$; the statements about normal matrices
involve only $q$, the spectrum, the norms $\norm{p(V)}$, and the algorithm.
\end{proof}

Condition (a) involves only the spectrum: whether AA$_H$($m$) beats
$\rho(W)$ from every start is the same for all inner products with
$\norm W_H=\rho(W)$, and \cref{main:thm:cusps,main:prop:peripheralcount}
apply to them with $q$ replaced by $\rho(W)$. 

For the block $E$ of \cref{main:prop:notbenchmark}, whose spectrum is a single conjugate pair with $\cos\theta=3/4\ne\rho(E)$, AA(1) therefore beats $\rho(E)$ from every start in every inner product in which $E$ is normal, although it is slower than fixed-point iteration in the Euclidean one.

Finding a suitable inner product requires information about $W$, but in several common cases one is known in closed form. If $W=D^{-1}B$ with $B$ symmetric and $D$ symmetric positive definite, as for the Jacobi iteration matrix $I-D^{-1}A$ of a symmetric $A$ with positive diagonal $D$, then $W$ is self-adjoint in the inner product of $D$; with $H=D$ the theory reduces to the symmetric case of~\cite{garner2026}. For a reversible Markov chain, the iteration matrix is self-adjoint in the inner product of the stationary distribution. For the discounted periodic chains of \cref{main:ex:markov}, the same inner product makes $W$ normaloid, though in general not normal, which suffices for \cref{si:cor:spectralbenchmark}(a). For a general matrix, an adapted inner product has to be built from spectral information, as in the proof of \cref{si:lem:adaptedexists}; once it is in place, the corollary shows that the only remaining obstruction is the moment equations on the phases of the eigenvalues of modulus $\rho(W)$.

\subsection{Existence of adapted inner products}\label{si:proof:adaptedexists}

\begin{proof}
(a) Let $W=XJX^{-1}$ be a real Jordan form~\cite[Section~3.4]{horn2013},
$J=D+E$, where $D$ is block diagonal with the real eigenvalues as $1\times1$
blocks and the blocks $\begin{psmallmatrix}\operatorname{Re}\lambda&-\operatorname{Im}\lambda\\\operatorname{Im}\lambda&\operatorname{Re}\lambda\end{psmallmatrix}$ for the pairs $\lambda,\bar\lambda$, and $E$ carries the identity blocks on the block
superdiagonal of each Jordan block. Then $D$ is normal with $\norm D=\rho$,
and $\norm E\le1$ because $E$ has at most one entry equal to one in every
row and column and no other nonzero entries. Let $\Delta_t$ be the positive
diagonal matrix that multiplies the $k$-th coordinate, or coordinate pair,
of each Jordan block by $t^k$. Then $\Delta_t^{-1}J\Delta_t=D+tE$, which has
norm at most $\rho+t$, and $S=(X\Delta_t)^{-1}$, $H=S^TS$, with
$t=\varepsilon$ give (a).

(b) If $\norm W_H=\rho$, then \cref{si:lem:normattained} applied to
$V=SWS^{-1}$ shows that its eigenvalues of modulus $\rho$, which are those
of $W$, are semisimple. Conversely, order the real Jordan form so that the
blocks of the eigenvalues of modulus $\rho$ come first. Since these
eigenvalues are semisimple, their blocks are $1\times1$ blocks $\pm\rho$ or
$2\times2$ blocks $\rho R_\theta$, so $J=\diag(\rho U_0,J_2)$ with $U_0$
orthogonal and $\rho(J_2)<\rho$. By the construction in (a), applied to $J_2$, there is $S_2$ with
$\norm{S_2J_2S_2^{-1}}<\rho$, and $S=\diag(I,S_2)X^{-1}$ gives
$SWS^{-1}=\diag(\rho U_0,S_2J_2S_2^{-1})$, so $H=S^TS$ is adapted by the
characterization after \cref{si:def:adapted}. If $W$ is diagonalizable,
$J_2$ is itself block diagonal with normal blocks, so $S=X^{-1}$ makes
$SWS^{-1}=J$ normal. Finally, if $W$ is normal in the inner product of $H$,
then $\norm W_H=\rho$, and $\norm{Wx}_H=\rho\norm x_H$ exactly on the
peripheral subspace; so $H$ is adapted.
\end{proof}

\subsection{The norm-attainment hypothesis}
\label{si:rem:normhypothesis}The last part of \cref{main:thm:normalalternative} assumes that the block $W_0$ in \eqref{main:eq:peripheralsplit} is absent or has norm below $q$, that is, that the norm of $W$ is attained only on the peripheral subspace (\cref{si:lem:normattained}), or equivalently that the working inner product is adapted in the sense of \cref{si:def:adapted}. This hypothesis is needed for the undelayed certificate (v), but not for the moment criterion. For the matrix $W=\diag(q,qJ)$ discussed after \cref{main:prop:delayedcertificate}, condition (ii) holds at $m=1$, but no certificate of degree one exists. Nevertheless AA(1) terminates after at most five steps from every initialization: $Wy=(qy_1,qy_3,0)^T$ for every $y$, so the third coordinate of every residual from $r_1$ on is zero; hence from $r_3$ on every residual is a multiple of $\mathbf e_1$, and $r_3\ne r_4$ unless both vanish, since $\norm{r_4}\le q\norm{r_3}$ by \eqref{main:eq:basiccontraction}; the affine hull of two distinct multiples of $\mathbf e_1$ contains the origin, giving $r_5=0$ at the latest. Thus (iv) holds while (v) fails; only delayed certificates exist, for instance $z^2(z-q)/(1-q)$, which annihilates $W$ because $J^2=0$.

\subsection{Invariant subspaces in a fixed inner product}\label{si:sec:invariantmoved}
\label{si:sec:euclideanobstruction}
We now fix the inner product, which we take to be the Euclidean one, and
consider matrices that are not normaloid in it. The following observation
underlies the constructions of \cref{main:sec:everymemoryobstruction}.

\begin{proposition}[Invariant subspaces]\label{si:prop:invariantsubspace}
Let $M$ be a $W$-invariant subspace. If $r_0\in M$, then every residual of
AA($m$) lies in $M$, and the residual sequence is the one that
AA($m$) generates for the restriction $W|_M$ with the Euclidean
inner product of $M$. Consequently, if $\norm{Wx}=\rho(W)\norm x$ for all
$x\in M$, and the phases of the eigenvalues of $W|_M$ carry weights with the
moment property \eqref{main:eq:momentcondition} of degree $m$ at
$q=\rho(W)$, then some initialization makes AA($m$) follow
fixed-point iteration exactly, with $\norm{r_k}=\rho(W)^k\norm{r_0}$.
\end{proposition}
\begin{proof}
If $r_{k-m_k},\ldots,r_k\in M$, their affine hull lies in $M$, so its
minimum-norm point $y_{k+1}$, which is the same whether computed in
$\R^n$ or in $M$, lies in $M$, and $r_{k+1}=Wy_{k+1}\in M$. For the second
claim, $W|_M=\rho(W)U$ with $U$ an isometry of $M$, hence orthogonal, so
the construction in the proof of \cref{main:thm:normalalternative} applies
to $W|_M$ and gives $r_0\in M$ from which AA($m$) for $W|_M$ follows
fixed-point iteration. By the first claim this is the AA($m$) sequence for
$W$.
\end{proof}

\subsection{Verification of \texorpdfstring{\cref{main:prop:notbenchmark}}{the slow-AA example}}\label{si:sec:notbenchmark}
This section verifies \cref{main:prop:notbenchmark} from the specified start.

\begin{proof}[Verification of \cref{main:prop:notbenchmark}]
The eigenvalues of $W_\lambda$ are $\lambda$ and $(3\pm i\sqrt7)/6$, of modulus $2/3$, so $\rho(W_\lambda)=\lambda$ is simple, and $\norm{W_\lambda}=\max(\lambda,\norm E)=2\sqrt2/3$, since $E^TE=(1/9)\begin{psmallmatrix}5&3\\3&5\end{psmallmatrix}$ has eigenvalues $8/9$ and $2/9$. With a single eigenvalue of modulus $\lambda$, at phase zero, the moment equations read $1=\lambda^\ell$. The subspace $\{0\}\times\R^2$ is invariant and contains $r_0$, so every residual lies in it and each least-squares step is that of the restriction $E$; we write the residuals in $\R^2$, where $r_0=\mathbf e_1$. Then $r_1=Er_0=(1,2)^T/3$. The point of the line through $r_0$ and $r_1$ nearest the origin is $y_2=(1,1)^T/2$, so $r_2=Ey_2=(0,2)^T/3$. The line through $r_1$ and $r_2$ is horizontal, so its point nearest the origin is $r_2$ itself, and $r_3=Er_2=(-2,4)^T/9$. The point of the line through $r_2$ and $r_3$ nearest the origin is $y_4=(-1,1)^T/3$, so $r_4=Ey_4=-(2/9)\mathbf e_1=\mu r_0$ with $\mu=-2/9$. The line through $r_3$ and $r_4$ is vertical, so its point nearest the origin is $r_4$, and $r_5=Er_4=\mu r_1$. Consecutive residuals are distinct, so each retained affine hull is a line and every minimizing coefficient vector is unique. The minimum-norm point of a scaled window is the scaled minimum-norm point, so if $(r_{k+4},r_{k+5})=\mu(r_k,r_{k+1})$, then $y_{k+6}=\mu y_{k+2}$ and $r_{k+6}=\mu r_{k+2}$. Induction from $k=0$ gives $r_{k+4}=\mu r_k$ for all $k\ge0$, hence $\rhAA(x_0)=(2/9)^{1/4}$ by \eqref{main:eq:errorequivalence}. Since the eigenvalues of $E$ are not real, fixed-point iteration from any nonzero start in this subspace converges at the rate $\rho(E)=2/3$, and $2/3<\lambda<(2/9)^{1/4}$ by the choice of $\lambda$.
\end{proof}

\subsection{Attaining the ceiling}\label{si:sec:attaining}
The next theorem gives a sufficient condition for attaining $w(W)$, in the same way as the norm in \cref{main:thm:normalalternative}, on an invariant subspace on which $W$ acts as $w(W)$ times an orthogonal map, where AA sees a normaloid matrix of norm $w(W)$; \cref{si:prop:spectralmeasure} below gives a necessary condition.

\begin{theorem}[Attaining the ceiling]\label{si:thm:numradequality}
Let $m\ge1$ and $w=w(W)>0$. If $W$ has an invariant subspace $M_w$ with $\norm{Wx}=w\norm x$ for all $x\in M_w$, so that $W$ acts on $M_w$ as $wU$ with $U$ orthogonal, and the moment equations \eqref{main:eq:momentcondition} of degree $m$ at level $w$ have a solution on the phases of the eigenvalues of $U$, then some starting point makes AA($m$) reproduce fixed-point iteration exactly, with $\norm{r_k}=w^k\norm{r_0}$ and $\rhAA(x_0)=w$.
\end{theorem}

\begin{proof}
The restriction $W|_{M_w}$ is normaloid with norm $w$, and the moment equations at level $w$ have a solution, so \cref{main:thm:normalalternative}, applied to $W|_{M_w}$, provides $r_0\in M_w$ from which AA($m$) for $W|_{M_w}$ reproduces fixed-point iteration with $\norm{r_k}=w^k\norm{r_0}$; by \cref{si:prop:invariantsubspace} this is the AA($m$) run for $W$.
\end{proof}

For normaloid $W$, with $M_w$ the subspace $M$ of \eqref{main:eq:peripheralsplit}, this is the reproduction statement of \cref{main:thm:normalalternative}. A subspace as in \cref{si:thm:numradequality} is spanned by eigenvectors whose eigenvalues have modulus $w(W)$, and since $w(W)\ge\rho(W)$ it exists only when $W$ is spectraloid. Call a nonterminating run \emph{scaled-periodic} if from some step on its residuals are periodic up to a scalar factor, $r_{k+p}=\mu r_k$ for all large $k$, as in \cref{main:prop:notbenchmark}; its factor is then $|\mu|^{1/p}$. For such runs the converse holds: if $|\mu|^{1/p}=w(W)$, then from some step on the run reproduces fixed-point iteration, and the span of its mixed points is an invariant subspace on which $W$ acts as $w(W)$ times an orthogonal map of finite order whose phases carry a solution of the moment equations (\cref{si:thm:selfsimilarfpi}). In particular no scaled-periodic run reaches the ceiling when $W$ is not spectraloid, which sharpens part (c) of \cref{main:cor:spectraloid}. Beyond spectraloid matrices the moment equations constrain runs at the ceiling in a different way: they are posed on the touching set $\Theta(W)$ of \eqref{main:eq:touchingset}, the angles at which the numerical range touches the circle of radius $w(W)$, and for a probability measure $\nu$ on $\Theta(W)$ they take the form \eqref{main:eq:measuremoments}.

\begin{proposition}[Runs at the ceiling]\label{si:prop:spectralmeasure}
Let $m\ge1$ and $w=w(W)>0$. If an ordinary AA($m$) run satisfies $\sum_kw^{-2k}\norm{y_k}^2=\infty$, as it does when $\norm{r_k}/w^k$ does not tend to zero, then some probability measure $\nu$ on $\Theta(W)$ satisfies \eqref{main:eq:measuremoments}. Consequently, if no probability measure on $\Theta(W)$ satisfies \eqref{main:eq:measuremoments}, then every ordinary AA($m$) run has $\sum_kw^{-2k}\norm{r_k}^2<\infty$, so $\norm{r_k}=o(w^k)$. This is the case when $m\ge2$ and $\Theta(W)$ has at most two points, and when $m=1$ and $\Theta(W)$ is a single real point or a pair $\{\theta_0,-\theta_0\}$ with $\cos\theta_0\ne w$.
\end{proposition}
The proof is in \cref{si:sec:spectralmeasureproof}; for normaloid $W$ the proposition is the necessity half of \cref{main:thm:normalalternative}. When the moment equations have no solution on $\Theta(W)$, in particular when $m\ge2$ and the numerical range touches the circle of radius $w(W)$ at a single point or a symmetric pair, it shows that $\norm{r_k}/w(W)^k\to0$ along every run, which sharpens the uniform bound $\norm{r_k}\le e\norm W\sqrt{k-1}\,w(W)^{k-1}\norm{r_0}$ that \eqref{main:eq:numradbound} gives with $\delta=1/(k-1)$. By itself this leaves room for a run that decays only slightly faster than $w(W)^k$, for instance like $w(W)^k/k$. \Cref{main:thm:uniformceiling} closes that gap: under the same hypothesis the factor is bounded below $w(W)$ uniformly over starts.

What \cref{main:thm:uniformceiling} leaves open is the case in which the touching set carries a solution of the moment equations. Of particular interest are nonspectraloid matrices, for instance with $m=1$ and $\Theta(W)=\{0,\pi\}$: no scaled-periodic run then reaches the ceiling (\cref{si:thm:selfsimilarfpi}), but whether some run keeps pace with $w(W)^k$ is unknown (\cref{si:rem:remaining}). \Cref{si:tab:scales} summarizes where AA($m$) stands relative to fixed-point iteration.

\begin{table}[t]
\centering\footnotesize\setlength{\tabcolsep}{4pt}
\begin{tabular}{@{}p{0.25\linewidth}p{0.13\linewidth}p{0.27\linewidth}p{0.29\linewidth}@{}}
\toprule
\raggedright Class of $W$ & \raggedright Scales & \raggedright Versus worst-case fixed-point factor & \raggedright Factor $w$ attained\tabularnewline
\midrule
\raggedright normaloid & \raggedright $\rho=w=\norm W$ & \raggedright never slower; faster from every start if and only if the moment equations on the eigenvalues of modulus $\norm W$ have no solution (\cref{main:thm:normalalternative}) & \raggedright if and only if these equations have a solution (\cref{main:thm:normalalternative})\tabularnewline[3pt]
\raggedright spectraloid, not normaloid & \raggedright $\rho=w<\norm W$ & \raggedright never slower (\cref{main:cor:spectraloid}) & \raggedright when $W$ has an invariant block $wU$ whose phases solve the equations (\cref{si:thm:numradequality}); every scaled-periodic run at the ceiling is of this kind (\cref{si:thm:selfsimilarfpi})\tabularnewline[3pt]
\raggedright not spectraloid & \raggedright $\rho<w<\norm W$ & \raggedright can be slower; factor at most $w$, or $w/\rho$ times the benchmark if $\rho>0$ (\cref{main:cor:spectraloid}) & \raggedright never by a scaled-periodic run (\cref{si:thm:selfsimilarfpi}); never, with a factor uniformly below $w$, if the equations have no solution on $\Theta(W)$ (\cref{main:thm:uniformceiling}); open otherwise\tabularnewline
\bottomrule
\end{tabular}
\caption{AA($m$) and fixed-point iteration in a fixed inner product, with $\rho=\rho(W)$ and $w=w(W)$. No run exceeds the ceiling $w$ (\cref{main:thm:numradceiling}), and the factor is uniformly below $w$ whenever the moment equations have no solution on $\Theta(W)$ (\cref{main:thm:uniformceiling}); the last column says when some run has factor exactly $w$.}
\label{si:tab:scales}
\end{table}

\subsection{Scaled-periodic runs at the ceiling}\label{si:sec:numradproofs}
A nonterminating run is \emph{scaled-periodic} if its residuals are eventually periodic up to a scalar factor, $r_{k+p}=\mu r_k$ for all $k\ge k_0$, for some period $p\ge1$ and some real $\mu\ne0$; its factor is $|\mu|^{1/p}$. The slow orbit of \cref{main:prop:notbenchmark} is of this kind. This subsection shows that a scaled-periodic run at the ceiling, $|\mu|^{1/p}=w(W)$, is fixed-point iteration on an invariant block of the kind in \cref{si:thm:numradequality}. Throughout, $w=w(W)>0$, $A=W/w$, so that $w(A)=1$, $\tilde y_k=w^{-k}y_k$, and $\ip ab=a^*b$ on $\C^n$. The first step records what the orthogonality relations of the least-squares step alone say about the Fourier modes of such a run; the theorem that follows adds the membership of each mixed point in the affine hull of its window, which the orthogonality relations do not see.

\begin{proposition}[Modes of a scaled-periodic run]\label{si:prop:selfsimilar}
Let $m\ge1$, and let an AA($m$) run be scaled-periodic, $r_{k+p}=\mu r_k$ for all $k\ge k_0$, with $|\mu|^{1/p}=w$. Put $p'=p$ if $\mu>0$ and $p'=2p$ otherwise, and $\omega_{p'}=e^{2\pi i/p'}$, in the convention $\omega_k=e^{2\pi i/k}$ of $\omega_{3}$. Then $\tilde y_{k+p'}=\tilde y_k$ for all $k\ge k_1:=k_0+m+1$, and the modes
\[
 Z_s=\sum_{k=k_1+1}^{k_1+p'}\omega_{p'}^{sk}\,\tilde y_k,\qquad s\in\mathbb Z/p'\mathbb Z,
\]
satisfy
\begin{equation}\label{si:eq:modeadjoint}
 \ip{Z_s}{AZ_s}=\omega_{p'}^{-s}\norm{Z_s}^2\qquad\text{and}\qquad A^TZ_s=2\omega_{p'}^{s}Z_s-\omega_{p'}^{2s}AZ_s .
\end{equation}
Thus every nonzero mode attains the numerical radius, $Z_s^*WZ_s=w\,\omega_{p'}^{-s}\norm{Z_s}^2$, and the normalized weights $\norm{Z_s}^2/\sum_t\norm{Z_t}^2$ at the phases $2\pi s/p'$ solve the moment equations \eqref{main:eq:momentcondition} of degree $m$ at level $w$.
\end{proposition}
\begin{proof}
Since $y_{k+1}$ is determined by the window $(r_{k-m},\dots,r_k)$, and the windows scale by $\mu$ after $p$ steps once all their residuals have index at least $k_0$, we have $y_{k+p}=\mu y_k$ for $k\ge k_1$; hence $\tilde y_{k+p'}=\tilde y_k$ for $k\ge k_1$, with $p'$ as in the statement. For $k\ge k_1+m$ and $0\le j\le m$ the residual $r_{k-j}=Wy_{k-j}$ lies in the window, and $y_{k+1}$ is the orthogonal projection of the origin onto its affine hull, so $\ip{y_{k+1}}{r_{k-j}}=\norm{y_{k+1}}^2$, that is,
\begin{equation}\label{si:eq:rescaledwindow}
 \ip{\tilde y_{k+1}}{A\tilde y_{k-j}}=w^j\norm{\tilde y_{k+1}}^2\qquad(k\ge k_1+m,\ 0\le j\le m).
\end{equation}
Extend $(\tilde y_k)_{k\ge k_1}$ periodically to all $k\in\mathbb Z$. Both sides of \eqref{si:eq:rescaledwindow} are then $p'$-periodic in $k$, so the relation holds for every $k\in\mathbb Z$; the modes invert to $\tilde y_k=p'^{-1}\sum_s\omega_{p'}^{-sk}Z_s$ for all $k$, and $Z_{-s}=\overline{Z_s}$ because the $\tilde y_k$ are real. Below, sums over $k$ run over one period and sums over $s$ over $\mathbb Z/p'\mathbb Z$. Orthogonality of the characters gives, for all integers $a$ and $b$,
\begin{equation}\label{si:eq:fourierbilinear}
 \sum_k\ip{\tilde y_{k+a}}{A\tilde y_{k+b}}=\frac1{p'}\sum_s\omega_{p'}^{s(a-b)}\ip{Z_s}{AZ_s},\qquad
 \sum_k\norm{\tilde y_k}^2=\frac1{p'}\sum_s\norm{Z_s}^2 .
\end{equation}
Summing \eqref{si:eq:rescaledwindow} over a period and using \eqref{si:eq:fourierbilinear} with $a=1$ and $b=-j$ gives
\begin{equation}\label{si:eq:modemoments}
 \sum_s\omega_{p'}^{s(j+1)}\ip{Z_s}{AZ_s}=w^j\sum_s\norm{Z_s}^2\qquad(0\le j\le m).
\end{equation}
For $j=0$ the left side is a sum of complex numbers of modulus at most $\norm{Z_s}^2$, because $w(A)=1$, and it equals the sum of these moduli; hence $\omega_{p'}^s\ip{Z_s}{AZ_s}=\norm{Z_s}^2$ for every $s$, which is the first relation in \eqref{si:eq:modeadjoint} and the attaining property $Z_s^*WZ_s=w\omega_{p'}^{-s}\norm{Z_s}^2$. Inserting it into \eqref{si:eq:modemoments} for $1\le j\le m$ gives $\sum_s\omega_{p'}^{sj}\norm{Z_s}^2=w^j\sum_s\norm{Z_s}^2$, whose real part is the moment equation of degree $j$ at level $w$ with the normalized weights $\norm{Z_s}^2/\sum_t\norm{Z_t}^2$ at the phases $2\pi s/p'$; the imaginary part vanishes because $\norm{Z_{-s}}=\norm{Z_s}$. Finally, let $H_s=(\omega_{p'}^sA+\omega_{p'}^{-s}A^T)/2$. For every unit vector $z$, $z^*H_sz=\operatorname{Re}(\omega_{p'}^sz^*Az)\le1$, so $H_s\preceq I$, while $Z_s^*H_sZ_s=\operatorname{Re}(\omega_{p'}^sZ_s^*AZ_s)=\norm{Z_s}^2$; hence $(I-H_s)Z_s=0$, which is the second relation in \eqref{si:eq:modeadjoint}.
\end{proof}

\begin{theorem}[Scaled-periodic runs at the ceiling are fixed-point iteration]\label{si:thm:selfsimilarfpi}
In the setting of \cref{si:prop:selfsimilar}, $r_k=y_{k+1}$ for all $k\ge k_1$: from step $k_1$ on, the run reproduces fixed-point iteration. The span $V$ of the mixed points $y_k$, $k\ge k_1$, is $W$-invariant, and $W|_V=wU$ with $U$ orthogonal and $U^{p'}=I$; the nonzero modes $Z_s$ are pairwise orthogonal eigenvectors of $U$ for the eigenvalues $\omega_{p'}^{-s}$, and the weights obtained by normalizing $\norm{Z_s}^2$ solve the moment equations. Thus $V$ is a subspace as in \cref{si:thm:numradequality}. In particular $W$ is spectraloid, and when $\rho(W)<w(W)$ no scaled-periodic run has factor $w(W)$.
\end{theorem}

\begin{proof}
We work with the periodic extension of $(\tilde y_k)_{k\ge k_1}$ to all $k\in\mathbb Z$ from the proof of \cref{si:prop:selfsimilar}, for which \eqref{si:eq:rescaledwindow} holds for every $k$, and write
\[
 g(a,b)=\ip{\tilde y_a}{\tilde y_b},\qquad \beta(a,b)=\ip{\tilde y_a}{A\tilde y_b},\qquad n_k=\norm{\tilde y_k}\qquad(a,b,k\in\mathbb Z).
\]
Only two properties of the run enter. The first is \eqref{si:eq:rescaledwindow}, which reads $\beta(k+1,k-j)=w^jn_{k+1}^2$ for $0\le j\le m$. The second is the membership of $y_{k+1}$ in the affine hull of its window, $y_{k+1}=\sum_{j=0}^m\alpha_j^{(k)}r_{k-j}$ with $\sum_j\alpha_j^{(k)}=1$, which after rescaling reads
\begin{equation}\label{si:eq:affine}
 \tilde y_{k+1}=\sum_{j=0}^m\alpha_j^{(k)}\,w^{-j}A\tilde y_{k-j},\qquad \sum_{j=0}^m\alpha_j^{(k)}=1 ;
\end{equation}
it holds for the run at every step $k\ge k_1+m$, hence for every $k\in\mathbb Z$ in the periodic extension, since the vectors involved are $p'$-periodic. Since the run does not terminate, $n_k>0$.

\emph{Step 1: the adjoint relation.} Inserting the second relation of \eqref{si:eq:modeadjoint} into $\tilde y_k=p'^{-1}\sum_s\omega_{p'}^{-sk}Z_s$ gives, for every $k\in\mathbb Z$,
\begin{equation*}
 A^T\tilde y_k=2\tilde y_{k-1}-A\tilde y_{k-2},
\end{equation*}
and hence, since $\beta(a,b)=\ip{A^T\tilde y_a}{\tilde y_b}$ and the vectors are real,
\begin{equation}\label{si:eq:betaflip}
 \beta(a,b)=2g(a-1,b)-\beta(b,a-2)\qquad(a,b\in\mathbb Z).
\end{equation}
This is the only place where $w(A)=1$ enters; the rest is algebra.

\emph{Step 2: the Gram band.} Taking $a=b+1$ in \eqref{si:eq:betaflip} and using \eqref{si:eq:rescaledwindow} with $j=0$ on both sides gives $n_{b+1}^2=2n_b^2-n_b^2$, so $n_k=c>0$ is constant; in particular $\norm{y_{k+1}}=w\norm{y_k}$ for $k\ge k_1$. Taking $a=b+2$ and $a=b$ in \eqref{si:eq:betaflip}, with \eqref{si:eq:rescaledwindow} for $j=1$, gives
\[
 wc^2=2g(b+1,b)-\beta(b,b),\qquad \beta(b,b)=2g(b-1,b)-wc^2,
\]
so $g(b+1,b)=g(b,b-1)$ for all $b$: the function $g(k,k-1)$ is a constant $\gamma_1$, and $\beta(b,b)=2\gamma_1-wc^2$. For $3\le\delta\le m+1$, taking $a=b+\delta$ in \eqref{si:eq:betaflip}, applying \eqref{si:eq:betaflip} once more to the reversed term $\beta(b,b+\delta-2)$, and using \eqref{si:eq:rescaledwindow} with $j=\delta-1$ twice,
\[
 w^{\delta-1}c^2=2g(b+\delta-1,b)-\bigl(2g(b-1,b+\delta-2)-w^{\delta-1}c^2\bigr),
\]
so $g(k,k-j)$ is a constant $\gamma_j$ for $2\le j\le m$ as well. The values follow from the moment equations of \cref{si:prop:selfsimilar}, $\sum_s\omega_{p'}^{sj}\norm{Z_s}^2=w^j\sum_s\norm{Z_s}^2$: as in \eqref{si:eq:fourierbilinear}, $\sum_s\omega_{p'}^{sj}\norm{Z_s}^2=p'\sum_kg(k,k-j)=p'^2\gamma_j$ and $\sum_s\norm{Z_s}^2=p'\sum_kn_k^2=p'^2c^2$. Hence
\begin{equation}\label{si:eq:gramband}
 g(k,k-j)=w^jc^2\quad(0\le j\le m),\qquad \beta(k,k)=wc^2\qquad(k\in\mathbb Z).
\end{equation}

\emph{Step 3: the defects are orthogonal to a band of mixed points.} Put $d_i=A\tilde y_i-\tilde y_{i+1}$, so that $\ip{d_i}{\tilde y_a}=\beta(a,i)-g(a,i+1)$. For $1\le a-i\le m+1$ both terms equal $w^{a-i-1}c^2$, by \eqref{si:eq:rescaledwindow} and \eqref{si:eq:gramband}; for $a=i$ both equal $wc^2$; and for $a=i-\delta$ with $1\le\delta\le m-1$, \eqref{si:eq:betaflip} gives $\beta(i-\delta,i)=2g(i-\delta-1,i)-\beta(i,i-\delta-2)=2w^{\delta+1}c^2-w^{\delta+1}c^2$, which again equals $g(i-\delta,i+1)$. Therefore
\begin{equation}\label{si:eq:bandorth}
 d_i\perp\tilde y_a\qquad\text{whenever }-(m-1)\le a-i\le m+1 .
\end{equation}

\emph{Step 4: the mixing coefficients vanish.} Fix $k$ and write $\alpha_j=\alpha_j^{(k)}$. Substituting $A\tilde y_{k-j}=\tilde y_{k-j+1}+d_{k-j}$ in \eqref{si:eq:affine} and using $\alpha_0=1-\sum_{j\ge1}\alpha_j$,
\begin{equation*}
 \sum_{j=1}^m\alpha_j\bigl(\tilde y_{k+1}-w^{-j}\tilde y_{k-j+1}\bigr)=\sum_{j=0}^m\alpha_jw^{-j}d_{k-j}.
\end{equation*}
Take the inner product with $\tilde y_{k+1-\ell}$ for $\ell\in\{1,\dots,m\}$. The right side vanishes by \eqref{si:eq:bandorth}, because $(k+1-\ell)-(k-j)=j+1-\ell$ lies in $[1-m,m+1]$ for $0\le j\le m$. On the left, \eqref{si:eq:gramband} gives $g(k+1,k+1-\ell)=w^\ell c^2$ and $g(k-j+1,k+1-\ell)=w^{|\ell-j|}c^2$, so, after division by $c^2$,
\[
 \sum_{j=1}^{\ell-1}\alpha_j\,w^{\ell}\bigl(1-w^{-2j}\bigr)+\bigl(w^{\ell}-w^{-\ell}\bigr)\sum_{j=\ell}^{m}\alpha_j=0\qquad(1\le\ell\le m).
\]
For $\ell=1$ this reads $(w-w^{-1})\sum_{j=1}^m\alpha_j=0$, so $\sum_{j\ge1}\alpha_j=0$. If $\alpha_1=\dots=\alpha_{\ell-2}=0$ for some $2\le\ell\le m$, then $\sum_{j\ge\ell}\alpha_j=-\alpha_{\ell-1}$, and the equation for $\ell$ reduces to $\alpha_{\ell-1}w^{-\ell}(1-w^2)=0$, so $\alpha_{\ell-1}=0$. Thus $\alpha_1=\dots=\alpha_{m-1}=0$, then $\alpha_m=\sum_{j\ge1}\alpha_j=0$, and $\alpha_0=1$: $\tilde y_{k+1}=A\tilde y_k$ for every $k$, that is, $r_k=Wy_k=y_{k+1}$ for all $k\ge k_1$. The argument applies to every coefficient vector representing $y_{k+1}$, so it also shows that the residuals of each window are affinely independent; this is consistent with scaled periodicity, because a window containing two proportional residuals $r_{k-p}$ and $r_k$ would have the origin in its affine hull and the run would terminate, whereas $n_k>0$.

\emph{Step 5: the invariant block.} On $V=\operatorname{span}\{\tilde y_k:k\ge k_1\}$ we have $A\tilde y_k=\tilde y_{k+1}$, so $AV\subseteq V$ and $A^{p'}=I$ on $V$ by periodicity. Hence $A|_V$ is diagonalizable over $\C$ with eigenvalues among the $p'$-th roots of unity, and $AZ_s=\sum_k\omega_{p'}^{sk}\tilde y_{k+1}=\omega_{p'}^{-s}Z_s$. A unit eigenvector $z$ of $A$ for an eigenvalue $\lambda$ of modulus one attains the numerical radius, so the argument that gave the second relation in \eqref{si:eq:modeadjoint}, with $\bar\lambda$ in place of $\omega_{p'}^s$, yields $A^Tz=\bar\lambda z$; for eigenvectors $z,z'$ with distinct eigenvalues $\lambda\ne\lambda'$ of modulus one this gives $\lambda'\ip z{z'}=\ip{z}{Az'}=\ip{A^Tz}{z'}=\lambda\ip z{z'}$, so they are orthogonal. Thus $A|_V$ has an orthonormal eigenbasis with eigenvalues of modulus one, and being real it is orthogonal: $W|_V=wU$ with $U=A|_V$, $U^{p'}=I$, the nonzero modes are its pairwise orthogonal eigenvectors, and their weights solve the moment equations by \cref{si:prop:selfsimilar}. Since $U$ has eigenvalues of modulus one, $\rho(W)\ge w$, so $W$ is spectraloid, and no such run exists when $\rho(W)<w$.
\end{proof}

Together with \cref{si:thm:numradequality}, this characterizes the scaled-periodic runs at the ceiling: one exists if and only if $W$ has an invariant subspace on which it acts as $w(W)$ times an orthogonal map of finite order whose phases carry a solution of the moment equations of degree $m$ at level $w(W)$. The finite order is what scaled periodicity adds to the hypothesis of \cref{si:thm:numradequality}. The proof uses neither the dimension of $V$ nor any hypothesis on the attaining vectors of $W$, and \cref{si:prop:selfsimilar} alone could not give it: the orthogonality relations are also satisfied by sequences at the ceiling that are not runs (\cref{si:rem:remaining}).

\subsection{The spectral measure of a tail and the proof of \texorpdfstring{\cref{si:prop:spectralmeasure}}{the spectral-measure proposition}}\label{si:sec:spectralmeasureproof}
As before, $w=w(W)>0$, $A=W/w$, $\tilde y_k=w^{-k}y_k$ and $n_k=\norm{\tilde y_k}$. For $\theta\in\R$ let
\[
 H_\theta=\tfrac12\bigl(e^{i\theta}A+e^{-i\theta}A^T\bigr),\qquad
 \hat h(\theta)=\lambda_{\max}(H_\theta)=\max_{\norm z=1}\operatorname{Re}\bigl(e^{i\theta}z^*Az\bigr),
\]
so that $\zeta^*H_\theta\zeta=\operatorname{Re}(e^{i\theta}\zeta^*A\zeta)$ for $\zeta\in\C^n$. The function $\hat h$ is continuous, $\hat h\le1$ because $w(A)=1$, and $\{\hat h=1\}=\Theta(W)$: a point $\xi\in\mathcal W(A)$ with $\operatorname{Re}(e^{i\theta}\xi)=1$ and $|\xi|\le1$ equals $e^{-i\theta}$, and $e^{-i\theta}\in\mathcal W(A)$ means $we^{-i\theta}\in\mathcal W(W)$. Two facts about an ordinary AA($m$) run are used. First, $n_{t+1}\le\norm A\,n_t\le2n_t$, because $\norm{y_{t+1}}\le\norm{r_t}=\norm{Wy_t}$ and $\norm A\le2w(A)$. Second, for $b\ge1$ and $0\le j\le m$ the residual $r_b=Wy_b$ is retained in the window that produces $y_{b+j+1}$, and the orthogonality of the projection gives $\ip{y_{b+j+1}}{Wy_b}=\norm{y_{b+j+1}}^2$, that is,
\begin{equation}\label{si:eq:Rjtail}
 \ip{\tilde y_{b+j+1}}{A\tilde y_b}=w^j\,n_{b+j+1}^2\qquad(b\ge1,\ 0\le j\le m);
\end{equation}
for $j=0$ this is \eqref{main:eq:growthidentity}. Nothing else about the run enters this subsection and the next, so their results also hold for a run supplied with an arbitrary history, from the first index whose residual is computed from a mixed point.

For $1\le k\le L$ with $y_k\ne0$ let
\[
 \zeta_{k,L}(\theta)=\sum_{t=k}^Le^{it\theta}\tilde y_t,\qquad T_{k,L}=\sum_{t=k}^Ln_t^2,\qquad
 d\nu_{k,L}=\frac{\norm{\zeta_{k,L}(\theta)}^2}{2\pi\,T_{k,L}}\,d\theta ,
\]
so that $\nu_{k,L}$ is a probability measure on $(-\pi,\pi]$ by Parseval's identity, the \emph{spectral measure} of the tail of the run from index $k$ to index $L$; where $\zeta_{k,L}\ne0$ let $f_{k,L}=\zeta_{k,L}^*A\zeta_{k,L}/\norm{\zeta_{k,L}}^2$, a point of $\mathcal W(A)$, and $\delta_{k,L}=1-e^{i\theta}f_{k,L}$. The subscripts are dropped when $k$ and $L$ are fixed.

\begin{lemma}[Tail identities]\label{si:lem:tail}
For $1\le k\le L$ with $y_k\ne0$,
\begin{align}
 \frac1{2\pi}\int_{-\pi}^{\pi}\zeta^*\bigl(I-H_\theta\bigr)\zeta\,d\theta&=n_k^2,\label{si:eq:tailenergy}\\
 \frac1{2\pi}\int_{-\pi}^{\pi}e^{i(j+1)\theta}\,\zeta^*A\zeta\,d\theta&=w^j\bigl(T-e_j\bigr),\qquad
 e_j=\sum_{t=k}^{\min(k+j,L)}n_t^2\le\frac{4^{j+1}-1}3\,n_k^2\qquad(0\le j\le m).\label{si:eq:tailmoments}
\end{align}
Consequently $\operatorname{Re}\delta\ge1-\hat h\ge0$, $|\delta|^2\le2\operatorname{Re}\delta$, and
\begin{equation}\label{si:eq:tailestimates}
 \begin{gathered}
 \int(1-\hat h)\,d\nu\le\int\operatorname{Re}\delta\,d\nu=\frac{n_k^2}{T},\qquad
 \int|\delta|\,d\nu\le\frac{\sqrt2\,n_k}{\sqrt T},\\
 \Bigl|\int\cos(j\theta)\,d\nu-w^j\Bigr|\le\frac{e_j}{T}+\frac{\sqrt2\,n_k}{\sqrt T}\qquad(1\le j\le m).
 \end{gathered}
\end{equation}
\end{lemma}
\begin{proof}
Since the $\tilde y_t$ are real, $\zeta^*A\zeta=\sum_{a,b=k}^Le^{i(b-a)\theta}\ip{\tilde y_a}{A\tilde y_b}$, and averaging against $e^{i(j+1)\theta}$ keeps the pairs $a=b+j+1$, which by \eqref{si:eq:Rjtail} contribute $w^j\sum_{b=k}^{L-j-1}n_{b+j+1}^2=w^j(T-e_j)$; the bound on $e_j$ follows from $n_t\le2^{t-k}n_k$. This is \eqref{si:eq:tailmoments}. For \eqref{si:eq:tailenergy}, $\frac1{2\pi}\int\norm\zeta^2d\theta=T$ by Parseval's identity, and $\frac1{2\pi}\int\zeta^*H_\theta\zeta\,d\theta=\operatorname{Re}\frac1{2\pi}\int e^{i\theta}\zeta^*A\zeta\,d\theta=T-n_k^2$ by \eqref{si:eq:tailmoments} with $j=0$, where $e_0=n_k^2$ exactly. Next, $\operatorname{Re}\delta=1-\zeta^*H_\theta\zeta/\norm\zeta^2\ge1-\hat h$, and writing $e^{i\theta}f=a+ib$ with $a^2+b^2\le1$, $|\delta|^2=(1-a)^2+b^2\le(1-a)(1-a+1+a)=2\operatorname{Re}\delta$. The first estimate in \eqref{si:eq:tailestimates} is \eqref{si:eq:tailenergy} divided by $T$, the second follows from the Cauchy--Schwarz inequality, and the third from $\int e^{i(j+1)\theta}f\,d\nu=w^j(1-e_j/T)$, which is \eqref{si:eq:tailmoments} divided by $T$, and $e^{ij\theta}=e^{ij\theta}(e^{i\theta}f+\delta)$, by taking real parts.
\end{proof}

\begin{proof}[Proof of \cref{si:prop:spectralmeasure}]
The hypothesis $\sum_kw^{-2k}\norm{y_k}^2=\infty$ holds when $\norm{r_k}/w^k$ does not tend to zero, since $\norm{r_k}\le\norm W\norm{y_k}$. Take $k=1$ in \cref{si:lem:tail} and let $L\to\infty$, so that $T_{1,L}\to\infty$ while $n_1$ and the $e_j$ stay bounded. By \eqref{si:eq:tailestimates}, $\int(1-\hat h)\,d\nu_{1,L}\to0$ and $\int\cos(j\theta)\,d\nu_{1,L}\to w^j$ for $1\le j\le m$. The $\nu_{1,L}$ are probability measures on a compact set, so weak limits exist, and each weak limit $\nu$ satisfies $\int(1-\hat h)\,d\nu=0$, hence is supported on $\{\hat h=1\}=\Theta(W)$, and $\int\cos(j\theta)\,d\nu=w^j$ for $1\le j\le m$, which is \eqref{main:eq:measuremoments}. If no probability measure on $\Theta(W)$ satisfies \eqref{main:eq:measuremoments}, no run has $\sum_kw^{-2k}\norm{y_k}^2=\infty$, and then $\sum_kw^{-2k}\norm{r_k}^2\le\norm W^2\sum_kw^{-2k}\norm{y_k}^2<\infty$. Finally, let $\Theta(W)$ have at most two points; by symmetry it is $\{0\}$, $\{\pi\}$, $\{0,\pi\}$ or $\{\theta_0,-\theta_0\}$. In the first two cases $\int\cos\theta\,d\nu=\pm1\ne w$. In the third, $\int\cos2\theta\,d\nu=1\ne w^2$. In the fourth, $\int\cos\theta\,d\nu=\cos\theta_0$, and if this equals $w$ then $\int\cos2\theta\,d\nu=2w^2-1<w^2$. Hence \eqref{main:eq:measuremoments} fails for $m\ge2$ in all four cases, and for $m=1$ in the first two and in the fourth when $\cos\theta_0\ne w$. More generally, if $\Theta(W)$ is finite and $m\ge|\Theta(W)|$, no probability measure on $\Theta(W)$ satisfies \eqref{main:eq:measuremoments}, which is the first sufficient condition stated in \cref{main:thm:uniformceiling}: as in the proof of \cref{main:prop:peripheralcount}, let $Q$ be the product of $(x-\cos\theta)^2$ over the pairs $\{\theta,-\theta\}\subseteq\Theta(W)$ with $0<\theta<\pi$, multiplied by $1-x$ if $0\in\Theta(W)$ and by $1+x$ if $\pi\in\Theta(W)$. It has degree $|\Theta(W)|$, is nonnegative and not identically zero on $[-1,1]$, and vanishes at $\cos\theta$ for every $\theta\in\Theta(W)$, so $\int Q(\cos\theta)\,d\nu=0$ for every probability measure $\nu$ on $\Theta(W)$, whereas the moments $\int\cos(j\theta)\,d\nu=w^j$ for $0\le j\le m$ would give $\int Q(\cos\theta)\,d\nu=\frac1{2\pi}\int_0^{2\pi}Q(\cos\theta)P_w(\theta)\,d\theta>0$, because $\cos(j\theta)$ is a polynomial of degree $j$ in $\cos\theta$ and $P_w>0$.
\end{proof}

\subsection{A uniform bound below the ceiling, and explicit constants}\label{si:sec:uniformceilingproof}
This subsection proves \cref{main:thm:uniformceiling}, the analogue for the ceiling $w(W)$ of the uniform half of \cref{main:thm:normalalternative}: when the moment equations \eqref{main:eq:measuremoments} have no solution on the touching set, the factor is bounded below $w(W)$ uniformly over starting points, in the block form of \cref{main:thm:normalalternative}(iii) with $w(W)$ in place of $q$.

The sufficient conditions in the statement are the analogue of \cref{main:prop:peripheralcount} for the touching set; they are proved at the end of the proof of \cref{si:prop:spectralmeasure}. The proof of \cref{si:prop:spectralmeasure} passes to a weak limit and loses all quantitative information. The proof of \cref{main:thm:uniformceiling} keeps the spectral measure of the tail from an arbitrary index $k$ and tests it against the functional that separates the moment vector $(w^j)_{j=1}^m$ from the moment vectors of the probability measures on $\Theta(W)$. What makes the result uniform is that the energy identity \eqref{si:eq:tailenergy} of every tail has $n_k^2$, the energy at the start of that tail, on its right side. The constant $C$ comes from that separation and is not explicit; \cref{si:prop:touchingcertificate} below gives explicit constants when a polynomial vanishing on $\Theta(W)$ is available, and for the touching sets $\{0\}$ and $\{\pm\theta_0\}$ the polynomials $z-1$ and $z^2-2z\cos\theta_0+1$ serve under quadratic contact.

\begin{proof}[Proof of \cref{main:thm:uniformceiling}]
\emph{Separation.} The set
\[
 \mathcal K=\Bigl\{\Bigl(\int\cos(j\theta)\,d\mu\Bigr)_{j=1}^m:\ \mu\text{ a probability measure on }\Theta(W)\Bigr\}\subset\R^m
\]
is convex and compact, as the image of the weak-$*$ compact convex set of probability measures on the compact set $\Theta(W)$ under a continuous linear map, and $(w^j)_{j=1}^m\notin\mathcal K$ by hypothesis. Strict separation of a point from a compact convex set gives $h\in\R^m$ and $\eta>0$ with $\sum_jh_j(x_j-w^j)\le-2\eta$ for all $x\in\mathcal K$; for the point masses at $\theta\in\Theta(W)$ this says
\[
 \psi(\theta):=\sum_{j=1}^mh_j\bigl(\cos(j\theta)-w^j\bigr)\le-2\eta\qquad(\theta\in\Theta(W)).
\]
The set $N=\{\theta:\psi(\theta)<-\eta\}$ is open and contains $\Theta(W)$; its complement is compact and $\hat h<1$ on it, so the continuous function $1-\hat h$ has a positive minimum $c_N$ there. Also $|\psi|\le2\norm h_1$. The numbers $h$, $\eta$ and $c_N$ depend only on $W$ and $m$.

\emph{Two bounds for $\int\psi\,d\nu$.} Fix a run and $k\ge1$. If $y_k=0$ the run has terminated and \eqref{main:eq:uniformtail} is trivial, so let $y_k\ne0$, let $L\ge k$, and write $\nu=\nu_{k,L}$, $T=T_{k,L}$. The third estimate in \eqref{si:eq:tailestimates} and $e_j\le\frac{4^{m+1}}3n_k^2$ give
\[
 \int\psi\,d\nu=\sum_{j=1}^mh_j\Bigl(\int\cos(j\theta)\,d\nu-w^j\Bigr)\ge-\norm h_1\Bigl(\frac{4^{m+1}n_k^2}{3T}+\frac{\sqrt2\,n_k}{\sqrt T}\Bigr),
\]
while the definition of $N$ and the first estimate in \eqref{si:eq:tailestimates} give
\[
 \int\psi\,d\nu\le-\eta\,\nu(N)+2\norm h_1\,\nu(N^c)=-\eta+\bigl(\eta+2\norm h_1\bigr)\nu(N^c)\le-\eta+\bigl(\eta+2\norm h_1\bigr)\frac{n_k^2}{c_NT},
\]
because $\nu(N^c)\le c_N^{-1}\int_{N^c}(1-\hat h)\,d\nu\le c_N^{-1}n_k^2/T$. Put $x=n_k/\sqrt T$; then $0<x\le1$ because $T\ge n_k^2$, and the two bounds combine, using $x^2\le x$, to $\eta\le C_1x$ with
\[
 C_1=\norm h_1\Bigl(\frac{4^{m+1}}3+\sqrt2\Bigr)+\frac{\eta+2\norm h_1}{c_N}.
\]
Hence $T_{k,L}\le(C_1/\eta)^2n_k^2$ for every $L$, which is \eqref{main:eq:uniformtail} with $C=(C_1/\eta)^2$, and $C\ge1$ because $T_{k,k}=n_k^2$.

\emph{Consequences.} Let $S_k=\sum_{j\ge k}n_j^2$. Then $S_{k+1}=S_k-n_k^2\le(1-1/C)S_k$, so $n_{k+j}^2\le S_{k+j}\le(1-1/C)^jS_k\le C(1-1/C)^jn_k^2$, that is,
\[
 \norm{y_{k+j}}\le\sqrt C\,(1-1/C)^{j/2}\,w^j\norm{y_k}\qquad(k\ge1,\ j\ge0).
\]
With $k=1$, $y_1=r_0$ and $\norm{r_{j+1}}\le\norm W\norm{y_{j+1}}$ this is \eqref{main:eq:uniformrate}. For the block bound, $\norm{r_{k+k_0}}\le\norm W\norm{y_{k+k_0}}\le\norm W\sqrt C(1-1/C)^{(k_0-1)/2}w^{k_0-1}\norm{y_{k+1}}$ and $\norm{y_{k+1}}\le\norm{r_k}$ give $\norm{r_{k+k_0}}\le c\,w^{k_0}\norm{r_k}$ for all $k\ge0$ with $c=\frac{\norm W}{w}\sqrt C(1-1/C)^{(k_0-1)/2}\le2\sqrt C\,e^{-(k_0-1)/(2C)}$, which is below one as soon as $k_0-1>C\ln(4C)$.
\end{proof}

For normaloid $W$ the hypothesis of \cref{main:thm:uniformceiling} is condition (ii) of \cref{main:thm:normalalternative}, and the theorem gives (iii) by a route that does not pass through a certificate, at the price of an inexplicit block length. The constant can be made explicit when a polynomial of the kind the paper calls a certificate is available on the circle, where it is measured against the gap $1-\hat h$ rather than against the norm inside the disk. The first ingredient is an inequality that compares a mixed point with every affine combination of its window.

\begin{proposition}[Touching certificates]\label{si:prop:touchingcertificate}
Let $Q(z)=\sum_{j=0}^mq_jz^j$ be a real polynomial of degree at most $m$ with $Q(w)\ne0$, and write $Q(S^{-1})\tilde y_t=\sum_{j=0}^mq_j\tilde y_{t-j}$. Then every ordinary AA($m$) run satisfies the window inequality
\begin{equation}\label{si:eq:windowinequality}
 n_{t+1}\le\frac{\norm{A\,Q(S^{-1})\tilde y_t}}{|Q(w)|}\qquad(t\ge m+1).
\end{equation}
If moreover
\begin{equation}\label{si:eq:touchingcertificate}
 C_Q:=\sup_{\theta}\frac{|Q(e^{i\theta})|^2}{1-\hat h(\theta)}<\infty,
\end{equation}
then \eqref{main:eq:uniformtail}, and with it the conclusions of \cref{main:thm:uniformceiling}, hold with the explicit constant
\[
 C=\frac{4^{m+1}-1}{3}+\frac{4\,C_Q}{Q(w)^2}.
\]
In the quotient, set $0/0=0$ and $a/0=+\infty$ for $a>0$. Condition \eqref{si:eq:touchingcertificate} forces $Q(e^{i\theta})=0$ on $\Theta(W)$. It holds for $Q(z)=\prod_{\theta\in\Theta(W)}(z-e^{i\theta})$, a real polynomial of degree $|\Theta(W)|$, admissible when $m\ge|\Theta(W)|$, whenever $\Theta(W)$ is finite and the contact is quadratic, $1-\hat h(\theta)\ge c\operatorname{dist}(\theta,\Theta(W))^2$ for some $c>0$; if the contact has order $p$, $1-\hat h(\theta)\ge c\operatorname{dist}(\theta,\Theta(W))^{2p}$, the $p$-th power of this polynomial, of degree $p|\Theta(W)|$, serves instead when $m\ge p|\Theta(W)|$. A finite contact order exists whenever $\Theta(W)$ is finite.
\end{proposition}
\begin{proof}
For $t\ge m+1$ the window that produces $y_{t+1}$ is $r_{t-m},\dots,r_t$ with $r_{t-j}=Wy_{t-j}$, and $y_{t+1}$ is the minimum-norm point of its affine hull, so $\norm{y_{t+1}}\le\norm{\sum_jc_jr_{t-j}}$ for every $c\in\R^{m+1}$ with $\sum_jc_j=1$. Take $c_j=q_jw^j/Q(w)$. Since $w^jy_{t-j}=w^t\tilde y_{t-j}$, we get $\sum_jc_jr_{t-j}=\frac{w^t}{Q(w)}WQ(S^{-1})\tilde y_t=\frac{w^{t+1}}{Q(w)}AQ(S^{-1})\tilde y_t$, and dividing by $w^{t+1}$ gives \eqref{si:eq:windowinequality}.

Now fix $k\ge1$ with $y_k\ne0$ and $L\ge k+m$, and let $\zeta=\zeta_{k,L}$. The Fourier coefficient of $e^{is\theta}$ in $Q(e^{i\theta})\zeta(\theta)$ equals $Q(S^{-1})\tilde y_s$ for $k+m\le s\le L$, since all indices $s-j$ then lie in $[k,L]$; dropping the other coefficients, Parseval's identity, $H_\theta\preceq\hat h(\theta)I$ and \eqref{si:eq:tailenergy} give
\[
 \begin{aligned}
 \sum_{s=k+m}^{L}\norm{Q(S^{-1})\tilde y_s}^2
 &\le\frac1{2\pi}\int|Q(e^{i\theta})|^2\norm{\zeta}^2d\theta
 \le C_Q\frac1{2\pi}\int\bigl(1-\hat h(\theta)\bigr)\norm\zeta^2d\theta\\
 &\le C_Q\frac1{2\pi}\int\zeta^*(I-H_\theta)\zeta\,d\theta=C_Qn_k^2 .
 \end{aligned}
\]
By \eqref{si:eq:windowinequality} and $\norm A\le2$, $\sum_{s=k+m+1}^{L+1}n_s^2\le\frac{4}{Q(w)^2}\sum_{s=k+m}^L\norm{Q(S^{-1})\tilde y_s}^2\le\frac{4C_Q}{Q(w)^2}n_k^2$, and $\sum_{s=k}^{k+m}n_s^2\le\frac{4^{m+1}-1}3n_k^2$ because $n_{s+1}\le2n_s$. Letting $L\to\infty$ gives \eqref{main:eq:uniformtail} with the stated $C$, and the consequences follow as in the proof of \cref{main:thm:uniformceiling}. If $C_Q<\infty$ then $Q$ vanishes where $1-\hat h$ does, that is, on $\Theta(W)$. Conversely, the stated $Q$ satisfies $|Q(e^{i\theta})|\le2^{|\Theta(W)|-1}\operatorname{dist}(\theta,\Theta(W))$, because $|e^{i\theta}-e^{i\theta'}|\le|\theta-\theta'|$ for the distance on the circle and every factor is at most two; hence $|Q(e^{i\theta})|^2\le C'\operatorname{dist}(\theta,\Theta(W))^2\le(C'/c)(1-\hat h(\theta))$ under quadratic contact, and $|Q(e^{i\theta})|^{2p}\le(C'^p/c)(1-\hat h(\theta))$ under contact of order $p$. Finally, the eigenvalues of the real-analytic Hermitian family $\theta\mapsto H_\theta$ are real-analytic functions of $\theta$ \cite[Chapter~II, \S6]{kato1995}; for a branch $\lambda$ with $\lambda(\theta_0)=1$, $\theta_0\in\Theta(W)$, the nonnegative analytic function $1-\lambda$ either vanishes identically, which would make $\Theta(W)$ the whole circle, or satisfies $1-\lambda(\theta)\ge c(\theta-\theta_0)^{2p}$ near $\theta_0$ for some $c>0$ and integer $p\ge1$. Taking the largest such $p$ over the finitely many branches and touching points gives a finite contact order when $\Theta(W)$ is finite.
\end{proof}

\begin{remark}[The remaining case]\label{si:rem:remaining}
\Cref{main:thm:uniformceiling} answers the question whether the factor is uniformly below $w(W)$ in exactly the cases where \cref{si:prop:spectralmeasure} applies. The question remains open for a $W$ that is not spectraloid whose touching set carries a solution of the moment equations: for $m\ge2$ this requires at least three touching points, and for $m=1$ it includes $\Theta(W)=\{0,\pi\}$ and $\Theta(W)=\{\pm\theta_0\}$ with $\cos\theta_0=w(W)$. An example of the first kind is $W=B\oplus(-B)$ with $B=\bigl[\begin{smallmatrix}0.36&0.72\\0&0.36\end{smallmatrix}\bigr]$: $\norm W=0.36(1+\sqrt2)<1$, and the numerical range of $B$ is the disk of radius $0.36$ centered at $0.36$, so $w(W)=0.72$ is attained at the angles $0$ and $\pi$, by $(1,1,0,0)^T/\sqrt2$ and $(0,0,1,1)^T/\sqrt2$, which are not eigenvectors, whereas $\rho(W)=0.36$. There \cref{si:thm:selfsimilarfpi} still excludes scaled-periodic runs at the ceiling, but nothing excludes a run with $\norm{r_k}/w(W)^k\not\to0$. The tail-energy and separation argument uses the run only through the orthogonality relations \eqref{si:eq:Rjtail}, and in the remaining case these relations admit exact sequences at the ceiling. For $W=B\oplus(-B)$ and $m=1$, the sequence $\tilde y_t=z_0+(-1)^tz_\pi$ with $z_0=\sqrt{(1+w)/2}\,(1,1,0,0)^T/\sqrt2$ and $z_\pi=\sqrt{(1-w)/2}\,(0,0,1,1)^T/\sqrt2$ satisfies \eqref{si:eq:Rjtail} for $j=0,1$ with $n_t\equiv1$, and it fails to be a run only because $\tilde y_{t+1}$ is not on the line through $A\tilde y_t$ and $w^{-1}A\tilde y_{t-1}$: the attaining vector $z_0$ is not an eigenvector, so $Az_0$ is not parallel to $z_0$. An extension of this argument must therefore use additional information, such as membership of $y_{k+1}$ in the affine hull of its window, quantitatively, which Steps 3 and 4 of the proof of \cref{si:thm:selfsimilarfpi} do exactly but not approximately. For $m=1$ and distinct window points, the least-squares projection has the closed form $n_{t+1}^2=\bigl(\norm{p_0}^2\norm{p_1}^2-\ip{p_0}{p_1}^2\bigr)/\norm{p_0-p_1}^2$ with $p_0=A\tilde y_t$ and $p_1=w^{-1}A\tilde y_{t-1}$, the squared distance from the origin to the window line, which may be the right starting point.
\end{remark}

\section{Explicit bounds for special polynomials}\label{si:sec:sispecialized}
This section collects explicit constants. It refines the constant $K_p$
of the energy theorem, finds the best normalized two-term cubic on the
deltoid, shows that certificates of fixed degree certify gains of order
at most $1-q$ on the deltoid, and proves a rate bound for two-term
certificates of any degree, with applications to symmetric matrices and
to the deltoid.

\begin{remark}[Comparison of the bounds]\label{si:rem:compare}
Neither \cref{main:thm:certificate} nor \cref{main:thm:generalmargin} dominates the other. For given $W$ and $p$, \cref{main:thm:generalmargin} gives the smaller block factor exactly when $\tau_p>1-c^2$. For a fixed polynomial, with $q$ bounded away from zero and margin $1-\eta$ tending to zero, the exact formula in \eqref{main:eq:constants} gives $1-c^2=O((1-\eta)^2)$, whereas \eqref{main:eq:generalnormalizedmargin} bounds $\tau_p$ below by a positive multiple of $1-\eta$ when $q$ stays bounded away from zero; for small margins, the regime of \cref{main:sec:limits}, \cref{main:thm:generalmargin} is therefore the sharper statement. For large margins it can be the weaker one, since $\tau_p\le1/(2(d+1))$ always. For $W=I/2$ and the certificate $p(z)=2z-1$, which annihilates $W$, the block factors are $c=2/\sqrt5\approx0.894$ from \cref{main:thm:certificate}, $(1-\tau_p)^{1/2}\approx0.991$ from \cref{main:thm:generalmargin} with $\tau_p=1/(55+2\sqrt2)$, and $(2/3)^{1/2}\approx0.816$ from the simpler bound in \cref{si:thm:twoterm}. For $p(z)=1.01z-0.01$ and the same $W$ the margin is small, and the first two reverse ($0.99995$ against $0.99956$), while its simpler two-term bound gives $0.99751$. AA(1) terminates after two steps in both examples; the numbers compare the bounds, not the iterations. In applications one takes the best available bound.
\end{remark}

\subsection{Refinements of the energy constants}\label{si:sec:energyrefinements}
For ordinary affine AA the constant $K_p$ of \cref{main:thm:generalmargin} can
be reduced. The second refinement is used in
\cref{si:cor:quadraticlinearmargin}.

\begin{remark}[A one-sided refinement for ordinary AA]
\label{si:rem:onesidedmoments}
For ordinary affine AA, both energy factors remain valid with the smaller
constant
\begin{equation}\label{si:eq:generalKonesided}
 K_p=(d+1)\left(\sum_{i=0}^d|\hat a_i|\right)^2
              +\sum_{h=1}^d\max\{C_h,0\}M_h.
\end{equation}
Indeed, polarization and the bound on $\norm{h_j}$ give
\[
 \mu_j-q^j
 \le\sqrt{1-q^{2j}}\norm{\psi_j}
 -\frac{q^j}{2}\left[1-\norm{t_j}^2+
       \left(\sqrt{1-\norm{v_j}^2}-\norm{\psi_j}\right)^2\right]
 \le\sqrt{(j+1)(1-q^{2j})L}.
\]
Thus negative lag sums incur no linear moment-error penalty. The lower
bound $K_p\ge(d+1)(\sum_i|\hat a_i|)^2$ is unchanged, so the range and comparison
proofs still apply. The capped estimates of \cref{main:sec:nonlinear} use the symmetric
constant \eqref{si:eq:generalK}.
\end{remark}

\begin{remark}[A sharper Gram-matrix penalty]\label{si:rem:sharpgram}
For ordinary affine AA, the energy factors remain valid with
\begin{equation}\label{si:eq:generalKsharp}
 K_p^\sharp=\max\left\{2(q^{2d}+E_p),
       \left(\sum_{i=1}^d |\hat a_i|\sqrt{i+1}\right)^2
       +\sum_{h=1}^d\max\{C_h,0\}M_h\right\}.
\end{equation}
This constant is no larger than the one-sided constant in
\eqref{si:eq:generalKonesided}. To see this, first observe that
\(\langle t_0,t_j\rangle=\mu_j\), so Gram errors with an index zero
vanish. For \(1\le i\le j\le d\), defect telescoping and
\eqref{si:eq:generalrawenergy} give
\[
 \bigl|\langle t_i,t_j\rangle-\mu_{j-i}\bigr|
 \le \sqrt{(1-\|t_i\|^2)(1-\|t_j\|^2)}
 \le\sqrt{(i+1)(j+1)}L.
\]
Summing these errors in the polynomial quadratic form replaces the
old Gram penalty by \((\sum_{i=1}^d|\hat a_i|\sqrt{i+1})^2L\).
The one-sided moment estimates are unchanged. Enlarging this constant
by the displayed maximum preserves the energy inequality, while
\(\Delta_p\le q^{2d}+E_p\) gives
\(0<\widehat\tau_p\le\Delta_p/K_p^\sharp\le1/2\).
Moreover, both terms inside the maximum are bounded by the previous
one-sided constant: use
\(q^{2d}+E_p\le(\sum_{i=0}^d|\hat a_i|)^2\), \(d+1\ge2\), and
\(\sum_{i=1}^d|\hat a_i|\sqrt{i+1}\le\sqrt{d+1}\sum_{i=0}^d|\hat a_i|\).
Thus the sharper constant preserves positivity, strict block contraction, and the comparison $\widehat\tau_p\ge\tau_p$; the corresponding upper bounds are $\widehat\tau_p\le1/2$ and $\tau_p\le1/4$.
\end{remark}

\subsection{The best two-term cubic polynomial}

The polynomial family of Cowal, Marshall, and Pollock~\cite{cowal2026deltoid}
starts with $P_0=1$, $P_1=z$, $P_2=z^2$, and obeys
$P_{j+1}(z)=\bigl(3zP_j(z)-P_{j-2}(z)\bigr)/2$ for $j\ge2$.
In particular, $P_3(z)=(3z^3-1)/2$, and \cite[Theorem~2]{cowal2026deltoid} bounds all these polynomials by one on $\Gamma_3$. Balancing the values of $z^3$ at the cusps and the arc midpoints gives the following bound.

\begin{lemma}[Sharp cubic deltoid bound]\label{si:lem:cubic}
For every $z\in\Gamma_3$, $|27z^3-13|\le14$. On the boundary, equality holds
exactly at the cusps, where $z^3=1$, and at the arc midpoints
$\gamma_3(\pi/3),\gamma_3(\pi),\gamma_3(5\pi/3)$, where $z^3=-1/27$.
\end{lemma}
\begin{proof}
Expanding the cube gives $27\gamma_3(\theta)^3-13=8e^{3i\theta}-1+6e^{-3i\theta}+e^{-6i\theta}$,
whose real and imaginary parts are $2\cos^23\theta+14\cos3\theta-2$ and
$2\sin3\theta\,(1-\cos3\theta)$. Therefore
\begin{align*}
 \tfrac14|27\gamma_3(\theta)^3-13|^2-49
 &=(\cos^23\theta+7\cos3\theta-1)^2+\sin^23\theta\,(1-\cos3\theta)^2-49\\
 &=(\cos^23\theta-1)(16\cos3\theta+47)\le0.
\end{align*}
Equality requires $\cos3\theta=\pm1$, which gives the stated points. The
maximum-modulus principle extends the bound to $\Gamma_3$.
\end{proof}

The bound $|3z^3-1|\le2$ follows from
$|3z^3-1|\le(|27z^3-13|+4)/9\le2$.
The bound on $27z^3-13$ yields the best normalized two-term cubic below.

\begin{theorem}[Explicit depth-three factor]\label{si:thm:cubic}
Let $W$ be real and normal, let $q=\norm{W}=\rho(W)\in(0,1)$, and
assume the full spectral enclosure $\spec(W)\subseteq q\Gamma_3$. Set $p_3(z)=(27z^3-13q^3)/(27-13q^3)$, $a_q=14/(27-13q^3)$, $\beta_q=27\sqrt3/(27-13q^3)$, and $c_q=(a_q+\beta_q\sqrt{1+\beta_q^2-a_q^2})/(1+\beta_q^2)$.
Then $p_3$ is the unique minimizer in
\begin{equation}\label{si:eq:cubicminimax}
 \inf_{\substack{a_0,a_3\in\C\\a_0+a_3=1}}
 \max_{z\in q\Gamma_3}|a_0+a_3z^3|
 =\frac{14q^3}{27-13q^3}.
\end{equation}
and for every $m\ge3$, $c_q<1$ and $\norm{r_{k+4}}\le c_q q^4\norm{r_k}$, $\rhAA(x_0)\le q c_q^{1/4}<q$.
\end{theorem}
\begin{proof}
$p_3(1)=1$, and normality with \cref{si:lem:cubic} gives $\norm{p_3(W)}=\max_{\lambda\in\spec(W)}|p_3(\lambda)| \le14q^3/(27-13q^3)<q^3$,
the last inequality because $13q^3<13$. The only nonconstant coefficient of
$p_3$ is $27/(27-13q^3)$ at degree three, so the number $\beta$ in \eqref{main:eq:constants} is $\sqrt3\cdot27/(27-13q^3)$.
Apply \cref{main:thm:certificate} with $d=3$.

For any $p(z)=a_0+a_3z^3$ with $p(1)=1$, interpolation at the boundary
points $q$ and $-q/3$ gives $28q^3=(27+q^3)p(q)-27(1-q^3)p(-q/3)$.
The triangle inequality yields the lower bound in
\eqref{si:eq:cubicminimax}, and the cubic lemma gives attainment by $p_3$.
Equality forces $p(q)=14q^3/(27-13q^3)$ and
$p(-q/3)=-14q^3/(27-13q^3)$; these two values determine $p$ uniquely.

\end{proof}

A simpler but weaker expression follows from \cref{main:cor:explicitconstant}:
since $1-a_q=13(1-q^3)/(27-13q^3)$ and $1+\beta_q^2=[(27-13q^3)^2+2187]/(27-13q^3)^2$,
\[
 \rhAA(x_0)\le q\left[1-
 \frac{169\,(1-q^3)^2}{8\bigl((27-13q^3)^2+2187\bigr)}\right]<q.
\]
This proves strict improvement, though not the square-root speedup
that the exterior growth of these polynomials yields for the power
method in~\cite{cowal2026deltoid}.

\subsection{Limits of fixed-degree full-deltoid certificates}
Whatever its coefficients, a certificate of fixed degree certifies a gain
of order at most $1-q$ on the deltoid. These fixed-degree estimates
do not settle the square-root question of \cref{main:sec:sqrtstatus}.

\begin{proposition}[A fixed-degree limit on polynomial margins]
\label{si:prop:fixeddegreecertificate}
Let $d\ge1$, $0<q<1$, and
\[
 A_d(q)=\inf_{\substack{p\in\C[z],\ \deg p\le d\\p(1)=1}}
               \max_{z\in q\Gamma_3}|p(z)|.
\]
If $T_d$ denotes the degree-$d$ Chebyshev polynomial, then $A_d(q)\ge 1/(T_d((3-q)/(2q)))$, $1-A_d(q)/q^d\le(3d^2/(2q)-d)(1-q)$.
For every fixed $d\ge2$, the maximal polynomial margin
$1-A_d(q)/q^d$ is $\Theta_d(1-q)$ as $q\uparrow1$, without any
assumption on the coefficients.
If a polynomial with $p(1)=1$ and degree at most $d=d(q)$ satisfies
$\max_{q\Gamma_3}|p|\le q^d(1-c\sqrt{1-q})$ for some $c>0$, then $d\ge(q/3)\bigl(1+\sqrt{1+6c/(q\sqrt{1-q})}\bigr)\ge\sqrt{2cq/3}\,(1-q)^{-1/4}$.
\end{proposition}
\begin{proof}
The real segment $[-q/3,q]$ lies in $q\Gamma_3$. Put
$\tilde p(t)=p(q(2t+1)/3)$ and $t_q=(3-q)/(2q)>1$; then $\tilde p(t_q)=1$.
Let $L_j$ be the Lagrange cardinal polynomials at
$t_j=\cos(j\pi/d)$, $0\le j\le d$.
For $t>1$, the sign of $L_j(t)$ is $(-1)^j$. Interpolating $T_d$,
whose value at $t_j$ is $(-1)^j$, therefore gives
$\sum_j|L_j(t)|=T_d(t)$. Interpolating $\tilde p$ and applying the triangle
inequality gives $1\le T_d(t_q)\max_{-1\le t\le1}|\tilde p(t)| \le T_d(t_q)\max_{z\in q\Gamma_3}|p(z)|$.
This argument applies to complex coefficients and proves the first bound.

Since $T_d(t_q)=\cosh(d\operatorname{arcosh}t_q)$, $\log\cosh x\le x^2/2$
for $x\ge0$, and $(\operatorname{arcosh}t_q)^2\le2(t_q-1)=3(1-q)/q$, we have $\log(q^dT_d(t_q))\le d\log q+3d^2(1-q)/(2q)\le(3d^2/(2q)-d)(1-q)$.
The inequality $1-e^{-u}\le u$ proves the second bound.
For the matching lower order when $d\ge2$, use
$p(z)=z^{d-2}p_{2,t_{\mathrm{op}}}(z)$. The cusp contacts in
\cref{si:prop:quadraticminimax} have modulus $q$, so
$\max_{q\Gamma_3}|p|=q^d a_{2,\mathrm{op}}$. Thus
$1-A_d(q)/q^d\ge (2/\sqrt3-1)(1-q)+O((1-q)^2)$.
The asserted degree requirement follows by solving
$3d^2/(2q)-d\ge c/\sqrt{1-q}$.
\end{proof}

\begin{remark}[A limit of the fixed-degree energy guarantee]
\label{si:rem:fixeddegreeenergylimit}
Fix $d$, and suppose the energy estimate is applied with $\norm{p(W)}$
replaced by the full-deltoid bound $\max_{q\Gamma_3}|p|\le q^d$, with $\Delta_p>0$.
Its certified improvement over $q$ is $O_d(1-q)$, uniformly over the
real normalized polynomials. Indeed, the disk $|z|\le q/3$ lies in
$q\Gamma_3$, and Cauchy's coefficient estimate gives
$|a_j|\le(3/q)^j\max_{q\Gamma_3}|p|\le3^j$. Consequently $E_p=(1-q^2)\sum_{\ell=0}^{d-1}q^{2\ell}\bigl(\sum_{j=0}^{\ell}a_j\bigr)^2\le2(1-q)\sum_{\ell=0}^{d-1}\bigl(\sum_{j=0}^{\ell}3^j\bigr)^2$.
Also \cref{si:prop:fixeddegreecertificate} gives
$q^{2d}-(\max_{q\Gamma_3}|p|)^2=O_d(1-q)$, uniformly for $q\ge1/2$.
Thus $\Delta_p=q^{2d}-(\max_{q\Gamma_3}|p|)^2+E_p=O_d(1-q)$.
The constant in \eqref{si:eq:generalKsharp} has $K_p^\sharp\ge2q^{2d}$, so
$\widehat\tau_p\le\Delta_p/K_p^\sharp=O_d(1-q)$.
Finally, $0\le q-q(1-\widehat\tau_p)^{1/(2(d+1))} \le\widehat\tau_p=O_d(1-q)$.
This is a limitation of the stated fixed-degree energy guarantee. It
places no lower bound on the actual AA factor, and allows faster
certificates that use a smaller matrix-specific spectral set.
\end{remark}

\subsection{Two-term certificates at every depth}\label{si:sec:linearmargin}
For a two-term polynomial, the last projection separates the certificate
into an orthogonal component and a squared norm. Keeping this square gives
a linear certificate margin and a sharper scalar bound without the
coefficient constants needed for general polynomials.

\begin{theorem}[Two-term certificate]\label{si:thm:twoterm}
Let $W\in\R^{n\times n}$ with $q=\norm W\in(0,1)$, and let $m\ge d\ge1$.
Suppose $p(z)=a_0+a_dz^d$, $a_0+a_d=1$, $\norm{p(W)}<q^d$,
where the coefficients are real. These hypotheses force $a_d>0$.
Let $\tau\in(0,1)$ be the unique solution of
\begin{equation}\label{si:eq:twotermtau}
 (1+da_d)\tau-
 \left[\frac{|a_0|}{q^d}\sqrt{1-q^{2d}}
            -a_d\sqrt{(d+1)\tau}\right]_+^2
 =1-\frac{\norm{p(W)}^2}{q^{2d}},
\end{equation}
where $[t]_+=\max\{t,0\}$. Then exact unconstrained AA($m$) satisfies
\begin{equation}\label{si:eq:twotermblock}
 \norm{r_{k+d+1}}^2\le(1-\tau)q^{2(d+1)}\norm{r_k}^2
 \qquad(k\ge0).
\end{equation}
Consequently,
\begin{equation}\label{si:eq:twotermrate}
 \rhAA(x_0)\le q(1-\tau)^{1/[2(d+1)]}
 \le q\left(1-\frac{\tau}{2(d+1)}\right)<q.
\end{equation}
The simpler explicit lower bound
\begin{equation}\label{si:eq:twotermlinearmargin}
 \tau\ge\frac{1-\norm{p(W)}^2/q^{2d}}{1+da_d}
       \ge\frac{1-\norm{p(W)}/q^d}{1+da_d}
\end{equation}
can be used in both bounds.
\end{theorem}
This explicit bound is the one quoted at the end of \cref{main:sec:lineargain}.
\begin{proof}
If $a_d\le0$, then $a_0=1+|a_d|$ and $\norm{p(W)}\ge a_0-|a_d|\norm{W^d} \ge1+|a_d|(1-q^d)\ge1>q^d$,
which contradicts the certificate. Thus $a_d>0$.

Fix $k$. The case $r_k=0$ is immediate; otherwise normalize so that
$\norm{r_k}=1$. Use the rescaled residuals and corrections from the proof
of \cref{main:thm:generalmargin}: $T=W/q$, $u_j=q^{-j}r_{k+j}$, $v_j=q^{-j}y_{k+j+1}$, $f_j=u_j-v_j$.
They satisfy $u_{j+1}=Tv_j$ and $f_j\perp v_j$. Put
$L=1-\norm{u_{d+1}}^2$. The energy identity
\eqref{si:eq:generalenergy} gives $L\in[0,1]$ and
$1-\norm{v_d}^2\le L$. The joint raw-power estimate
\eqref{si:eq:generalrawenergy} gives $1-\norm{T^du_0}^2+\norm{T^du_0-v_d}^2 \le(d+1)(1-\norm{v_d}^2)$.
Set $w=v_d$ and $\psi=T^du_0-w$. Polarization therefore yields
$2\ip w\psi\ge-d(1-\norm w^2)\ge-dL$ and $\norm\psi^2\le(d+1)L$.

The last projection retains $r_k=u_0$. Thus
$\ip w{u_0}=q^d\norm w^2$, and, with $h=u_0-q^dw$, $h\perp w$, $\norm h^2=1-q^{2d}\norm w^2\ge1-q^{2d}$.
Since $a_0+a_d=1$,
$p(W)u_0=q^dw+a_0h+a_dq^d\psi$. Keeping the squared norm gives
\[
 \begin{split}
 \norm{p(W)u_0}^2
 &=q^{2d}\norm w^2+\norm{a_0h+a_dq^d\psi}^2
                         +2a_dq^{2d}\ip w\psi\\
 &\ge q^{2d}\bigl[1-(1+da_d)L\bigr]
   +\left[|a_0|\sqrt{1-q^{2d}}
                  -a_dq^d\sqrt{(d+1)L}\right]_+^2.
 \end{split}
\]
The certificate bound therefore says that the left side of
\eqref{si:eq:twotermtau}, evaluated at $L$, is at least its right side.
This scalar left side is continuous and strictly increasing on $[0,1]$:
its first term strictly increases and its subtracted squared positive
part decreases. At zero it is nonpositive, while the right side is
positive. At one, use $\bigl[|a_0|\sqrt{1-q^{2d}}/q^d-a_d\sqrt{d+1}\bigr]_+\le[|a_0|/q^d-a_d]_+\le\norm{p(W)}/q^d$,
where the last inequality follows from
$\norm{p(W)}\ge|a_0|-a_dq^d$.
Its left side at one consequently exceeds the right side by at least
$da_d>0$. This proves existence and uniqueness of $\tau\in(0,1)$ and
$L\ge\tau$. Restoring the scale proves \eqref{si:eq:twotermblock}; block
iteration and concavity prove \eqref{si:eq:twotermrate}.
Finally, dropping the nonnegative square in \eqref{si:eq:twotermtau}
gives \eqref{si:eq:twotermlinearmargin}.
\end{proof}

\begin{corollary}[Symmetric spectra in an asymmetric interval]\label{si:cor:lineartwoeigen}
Let $W$ be real symmetric with $\norm W=q\in(0,1)$ and
$\spec(W)\subseteq[-q',q]$, where $0\le q'<q$. Set $s=(q-q')/2$.
Then exact AA($m$), for every $m\ge1$ and every initialization, satisfies
\[
 \rhAA(x_0)
 \le q\left[1-
 \frac{s(1-q)[2q-s(1+q)]}{q^2(1-s)(2-s)}\right]^{1/4}
 \le q-\left[1-(5/8)^{1/4}\right](1-q)(q-q').
\]
\end{corollary}
\begin{proof}
Take $p(z)=(z-s)/(1-s)$. Then $a_0=-s/(1-s)$, $a_1=1/(1-s)$,
and $\norm{p(W)}\le(q-s)/(1-s)<q$. The bound
\eqref{si:eq:twotermlinearmargin} decreases as $\norm{p(W)}$ increases, and
substituting $(q-s)/(1-s)$ for $\norm{p(W)}$ gives $\tau\ge s(1-q)[2q-s(1+q)]/(q^2(1-s)(2-s))$,
which proves the first rate bound. Set
$z=(1-q)(q-q')/q=2s(1-q)/q\in(0,1)$. Since $s\le q/2$, $s(1-q)[2q-s(1+q)]/(q^2(1-s)(2-s)) \ge z/2-z^2/8$.
After cancellation of positive factors, this inequality is equivalent to $2q-s(1+q)-[q-s(1-q)/2](1-s)(2-s)=s\bigl[2q-3s+qs+s^2(1-q)\bigr]/2\ge0$.
The function $(1-z/2+z^2/8)^{1/4}$ is convex on $[0,1]$, since its
second derivative is
$[4+4z-z^2]/[256(1-z/2+z^2/8)^{7/4}]>0$.
Its endpoint secant therefore gives $\rhAA(x_0)\le q(1-z/2+z^2/8)^{1/4} \le q[1-(1-(5/8)^{1/4})z]$,
which proves the second bound. The coefficient
$1-(5/8)^{1/4}$ is the largest uniform coefficient supplied by the first
rate bound: take $q'=0$ and $q\downarrow0$. This optimality statement
concerns the explicit bound.
\end{proof}

\begin{corollary}[Linear-margin cubic deltoid rate]\label{si:cor:cubiclinearmargin}
Let $W$ be real normal, $q=\norm W\in(0,1)$, and
$\spec(W)\subseteq q\Gamma_3$, where $\Gamma_3$ is the deltoid of \eqref{main:eq:hypocycloid}. Define $\tau_3(q)= 13(1-q^3)(41-13q^3)/((27-13q^3)(108-13q^3))$.
For every $m\ge3$, exact unconstrained AA($m$) satisfies $\norm{r_{k+4}}^2\le(1-\tau_3(q))q^8\norm{r_k}^2$, $\rhAA(x_0)\le q(1-\tau_3(q))^{1/8}<q$.
The same conclusion holds for any real matrix $W$ with $\norm W=q<1$
satisfying $\norm{27W^3-13q^3I}\le14q^3$.
\end{corollary}
\begin{proof}
\Cref{si:lem:cubic,si:thm:cubic} give the certificate $a_0=-13q^3/(27-13q^3)$, $a_3=27/(27-13q^3)$, $\norm{p(W)}\le14q^3/(27-13q^3)<q^3$.
Substituting this bound for $\norm{p(W)}$ in
\eqref{si:eq:twotermlinearmargin} gives $(1-196/(27-13q^3)^2)/(1+3a_3) =13(1-q^3)(41-13q^3)/((27-13q^3)(108-13q^3))$.
Apply \cref{si:thm:twoterm}. Normality is used only to obtain the
operator certificate from spectral containment.
\end{proof}

The gain $q-q(1-\tau_3(q))^{1/8}$ is of linear order in $1-q$:
its ratio to $1-q$ tends to $273/2660$ as $q\uparrow1$.
The three-term quadratic is covered by \cref{main:thm:generalmargin}.

\section{Further results for normal deltoid matrices}\label{si:sec:sinormal}
This section states and proves \cref{si:thm:deltoidactivelinear}, gives the normal counterexamples mentioned after \cref{main:thm:everymemoryquadratic}, treats polynomial approximation on the deltoid and the hypocycloids, and proves \cref{main:thm:contact}.

\subsection{A quantitative improvement over the active fixed-point factor}
\label{si:sec:deltoidpicard}
The proof of \cref{si:thm:deltoidactivelinear} uses the radial structure of the deltoid and the restriction of AA to the spectral subspaces excited by the initial point.

\begin{lemma}[Radial structure of the deltoid]\label{si:lem:deltoidangular}
Put $\Psi(z)=3|z|^4+6|z|^2-8\operatorname{Re}(z^3)-1$ and $C(r)=(3r^4+6r^2-1)/(8r^3)$ for $r>0$. Then $\Psi=0$ on the boundary of $\Gamma_3$, $\Gamma_3=\{z:\Psi(z)\le0\}$, and $\Gamma_3$ is star-shaped with respect to the origin; in particular $q\Gamma_3\subseteq\Gamma_3$ for $0<q\le1$. Every point $re^{i\theta}\in\Gamma_3$ with $r>0$ satisfies
\begin{equation}\label{si:eq:deltoidangularrestriction}
 \cos3\theta\ge C(r),\qquad C'(r)=\frac38(1-r^{-2})^2\ge0,\qquad
 1-C(r)=\frac{(1-r)^3(1+3r)}{8r^3}.
\end{equation}
\end{lemma}
\begin{proof}
On the boundary, $|\gamma_3(\theta)|^2=(5+4\cos3\theta)/9$ and $27\operatorname{Re}\gamma_3(\theta)^3=2\cos^23\theta+14\cos3\theta+11$, and substitution gives $\Psi(\gamma_3(\theta))=0$. On a ray, $\Psi(re^{i\theta})=3r^4+6r^2-8r^3\cos3\theta-1$, $\partial_r\Psi(re^{i\theta})=12r\bigl[(r-1)^2+2r(1-\cos3\theta)\bigr]\ge0$,
and for $r>0$ the derivative vanishes only at the unit cusps. Thus $\Psi$ is strictly increasing along every ray, with $\Psi(0)=-1$, and vanishes at most once on it. Since $0$ is an interior point of the bounded region $\Gamma_3$, every ray leaves $\Gamma_3$ and therefore meets its boundary exactly once, at the zero of $\Psi$, and the segment before that point lies in $\Gamma_3$. This proves $\Gamma_3=\{\Psi\le0\}$ and star-shapedness. The inequality $\Psi(re^{i\theta})\le0$ reads $\cos3\theta\ge C(r)$, and the formulas for $C'$ and $1-C$ follow by differentiation and expansion.
\end{proof}

\begin{lemma}[The active subspace]\label{si:lem:activesubspace}
Let $W$ be real and normal with $\norm W<1$, and let $\mathcal S$ be the orthogonal sum of the real spectral subspaces of $W$ on which $x_0-x_*$ has a nonzero component. Then $\mathcal S$ reduces $W$ and contains every AA($m$) residual, these residuals are those of AA($m$) applied to the normal restriction $W_{\mathcal S}$, and $\norm{W_{\mathcal S}}=\rho(W_{\mathcal S})=\rho_{\mathrm{Picard}}(x_0)$.
\end{lemma}
\begin{proof}
Since $W-I$ is invertible and preserves the spectral subspaces, $r_0=(W-I)(x_0-x_*)$ has nonzero components on exactly the same subspaces, so $r_0\in\mathcal S$, and \cref{si:prop:invariantsubspace} gives the first two claims. The restriction is normal, so its norm is its spectral radius, the largest modulus of an eigenvalue whose spectral subspace is excited by $x_0-x_*$; by orthogonality of the spectral subspaces this is $\rho_{\mathrm{Picard}}(x_0)$.
\end{proof}

\begin{theorem}[A quantitative improvement over the active fixed-point factor]
\label{si:thm:deltoidactivelinear}
There exists $\rho_0<1$ with the following property. Let $W$ be real
and normal, $q=\norm W<1$, and $\spec(W)\subseteq q\Gamma_3$.
For every memory $m\ge2$ and every initial point with
$\rho:=\rho_{\mathrm{Picard}}(x_0)\in(\rho_0,1)$, AA($m$) satisfies
\begin{equation}\label{si:eq:deltoidactivelinear}
 \rhAA(x_0)\le\rho-\frac{\alpha_*}{12}(1-\rho),
\end{equation}
where the positive constant $\alpha_*\approx1.27\times10^{-2}$, given in
\eqref{si:eq:deltoidactivelinearconstant}, is independent of the matrix,
dimension, memory, and initialization. More precisely, a common certified
factor has the expansion
$R(\rho)=\rho-\alpha_*(1-\rho)/6+O((1-\rho)^2)$ as $\rho\uparrow1$.
\end{theorem}

\begin{proof}[Proof of \cref{si:thm:deltoidactivelinear}]
We use the constants
\begin{equation}\label{si:eq:deltoidactivelinearconstant}
 K_*=\frac{324+54\sqrt6+42\sqrt2}{100},\qquad
 G_*=\frac{21+18\sqrt3}{25},\qquad
 \alpha_*=
 \frac{9/25}{\left(\sqrt{G_*^2+(6/5)K_*}+G_*\right)^2}.
\end{equation}

\emph{A fixed quadratic with one boundary contact.}
Set $p(z)=(9z^2+3z-2)/10$, $K_\rho=\Gamma_3\cap\{|z|\le\rho\}$, $A_\rho=\max_{z\in K_\rho}|p(z)|$.
For $z=\gamma_3(t)$ and $x=\cos t$, direct expansion gives $100\bigl(1-|p(z)|^2\bigr)=(1-x)Q(x)$, $Q(x)=60+172x+64x^2-32x^3+320x^4+256x^5$.
Here $Q(x)\ge28/5$ on $[-1,1]$. For an exact positivity
certificate, its coefficients in the degree-five Bernstein basis
after linear rescaling of each indicated interval to $[0,1]$ are
\[
\begin{array}{c|rrrrrr}
 {[-1,-1/2]}&48&214/5&128/5&12&28/5&6\\
 {[-1/2,0]}&6&32/5&68/5&136/5&214/5&60\\
 {[0,1]}&60&472/5&676/5&896/5&1436/5&840
\end{array}
\]
and each Bernstein basis function is nonnegative, with their sum equal
to one. It follows that $|p|\le1$ on the boundary of $\Gamma_3$, with
equality only at $z=1$. The maximum-modulus principle gives the same
conclusion on the filled region. Also $p(1)=1$ and $p'(1)=21/10$.

\emph{The maximum on the truncated deltoid.}
Put $\varepsilon=1-\rho$. We claim
\begin{equation}\label{si:eq:activequadraticmaximum}
 A_\rho=p(\rho)+O(\varepsilon^3),\qquad \rho\uparrow1.
\end{equation}
The lower bound $A_\rho\ge p(\rho)$ holds for $\rho$ near one,
because $\rho\in K_\rho$ and $p(\rho)>0$.
Every maximizing point $z_\rho\in K_\rho$ tends to $1$:
otherwise compactness, continuity, and $p(\rho)\to1$ would produce
another modulus-one contact of $p$ on $\Gamma_3$.
Write $z_\rho=r e^{i\phi}$ with $\phi\to0$ and $s=1-r\ge\varepsilon$.
By \eqref{si:eq:deltoidangularrestriction}, $0\le1-\cos3\phi \le s^3(1+3r)/(8r^3)$, $\phi^2=O(s^3)$.
All constants here and below are uniform in a fixed neighborhood of
$1$. Since $p$ is nonzero there and has real coefficients,
$|p(re^{i\phi})|$ is a smooth even function of $\phi$.
Taylor expansion thus gives
$|p(re^{i\phi})|\le p(r)+C s^3$.
Moreover $p'(r)$ is bounded below by some $c_1>0$ near $r=1$.
Consequently $A_\rho\le p(\rho)-c_1(s-\varepsilon)+Cs^3$.
For $\rho$ sufficiently close to one, $Cs^2\le c_1/2$.
Comparing with $A_\rho\ge p(\rho)$ gives $s\le2\varepsilon$ and
then $A_\rho\le p(\rho)+8C\varepsilon^3$, proving
\eqref{si:eq:activequadraticmaximum}. In particular, $\rho^2-A_\rho =\varepsilon/10+\varepsilon^2/10+O(\varepsilon^3)>0$
for all sufficiently small positive $\varepsilon$.

\emph{Transfer to the active residual subspace.}
Use the active subspace $\mathcal S$ of \cref{si:lem:activesubspace}. Its normal restriction has norm $\rho$ and spectrum in $q\Gamma_3\cap\{|z|\le\rho\}\subseteq K_\rho$ by \cref{si:lem:deltoidangular}, and every AA residual belongs to $\mathcal S$. Thus $\norm{p(W_{\mathcal S})}\le A_\rho<\rho^2$.
For the coefficients $(-1/5,3/10,9/10)$, let $E_p$ be as in \eqref{si:eq:generalgrambaseline}. Using $A_\rho$ as the operator-norm bound gives a certified excess $\Delta_\rho\le\Delta_p$, where $\Delta_p$ is defined in \eqref{main:eq:generalexcess}, with
\begin{align*}
 E_p&=\frac{(1-\rho^2)(4+\rho^2)}{100}
     =\frac{\varepsilon}{10}+O(\varepsilon^2),\\
 \rho^4+E_p-p(\rho)^2
 &=\frac{3\rho}{50}(1-\rho)(2+6\rho-3\rho^2),\\
 \Delta_\rho&:=\rho^4+E_p-A_\rho^2
          =\frac{3\varepsilon}{10}+O(\varepsilon^2).
\end{align*}
Use the sharper constant of \cref{si:rem:sharpgram}; its floor
is inactive near $\rho=1$, because $K_*>2$.
Direct substitution in \eqref{si:eq:generalG} and~\eqref{si:eq:generalKsharp} gives
$K_p=K_*+O(\varepsilon)$ and
$G_p=G_*\sqrt\varepsilon+O(\varepsilon^{3/2})$, with the constants
in \eqref{si:eq:deltoidactivelinearconstant}. Hence
\eqref{si:eq:generalsharptau}, with $\Delta_\rho$ in place of $\Delta_p$, gives a common certified gain $\widehat\tau_\rho=\alpha_*\varepsilon+O(\varepsilon^2)$ and $R(\rho):=\rho(1-\widehat\tau_\rho)^{1/6}=\rho-\alpha_*\varepsilon/6+O(\varepsilon^2)$.
These quantities depend only on $\rho$ and the fixed polynomial.
The corollary applies to every $m\ge2$, including ordinary startup,
so $\rhAA(x_0)\le R(\rho)$. Choosing a common $\rho_0$ sufficiently
close to one proves \eqref{si:eq:deltoidactivelinear}.
\end{proof}

\begin{remark}[The quantitative threshold]\label{si:rem:rhozero}
The proof establishes a universal $\rho_0<1$ without assigning it a
numerical value. Sampling $\partial K_\rho$ only gives a lower bound on
$A_\rho$; substituting that sampled maximum into the rate formula does
not certify an upper bound for AA. An explicit numerical threshold
would require a rigorous upper bound for $A_\rho$ over a whole interval
of $\rho$.
\end{remark}

\subsection{Proof of \texorpdfstring{\cref{main:thm:everymemoryquadratic}}{the fixed-certificate theorem} and further properties of its examples}\label{si:sec:everymemoryfurther}
\begin{proof}[Proof of \cref{main:thm:everymemoryquadratic}]
Put $q=1-\varepsilon$ and $\rho=1-2\varepsilon$, let $U\in\R^{(m+1)\times(m+1)}$ and $v$ be as in \cref{main:lem:momentmatrix} with $\rho$ in place of $q$, let
\[
 N=\begin{pmatrix}0&q\\ 1/(4q)&0\end{pmatrix},
\]
which has norm $q$ because $q>1/2$, square $N^2=I/4$, and eigenvalues $\pm1/2$, and set $W_\varepsilon=\diag(\rho U,N)$ and $x_{0,\varepsilon}=(W_\varepsilon-I)^{-1}r_0$ with $r_0=(v,0,0)^T$, so that $x_*=0$. Since $q>\rho>1/2$, the matrix $W_\varepsilon$ is
invertible and diagonalizable, has norm $q$ and spectral radius $\rho$, and
satisfies $\norm{W_\varepsilon^2}=\rho^2$; it is nonsymmetric because $N$ is.

Suppose inductively that $r_i=(\rho^iU^iv,0,0)^T$ through step $k$. For
every retained lag $0\le j\le m_k$, \eqref{main:eq:moments} gives
$\ip{r_k}{r_{k-j}}=\rho^{2k-j}\ip{U^jv}{v}=\rho^{2k}=\norm{r_k}^2$. Thus $r_k$
is orthogonal to every direction $r_{k-j}-r_k$ of the retained affine hull,
hence is its minimum-norm point, and $r_{k+1}=W_\varepsilon r_k$; this
includes the first step at $k=0$. Every retained set of residuals is the
image under an invertible power of $\rho U$ of a set
$\{v,\rho Uv,\ldots,(\rho U)^{m_k}v\}$, which is linearly independent by
\cref{main:lem:momentmatrix}, so the only minimizing coefficient vector
selects $r_k$. Consequently $x_k=W_\varepsilon^kx_{0,\varepsilon}$, its
auxiliary components vanish, and orthogonality of $U$ gives
$\norm{x_k}=\rho^k\norm{x_{0,\varepsilon}}$, which proves the last equation of 
\eqref{main:eq:everymemoryfactor}. Finally,
$(1-\rhAA(x_{0,\varepsilon}))/\sqrt{1-\norm{W_\varepsilon}}=2\sqrt\varepsilon\to0$.
\end{proof}

\emph{Capped AA and nonlinear maps.} Since the minimizing coefficients select the newest residual, their $\ell^1$ norm is one, so capped AA($m$) with any cap $\Lambda\ge1$ (\cref{main:sec:geometry}) produces the same run. It also arises near a fixed point of an analytic map: $G(x,t)=(W_\varepsilon x,t^2)$ has derivative $\diag(W_\varepsilon,0)$ at its fixed point $0$ and acts as $W_\varepsilon$ on the invariant hyperplane $t=0$, so the same runs occur from starts arbitrarily close to that fixed point. Neither the cap nor nonlinearity removes the obstruction.

\emph{The example is not an artifact of the Euclidean inner product.} With $S=\diag(I,1/(2q),1)$, the matrix $SW_\varepsilon S^{-1}=\diag\bigl(\rho U,\tfrac12\begin{psmallmatrix}0&1\\1&0\end{psmallmatrix}\bigr)$ is normal with norm $\rho=\rho(W_\varepsilon)$, so the inner product of $H=S^TS$ is adapted to $W_\varepsilon$. The weights carried by $v$ solve the moment equations of degree $m$ at the level $\rho$, so by \cref{si:cor:spectralbenchmark} no inner product with $\norm{W_\varepsilon}_H=\rho(W_\varepsilon)$ makes AA($m$) beat $\rho(W_\varepsilon)$ from every start, whereas AA$_H$($m+1$) does, because $U$ has $m+1$ distinct eigenvalues (\cref{main:prop:peripheralcount}). The certificate $z^2$ reflects only the auxiliary block.

\emph{The factor cannot be bounded away from $q$ uniformly.} Even for fixed $q$, no bound below $q$ holds uniformly over the matrices with norm $q$ for which $z^2$ is a strict certificate. For $0<\rho<q<1$, replacing $N$ by
\[
 \begin{pmatrix}0&q\\ \rho^2/(4q)&0\end{pmatrix},
\]
which has norm $q$, square $(\rho^2/4)I$, and eigenvalues $\pm\rho/2$, gives $\norm W=q$ and $\norm{W^2}=\rho^2<q^2$, so $z^2$ is a certificate, while $\rhAA(x_0)=\rho$; as $\rho\uparrow q$, the factor tends to $q$ and the margin to zero.

\emph{Normality does not restore square-root acceleration from individual starts.} \Cref{si:prop:normalallmemory} gives, at every memory, a normal $W$ with $\norm W=1-\varepsilon$ and a degree-one certificate, together with a start from which AA follows fixed-point iteration with factor $1-2\varepsilon$, and \cref{si:rem:depthoneobstruction} gives a normal example at memory one with a fixed linear certificate.

\subsection{A normal counterexample at every memory}
\label{si:sec:normalobstruction}

\begin{proposition}[A normal obstruction at every memory]
\label{si:prop:normalallmemory}
Fix $m\ge1$ and $0<\varepsilon<1/2$. There are a real normal invertible
matrix $W_\varepsilon\in\R^{(m+2)\times(m+2)}$ and a nonzero initial
point such that, with $q=1-\varepsilon$, $\norm{W_\varepsilon}=\rho(W_\varepsilon)=q$, $\rhAA(x_0)=\rho_{\mathrm{Picard}}(x_0)=1-2\varepsilon$.
Ordinary AA($m$) follows fixed-point iteration exactly, with unique minimizing
coefficients. Nevertheless, $p_\varepsilon(z)=\bigl(1+\varepsilon/4\bigr)z-\varepsilon/4$, $\norm{p_\varepsilon(W_\varepsilon)}=q-\varepsilon^2/4<q$.
\end{proposition}
\begin{proof}
Put $\rho=1-2\varepsilon$ and use the orthogonal matrix $U$ and vector
$v$ of \eqref{main:eq:momentmatrix} and~\eqref{main:eq:moments}; that construction
is valid for every $m\ge1$ and $0<\rho<1$. Set $W_\varepsilon=\diag(\rho U,q)$, $r_0=(v,0)^T$, $x_0=(W_\varepsilon-I)^{-1}r_0$.
The matrix is real, normal and invertible, with norm and spectral
radius $q$. The moment identity $\ip v{U^jv}=\rho^j$ for $0\le j\le m$
gives $\ip{r_k}{r_{k-j}}=\rho^{2k}=\norm{r_k}^2$
for the candidate Picard residuals $r_k=(\rho^kU^kv,0)^T$ and every
retained lag. Hence the latest residual is the minimum-norm point of
each retained affine hull. The independence in \eqref{main:eq:moments}
makes the minimizing coefficients unique. Ordinary startup therefore
produces this Picard sequence, and
$\norm{x_k}=\rho^k\norm{x_0}$.
On the active block, $\norm{p_\varepsilon(\rho U)}\le\rho+\varepsilon(1+\rho)/4=q-\varepsilon/2-\varepsilon^2/2<q-\varepsilon^2/4$.
The remaining scalar block has
$p_\varepsilon(q)=q-\varepsilon^2/4>0$, proving the exact norm and margin.
Finally,
$(1-\rhAA(x_0))/\sqrt{1-\norm{W_\varepsilon}}=2\sqrt\varepsilon\to0$.
\end{proof}

Thus even normal matrices with degree-one certificates need not have a uniform square-root gain over $q$, and need not improve strictly on fixed-point iteration from each individual start. Here the active fixed-point factor is below the spectral radius, the polynomial depends on $\varepsilon$, and its margin is quadratic in $\varepsilon$; the fixed-polynomial, linear-margin assertion of \cref{main:thm:everymemoryquadratic} is stronger in these two respects.

This example is not confined to a scaled deltoid and does not contradict
\cref{si:thm:deltoidactivelinear}. Its start is supported on a subspace
of radius $\rho<q$, so equality with its own Picard rate is compatible
with strict improvement over $q$. It also gives no conclusion about
the measure of slow starting points. At memory one the following
variant retains a fixed polynomial and a linear-order margin.

\begin{remark}[The obstruction at memory one]\label{si:rem:depthoneobstruction}
The same square-root obstruction holds at memory one with a fixed linear
certificate, even for a normal matrix. For $0<\varepsilon<1/4$, put
$q=1-\varepsilon$ and $\rho=1-2\varepsilon$, and take $W=\diag(\rho R_\theta,q)$, $\theta=\arccos\rho$, $r_0=(1,0,0)^T$, $x_0=(W-I)^{-1}r_0$.
Here $R_\theta$ denotes planar rotation. This real normal matrix has
$\norm W=\rho(W)=q$. The fixed polynomial
$p(z)=(3z-1)/2$ satisfies $p(1)=1$ and $\norm{p(W)}^2=\max\{(3\rho^2+1)/4,(3q-1)^2/4\}=1-3\varepsilon+3\varepsilon^2<q^2$.
For the candidate Picard sequence $r_k=W^kr_0$, $\ip{r_k}{r_{k-1}}=\rho^{2k-1}\cos\theta =\rho^{2k}=\norm{r_k}^2$ for $k\ge1$.
Thus the newest residual is the minimum-norm point of each retained
affine line, and ordinary AA(1) follows Picard exactly. Consecutive
residuals are linearly independent, so every minimizing coefficient
vector uniquely selects the newest residual. Both iterations have
factor $\rho$ on the active plane, while the ambient spectral radius is
$q$. In particular, $(1-\rhAA(x_0))/\sqrt{1-\norm W} =2\sqrt\varepsilon\to0$.
\end{remark}

\subsection{Realizing regional polynomial bounds}
The following finite-dimensional minimax argument justifies the lower bound for adaptive polynomial methods in \cref{main:sec:sqrtstatus}.

\begin{lemma}[A finite worst-case realization of the regional bound]
\label{si:lem:regionalworstcase}
Let $E\subseteq\mathbb C$ be compact and invariant under complex conjugation,
with $1\notin E$, and let $n\ge1$. Suppose $E$ contains at least $n+1$
distinct points, and set
\[
 A_n(E)=\min_{\substack{p\in\mathbb R[z],\ \deg p\le n\\p(1)=1}}
             \max_{z\in E}|p(z)|.
\]
There exist a real normal matrix $W$ of dimension at most $2n+2$, with
$\operatorname{spec}(W)\subseteq E$, and a real unit vector $r_0$ such that
\[
 \min_{\substack{p\in\mathbb R[z],\ \deg p\le n\\p(1)=1}}
       \|p(W)r_0\|=A_n(E).
\]
The same minimum is obtained if complex coefficients are allowed.
\end{lemma}
\begin{proof}
Write $p_c(z)=1+\sum_{j=1}^n c_j(z^j-1)$, $c\in\mathbb R^n$.
The maximum norm on $E$ is a norm on polynomials of degree at most $n$,
so a minimizer $c_*$ exists. Put $p_*=p_{c_*}$ and
$F(c)=\max_{z\in E}|p_c(z)|^2$.
At $c_*$ the origin belongs to the convex hull of the gradients
$\nabla_c|p_c(z)|^2|_{c=c_*}$ over the active points
$|p_*(z)|=A_n(E)$: otherwise strict separation gives a direction that
strictly decreases every active function and hence, by compactness,
their maximum for a sufficiently small step. Carath\'eodory's theorem
therefore gives active points $z_1,\ldots,z_s$, $s\le n+1$, and weights
$w_i\ge0$, $\sum_iw_i=1$, such that
\[
 \sum_iw_i\nabla_c|p_c(z_i)|^2\big|_{c=c_*}=0.
\]
The convex quadratic $\sum_iw_i|p_c(z_i)|^2$ consequently has minimum
$\sum_iw_i|p_*(z_i)|^2=A_n(E)^2$.
Replacing a nonreal $z_i$ by its conjugate does not change this quadratic.
For real $z_i$ use a scalar block $z_i$, and for a nonreal $z_i$ use
$B_i=\begin{psmallmatrix}\operatorname{Re}z_i&-\operatorname{Im}z_i\\
\operatorname{Im}z_i&\operatorname{Re}z_i\end{psmallmatrix}$.
Take $W$ to be the orthogonal direct sum of these blocks and let the
component of $r_0$ in block $i$ be $\sqrt{w_i}$ times its first coordinate
vector. Then $\|r_0\|=1$ and
$\|p(W)r_0\|^2=\sum_iw_i|p(z_i)|^2$ for real $p$, proving the assertion.
For complex $p=p_R+ip_I$ with $p(1)=1$, the real polynomial $p_R$ is
admissible and, because $W,r_0$ are real,
$\|p(W)r_0\|^2=\|p_R(W)r_0\|^2+\|p_I(W)r_0\|^2$; thus complex
coefficients cannot lower the minimum.
\end{proof}

\subsection{Polynomial approximation on the scaled deltoid}
Cowal, Marshall, and Pollock~\cite{cowal2026deltoid,cowal2026random} construct, from recurrences associated with random walks, explicit polynomials that are bounded on the deltoid and on the hypocycloids and grow outside them at the rates below; the following statements record the exact constants used in the main text.
\label{si:sec:deltoidapprox}
For $0<q<1$ let $A_d(q)=\inf\{\max_{z\in q\Gamma_3}|p(z)|:\ \deg p\le d,\ p(1)=1\}$,
the best uniform bound on $q\Gamma_3$ that a normalized polynomial of
degree at most $d$ can achieve. Real and complex coefficients give the
same value, because replacing $p$ by $(p+\bar p)/2$ keeps $p(1)=1$
and does not increase the maximum on the conjugation-symmetric set
$q\Gamma_3$. We use the exterior map of the deltoid,
$\psi_q(w)=(q/3)(2w+w^{-2})$, which maps $\{|w|>1\}$ conformally onto
the exterior of $q\Gamma_3$~\cite{cowal2026deltoid}. The upper bound below has
the order of the growth estimate of~\cite[Theorem~2]{cowal2026deltoid},
whose polynomials satisfy the same recurrence as $S_d$.

\begin{proposition}[Polynomial approximation on the scaled deltoid]
\label{si:prop:deltoidapprox}
Let $R>1$ be the unique solution of $q(2R+R^{-2})=3$. Then for
every $d\ge1$, $R^{-d}\le A_d(q)\le\min\{q^d,\ 3R^{-d}\}$,
so $\lim_{d\to\infty}A_d(q)^{1/d}=R^{-1}$, and with $\varepsilon=1-q$, $R=1+\sqrt\varepsilon+2\varepsilon/3+O(\varepsilon^{3/2})$, $R^{-1}=1-\sqrt\varepsilon+\varepsilon/3+O(\varepsilon^{3/2})$.
The upper bound is achieved by $p_d=S_d/S_d(1)$, in the sense that
$\max_{q\Gamma_3}|p_d|\le3R^{-d}$, where
$S_d(z)=\sum_{j=1}^3w_j(z)^d$ and $w_1,w_2,w_3$ are the roots, with
multiplicity, of $w^3-3zw^2/(2q)+1/2=0$. These are real
polynomials with $S_0=3$, $S_1=3z/(2q)$, $S_2=9z^2/(4q^2)$ and $S_j=3zS_{j-1}/(2q)-S_{j-3}/2$ for $j\ge3$,
and $S_j(1)>0$ for all $j$. If $W$ is real normal with
$\spec(W)\subseteq q\Gamma_3$, then $\norm{p(W)}\le\max_{q\Gamma_3}|p|$ for every
polynomial $p$, and in particular $\norm{p_d(W)}\le3R^{-d}$.
\end{proposition}
\begin{proof}
For an admissible $p$, the function $p(\psi_q(w))/w^d$ is analytic for
$|w|>1$, including at infinity, since $\psi_q(w)=2qw/3+O(w^{-2})$.
Its modulus on $|w|=1$ is at most $\max_{q\Gamma_3}|p|$, so by the maximum
principle the same bound holds at $w=R$, where $\psi_q(R)=1$ and the
function equals $R^{-d}$. This is the lower bound.

Writing $\psi_q(w)=z$ as $w^3-3zw^2/(2q)+1/2=0$ identifies the
$w_j(z)$. Newton's identities for the power sums of the roots of a monic
cubic with coefficients $-3z/(2q)$, $0$, $1/2$ give the stated
values of $S_0,S_1,S_2$ and the recurrence, so $S_d$ is a real polynomial
of degree $d$. If $z\in q\Gamma_3$, no root has $|w_j|>1$, because
$\psi_q$ maps $\{|w|>1\}$ onto the exterior of $q\Gamma_3$; hence
$|S_d(z)|\le3$ on $q\Gamma_3$. At $z=1$ the cubic
$f(w)=w^3-3w^2/(2q)+1/2$ has $f(-1)<0<f(0)$, $f(1)<0$ and
$f(w)\to\infty$, so its roots are $R>1$, some $a\in(0,1)$ and some
$-b\in(-1,0)$. The vanishing coefficient of $w$ gives $R(a-b)=ab>0$, so
$a>b$ and $S_d(1)=R^d+a^d+(-b)^d\ge R^d>0$. Thus $p_d=S_d/S_d(1)$ is
admissible with $\max_{q\Gamma_3}|p_d|\le3R^{-d}$, while $p(z)=z^d$ gives
$A_d(q)\le q^d$.

Expanding $q(2R+R^{-2})=3$ about $R=1$ gives
$q\bigl(1+(R-1)^2-4(R-1)^3/3+O((R-1)^4)\bigr)=1$, hence
$(R-1)^2-4(R-1)^3/3=\varepsilon+O(\varepsilon^2)$ and
$R-1=\sqrt\varepsilon+2\varepsilon/3+O(\varepsilon^{3/2})$, which
yields both expansions. For a real normal $W$, $\norm{p(W)}$ is the
largest value of $|p|$ on the spectrum.
\end{proof}

For each finite-dimensional system GMRES terminates, so its literal root
factor is zero; the proposition is a uniform statement over the spectral
region, and for real normal $W$ with $\spec(W)\subseteq q\Gamma_3$ it gives
$\norm{r_d^{\mathrm{GMRES}}}\le\min\{q^d,3R^{-d}\}\norm{r_0}$ for GMRES on
$(I-W)x=b$.

\begin{corollary}[Degree needed for a square-root factor]
\label{si:cor:degreeneeded}
Let $\varepsilon=1-q$ and $c>0$ with $c\sqrt\varepsilon<1$. If
$A_d(q)^{1/d}\le1-c\sqrt\varepsilon$, then $d\ge2qc/(3\sqrt\varepsilon)$.
Conversely, for every $c\in(0,1)$ and $C>\log3/(1-c)$ there is
$\varepsilon_0>0$ such that $d=\lceil C\varepsilon^{-1/2}\rceil$ gives
$A_d(q)^{1/d}\le1-c\sqrt\varepsilon$ whenever $\varepsilon<\varepsilon_0$.
\end{corollary}
\begin{proof}
By \cref{si:prop:fixeddegreecertificate} and its proof,
$A_d(q)\ge1/T_d\bigl((3-q)/(2q)\bigr)\ge\exp\bigl(-3d^2\varepsilon/(2q)\bigr)$,
so $d^{-1}\log A_d(q)\ge-3d\varepsilon/(2q)$.
Since $\log(1-c\sqrt\varepsilon)\le-c\sqrt\varepsilon$, the hypothesis
forces $3d\varepsilon/(2q)\ge c\sqrt\varepsilon$, which is the first
claim. For the second, \cref{si:prop:deltoidapprox} gives
$d^{-1}\log A_d(q)\le d^{-1}\log3-\log R$ with
$\log R=\sqrt\varepsilon+O(\varepsilon)$ and $d\ge C\varepsilon^{-1/2}$,
hence $d^{-1}\log A_d(q)\le-(1-C^{-1}\log3)\sqrt\varepsilon+O(\varepsilon)$,
which is below $\log(1-c\sqrt\varepsilon)=-c\sqrt\varepsilon+O(\varepsilon)$
for small $\varepsilon$ because $1-\log3/C>c$.
\end{proof}

The corollary concerns polynomial approximation. It does not give a
square-root rate for windowed AA with a window of that size, because the
constants of \cref{main:thm:certificate,main:thm:generalmargin} grow with
the coefficients of the certificate, and the coefficients of $p_d$ grow
exponentially in $d$. Conversely, a fixed window does not exclude a
square-root rate: the stationary $\nu$-step methods of
\cref{main:sec:sqrtstatus} attain $R_\nu^{-1}$ with $\nu$ retained iterates.

The argument transfers to the hypocycloids $\Gamma_\nu$ of
\eqref{main:eq:hypocycloid}, whose exterior map is also explicit.

\begin{proposition}[Polynomial approximation on the scaled hypocycloids]
\label{si:prop:hypocycloidapprox}
Let $\nu\ge2$, let $\Gamma_\nu$ be the $\nu$-cusped hypocycloid, and let $A_d^{(\nu)}(q)=\inf\bigl\{\max_{z\in q\Gamma_\nu}|p(z)|:\ \deg p\le d,\ p(1)=1\bigr\}$.
Let $R_\nu>1$ be the unique solution of
$(\nu-1)R_\nu+R_\nu^{-(\nu-1)}=\nu/q$. Then for every
$d\ge1$, $R_\nu^{-d}\le A_d^{(\nu)}(q)\le q^d$ and $A_d^{(\nu)}(q)\le\nu/(R_\nu^d-(\nu-1))$ when $R_\nu^d>\nu-1$,
so $\lim_{d\to\infty}A_d^{(\nu)}(q)^{1/d}=R_\nu^{-1}$, and with
$\varepsilon=1-q$, $R_\nu=1+\sqrt{2\varepsilon/(\nu-1)}+O(\varepsilon)$, $R_\nu^{-1}=1-\sqrt{2\varepsilon/(\nu-1)}+O(\varepsilon)$.
For $\nu=3$ this is \cref{si:prop:deltoidapprox} with a weaker upper
bound; for $\nu=2$ it is the Chebyshev rate on $[-q,q]$.
\end{proposition}
\begin{proof}
The map $\psi_\nu(w)=(q/\nu)\bigl((\nu-1)w+w^{-(\nu-1)}\bigr)$ is the
classical conformal map of $\{|w|>1\}$ onto the exterior of
$q\Gamma_\nu$ with $\psi_\nu(\infty)=\infty$~\cite{eiermann1993}; the
image of the unit circle is $q\gamma_\nu$, and $\psi_\nu'(w)=(q(\nu-1)/\nu)(1-w^{-\nu})$ does not vanish for $|w|>1$. On
$(1,\infty)$ the map increases from $\psi_\nu(1)=q<1$ to infinity, which
gives $R_\nu$. For an admissible $p$, the function $p(\psi_\nu(w))/w^d$ is
analytic for $|w|>1$, including at infinity, and its modulus on $|w|=1$
is at most $\max_{q\Gamma_\nu}|p|$; evaluating at $w=R_\nu$ gives the
lower bound. The bound $q^d$ comes from $p(z)=z^d$.

Writing $\psi_\nu(w)=z$ as
$w^\nu-aw^{\nu-1}+1/(\nu-1)=0$ with $a=\nu z/((\nu-1)q)$, let
$w_1,\ldots,w_\nu$ be the roots with multiplicity and
$S_d(z)=\sum_jw_j(z)^d$. Newton's identities for a monic polynomial whose
only nonzero coefficients besides the leading one are those of
$w^{\nu-1}$ and $w^0$ give $S_0=\nu$, $S_d=a^d\ (1\le d\le\nu-1)$, $S_\nu=a^\nu-\nu/(\nu-1)$, $S_d=aS_{d-1}-S_{d-\nu}/(\nu-1)\ (d>\nu)$,
so $S_d$ is a real polynomial of degree $d$; for $\nu=3$ this is the
recurrence of \cref{si:prop:deltoidapprox}. If $z\in q\Gamma_\nu$, no
root has $|w_j|>1$, because $\psi_\nu$ maps $\{|w|>1\}$ onto the exterior
of $q\Gamma_\nu$; hence $|S_d|\le\nu$ on $q\Gamma_\nu$. At $z=1$, which
lies outside $q\Gamma_\nu$, exactly one root lies in $\{|w|>1\}$, namely
$R_\nu$, and the others have modulus at most one, so
$S_d(1)\ge R_\nu^d-(\nu-1)$. When this is positive,
$p_d=S_d/S_d(1)$ is admissible with
$\max_{q\Gamma_\nu}|p_d|\le\nu/(R_\nu^d-(\nu-1))$.

For the expansion, $(\nu-1)R_\nu+R_\nu^{-(\nu-1)}=\nu+\nu(\nu-1)(R_\nu-1)^2/2+O((R_\nu-1)^3)$, so the defining
equation reads $q\bigl(1+(\nu-1)(R_\nu-1)^2/2+O((R_\nu-1)^3)\bigr)=1$,
and $R_\nu-1=\sqrt{2\varepsilon/(\nu-1)}+O(\varepsilon)$.
\end{proof}

For $d\ge1$ the polynomials $S_d$ are the Faber polynomials of
$q\Gamma_\nu$, and the recurrence for $S_d$ is closely related to the three-term recurrence of order $\nu$ that Cowal, Marshall, and Pollock derive from a random walk with steps $1$ and $1-\nu$~\cite[Section~1.4.2 and Remark~2.2]{cowal2026random}.
Their Theorem~1.2 gives growth at the rate $1+\sqrt{2\varepsilon/(\nu-1)}$
at $1+\varepsilon$, the leading term of $R_\nu$; the proposition adds the
matching lower bound for $A_d^{(\nu)}(q)$. As in
\cref{si:cor:degreeneeded}, degree of order $(1-q)^{-1/2}$ is necessary
and sufficient for a square-root factor on $q\Gamma_\nu$: the real axis
meets $\gamma_\nu$ only at $\theta=0$ and $\theta=\pi$ when $\nu\ge3$, because
$|\sin((\nu-1)\theta)|\le(\nu-1)|\sin\theta|$ with equality only when
$\sin\theta=0$, so $q\Gamma_\nu$ contains the segment
$[q\gamma_\nu(\pi),q]\supseteq[-(\nu-2)q/\nu,q]$, and the Chebyshev
argument applies on it. For $\nu=2$, use $\Gamma_2=[-1,1]$ directly.

\subsection{Proof of \texorpdfstring{\cref{main:thm:contact}}{the contact theorem}, and further remarks}\label{si:sec:contactproof}

\begin{proof}[Proof of \cref{main:thm:contact}]
(a) Cowal, Marshall, and Pollock~\cite[Theorem~1.4]{cowal2026random} show that a polynomial of degree at most $n$ bounded by one on $E_\delta$ has modulus at most $(1+3\varepsilon/(2\delta))^n$ at $1+\varepsilon$, for $0<\varepsilon<1$. Let $p$ have degree at most $n$ with $p(1)=1$, and put $M=\max_{qE_\delta}|p|$. Since $q>1/2$, the point $1/q=1+\varepsilon$ has $\varepsilon=(1-q)/q\in(0,1)$, and applying the bound to $p(qz)/M$ at $z=1/q$ gives $1/M\le(1+3(1-q)/(2\delta q))^n$. Hence, with $A_n$ as in \cref{si:lem:regionalworstcase},
\begin{equation}\label{si:eq:ellipselower}
 A_n(qE_\delta)\ \ge\ \Bigl(1+\frac{3(1-q)}{2\delta q}\Bigr)^{-n}.
\end{equation}
That lemma gives a real normal matrix $W'$ with $\spec(W')\subseteq qE_\delta$ and a unit vector $r_0'$ with $\min_p\norm{p(W')r_0'}=A_n(qE_\delta)$, the minimum taken over the polynomials $p$ of degree at most $n$ with $p(1)=1$. Let $W=\diag(W',q)$ and $r_0=(r_0',0)$. Then $\norm W=q$ and $\spec(W)\subseteq qE_\delta$, since $q\in qE_\delta$, and $p(W)r_0=(p(W')r_0',0)$ for every polynomial $p$. Every polynomial method therefore satisfies $\norm{r_n}=\norm{p_n(W')r_0'}\ge A_n(qE_\delta)$, and \eqref{si:eq:ellipselower} gives (a). The same argument applies whenever the region contains an ellipse tangent at the point of contact~\cite[after Theorem~1.4]{cowal2026random}.

(b) The infimum over real polynomials is at least the one in \cref{si:prop:hypocycloidapprox}, which therefore gives $A_n(q\Gamma_\nu)\ge R_\nu^{-n}$ and the expansion of $R_\nu^{-1}$. The construction of (a), with $q\Gamma_\nu$ in place of $qE_\delta$, gives the lower bound, since $q\in q\Gamma_\nu$ because $1$ is a cusp of $\Gamma_\nu$. The stationary $\nu$-step method has worst-case root-linear factor $R_\nu^{-1}$ on $q\Gamma_\nu$, the smallest asymptotic factor of any semi-iterative method on this region~\cite[Theorem~15 and Corollary~13]{eiermann1985}. The factor can also be verified directly: for an eigenvalue $\lambda$ of $W$, the error recurrence $e_{k+1}=\mu_0\lambda e_k+(1-\mu_0)e_{k+1-\nu}$ has characteristic equation $t^\nu-\mu_0\lambda t^{\nu-1}-(1-\mu_0)=0$. The substitution $t=w/R_\nu$, together with $\mu_0(\nu-1)q/\nu=R_\nu^{-1}$ and $\mu_0-1=R_\nu^{-\nu}/(\nu-1)$, the latter a rewriting of the equation defining $R_\nu$, turns it into $w^\nu-aw^{\nu-1}+1/(\nu-1)=0$ with $a=\nu\lambda/((\nu-1)q)$, the equation $\psi_\nu(w)=\lambda$ of \cref{si:prop:hypocycloidapprox}. For $\lambda\in q\Gamma_\nu$ its roots satisfy $|w|\le1$, so every characteristic root has modulus at most $R_\nu^{-1}$, and at the cusp $\lambda=q$ the root $w=1$ gives equality for the specified initialization. Indeed, for scalar $W=q$ the normalized error generating function is $[1-b\sum_{j=1}^{\nu-1}z^j]/[1-az+bz^\nu]$, with $a=\mu_0q$ and $b=\mu_0-1=R_\nu^{-\nu}/(\nu-1)$. Its denominator has a double zero at $R_\nu$, and its numerator there is $1-(\nu-1)^{-1}\sum_{j=1}^{\nu-1}R_\nu^{j-\nu}>0$, so the factor is $R_\nu^{-1}$. The regional bound also follows from~\cite[Theorem~5 and Corollaries~2 and~10]{niethammervarga1983} applied to the Euler function $\varphi\mapsto\mu_0\varphi/(1-(1-\mu_0)\varphi^\nu)$, whose reciprocal maps the disk of radius $R_\nu$ onto the exterior of $q\Gamma_\nu$; equivalently, the method is the semi-iterative method generated by generalized Faber polynomials of $q\Gamma_\nu$~\cite[Theorems~20 and~21]{eiermann1985}.
\end{proof}

\emph{The matrices of \cref{main:thm:everymemoryquadratic}.} Their spectra are not contained in the deltoid: the auxiliary eigenvalue $-1/2$ lies outside $q\Gamma_3$, since $q\Gamma_3\cap\R=[-q/3,q]$, and as $\rho\uparrow1$ the matrix \eqref{main:eq:momentmatrix} tends to the direct sum of $1$ and a signed cyclic shift whose eigenvalues solve $z^m=-1$, none of which is a cube root of unity.

\emph{The deltoid and the hypocycloids.} For $\nu=3$, \cref{si:prop:deltoidapprox} sharpens the bounds to $R^{-d}\le A_d(q)\le3R^{-d}$, where $R$ solves $q(2R+R^{-2})=3$ and $R^{-1}=1-\sqrt{1-q}+(1-q)/3+O((1-q)^{3/2})$. The reason is the cusp of $q\Gamma_3$ at $q$, which the exterior conformal map of the deltoid unfolds with a square root; the upper bound has the order of the growth estimate of~\cite[Theorem~2]{cowal2026deltoid}. Degree of order $(1-q)^{-1/2}$ is both necessary and sufficient for a square-root factor (\cref{si:cor:degreeneeded}). For general $\nu$, the factor $R_\nu^{-1}$ is the growth rate of the hypocycloid polynomials of~\cite[Theorem~1.2 and Remark~2.2]{cowal2026random}; it interpolates between the Chebyshev rate on an interval and the absence of acceleration on the disk. In the random-walk construction the constant $\sqrt{2/(\nu-1)}$ equals $\sqrt2/\sigma$, where $\sigma^2=\nu-1$ is the variance of the mean-zero walk that steps $+1$ with probability $(\nu-1)/\nu$ and $1-\nu$ with probability $1/\nu$: more cusps make the region rounder, the walk more spread out, and the achievable gain smaller. For $\nu=2$ the stationary method of \cref{main:thm:contact}(b) is the stationary second-order Richardson method, and~\cite[Example~4]{niethammervarga1983} treats $\nu=4$. It is also the fixed-point counterpart of the static momentum power method of Cowal, Marshall, and Pollock, which runs the same three-term recurrence, up to scaling, and converges at the same leading-order rate~\cite[Algorithm~2, Theorem~2.1, and Remark~2.2]{cowal2026random}.

\emph{AA on the deltoid.} The proof of \cref{si:thm:deltoidactivelinear} restricts AA to the spectral subspaces excited by $x_0$ and applies \cref{main:cor:excesscertificate} to the fixed quadratic $p(z)=(9z^2+3z-2)/10$, which satisfies $|p|\le1$ on $\Gamma_3$ with equality only at $z=1$. Its constant is small, $\rho_0$ is not explicit, and the bound neither limits the actual gain from above nor excludes square-root or faster convergence. Certificates of fixed degree cannot settle \eqref{main:eq:sqrttarget}: for every fixed $d\ge2$ the best degree-$d$ margin on $q\Gamma_3$ is of order $1-q$, so \cref{main:thm:generalmargin} gives a gain of order $1-q$ at most (\cref{si:prop:fixeddegreecertificate,si:rem:fixeddegreeenergylimit}).

\section{Capped energy bounds and a robust block estimate}\label{si:sec:cappedbounds}
This section proves the two estimates behind \cref{main:thm:nonlinear}.

\subsection{Bounding the mixing coefficients}

Let $g$ be defined near a fixed point $x_*$, write $r(x)=g(x)-x$, fix a finite cap $\Lambda>1$, and consider capped AA($m$), \eqref{main:eq:nonlinearcap}. Its minimizing residual combination $y_{k+1}$ is the unique minimum-norm point of the feasible residual set; the coefficients need not be unique, and all statements below cover every minimizing selection.

For any feasible residual combination $z$, convex projection gives $\ip{y_{k+1}}{z-y_{k+1}}\ge0$, $\norm{r_k-y_{k+1}}^2 \le\norm{r_k}^2-\norm{y_{k+1}}^2$, $\norm{y_{k+1}}\le\norm{r_k}$,
where $r_k=r(x_k)$ and the last two inequalities use $z=r_k$.
The additional feasible point
$(1+\delta)r_k-\delta r_i$, with
$\delta=(\Lambda-1)/2$, will recover the moments needed by the
polynomial energy argument whenever $r_i$ is retained.

\begin{remark}[The capped feasible set]\label{si:rem:cappedvertices}
The reflected coefficient vectors are exactly the extreme points of the coefficient set.
For a window of $v=m_k+1\ge2$ residuals, let $\mathbf e_1,\ldots,\mathbf e_v$
be the standard basis of $\R^v$. Then $F_\Lambda=\bigl\{\alpha\in\R^v:\ \sum_j\alpha_j=1,\ \sum_j|\alpha_j|\le\Lambda\bigr\} =\conv\bigl\{(1+\delta)\mathbf e_i-\delta\mathbf e_j:\ i\ne j\bigr\}$.
Indeed, $F_\Lambda$ is the section of the cross-polytope
$\{\norm\alpha_1\le\Lambda\}$ by the hyperplane $\sum_j\alpha_j=1$, which
contains none of the vertices $\pm\Lambda\mathbf e_j$ because $\Lambda>1$.
The vertices of the section are therefore the points where the hyperplane
meets an edge. Along an edge joining $\Lambda\mathbf e_i$ and $\Lambda\mathbf e_j$,
or $-\Lambda\mathbf e_i$ and $-\Lambda\mathbf e_j$, the coefficient sum is
constant and equal to $\pm\Lambda$. Along the edge from $\Lambda\mathbf e_i$ to
$-\Lambda\mathbf e_j$, it equals one only at
$(1+\delta)\mathbf e_i-\delta\mathbf e_j$. Consequently the feasible residual
set is $\conv\{(1+\delta)r_i-\delta r_j:\ i\ne j\}$, and $y_{k+1}$ is the
minimum-norm point of these $v(v-1)$ explicit vectors. It can be computed
by a finite active-set method for the nearest point of a polytope, such as
Wolfe's algorithm~\cite{wolfe1976}, without the $2^v$ sign inequalities
that describe the $\ell^1$ constraint. As $\Lambda\downarrow1$ the points
reduce to the retained residuals, and the feasible set to their convex hull.
\end{remark}

\begin{theorem}[A polynomial energy bound with any cap greater than one]
\label{si:thm:cappedgeneral}
Use the polynomial data of \cref{main:thm:generalmargin}, with
$\norm{p(W)}<q^d$ and $m\ge d$, and fix a cap $\Lambda>1$. Let
$\tau_{p,\Lambda}$ be the number defined in \eqref{si:eq:cappedgeneralK} in the proof;
it depends only on $q$, $\norm{p(W)}$, $\Lambda$, and the coefficients of $p$. Then
$0<\tau_{p,\Lambda}<1$, and exact capped affine AA($m$) satisfies $\norm{r_{k+d+1}}^2\le(1-\tau_{p,\Lambda})\,q^{2(d+1)}\norm{r_k}^2$,
and hence $\rhAA(x_0)\le q\,(1-\tau_{p,\Lambda})^{1/(2(d+1))}<q$.
No condition $\Lambda\ge\sum_i|a_i|$ is required. The estimate is
uniform over arbitrary earlier histories, and also holds for varying
caps bounded below by $\Lambda$. If $\Lambda\ge7/5$, then
$\tau_{p,\Lambda}=\tau_p$, so the certified energy factor is the ordinary
affine factor of \cref{main:thm:generalmargin}; this holds more precisely
whenever $(\Lambda+1)/(\Lambda-1)\le(1+\sqrt2)(1+\sqrt3)$, that is, for
caps at least approximately $1.3574$.
\end{theorem}
\begin{proof}
The cap enters only through $(\Lambda+1)/(\Lambda-1)$, and it replaces
$K_p$ by $K_{p,\Lambda}$ in
\begin{equation}\label{si:eq:cappedgeneralK}
 \begin{aligned}
 K_{p,\Lambda}&=K_p+\sum_{j=1}^d|C_j|q^j
 \max\left\{0,\frac1{\Lambda-1}+\frac{(j+1)(\Lambda-1)}{2((j+2)\Lambda-j)}
       -\frac{j+2}{2}-\sqrt{j+1}\right\},\\
 \tau_{p,\Lambda}&=\frac{\Delta_p}{2K_{p,\Lambda}+H_p^2},
 \end{aligned}
\end{equation}
where the $C_j$ are the lag sums from the proof of
\cref{main:thm:generalmargin}.
As in the proof of \cref{main:thm:generalmargin}, assume $\norm{r_k}=1$
and use $T$, $D$, $u_j$, $v_j$, $f_j$, and
$L=1-\norm{u_{d+1}}^2$. Convex projection gives
$\ip{v_j}{f_j}\ge0$, and the affine recurrence gives $\sum_{j=0}^d (\norm{f_j}^2+2\ip{v_j}{f_j}+\norm{Dv_j}^2)=L$.
Fix $1\le j\le d$, and put $t_j=T^ju_0$, $\psi_j=t_j-v_j$,
$h_j=u_0-q^jv_j$, and $b=1-\norm{v_j}^2\le L$.
The joint column estimate from \cref{main:thm:generalmargin} remains
valid because the additional inner products in the energy identity
are nonnegative. It gives $1-\norm{t_j}^2+\norm{\psi_j}^2\le(j+1)b$, $\norm{\psi_j}^2\le(j+1)\bigl(b-2\ip{v_j}{f_j}\bigr)$.

The window retains $r_k=u_0$. The feasible points $u_0$ and
$(1+\delta)q^ju_j-\delta u_0$, with $\delta=(\Lambda-1)/2$, give $0\le\ip{v_j}{h_j}\le q^j(\Lambda+1)\ip{v_j}{f_j}/(\Lambda-1)$.
Consequently $\norm{h_j}^2 =1-q^{2j}\norm{v_j}^2-2q^j\ip{v_j}{h_j} \le1-q^{2j}+q^{2j}b$.
Polarization gives $\ip{u_0}{t_j}-q^j=q^j\bigl(\norm{t_j}^2+\norm{v_j}^2-\norm{\psi_j}^2-2\bigr)/2+\ip{v_j}{h_j}+\ip{h_j}{\psi_j}$.
The first term is nonpositive; more precisely, since $\norm{t_j},\norm{v_j}\le1$ and $1-\norm{v_j}^2=b$, it equals $-q^j\bigl(b+(1-\norm{t_j}^2)+\norm{\psi_j}^2\bigr)/2$, which is at least $-q^j(j+2)b/2$ by the joint estimate and at most $-q^j\bigl(b+\norm{\psi_j}^2\bigr)/2$. The joint estimate and
$\ip{v_j}{h_j}\ge0$ bound the negative part of this expression by $q^j\bigl((j+2)/2+\sqrt{j+1}\bigr)b+\sqrt{(j+1)(1-q^{2j})b}$, because $|\ip{h_j}{\psi_j}|\le\norm{h_j}\norm{\psi_j}\le\bigl(\sqrt{1-q^{2j}}+q^j\sqrt b\bigr)\sqrt{(j+1)b}$.
For the positive part, keep the first term in the form $-q^j\bigl(b+\norm{\psi_j}^2\bigr)/2$, bound the second by
$\ip{v_j}{h_j}\le q^j\frac{\Lambda+1}{\Lambda-1}\ip{v_j}{f_j}\le q^j\frac{\Lambda+1}{2(\Lambda-1)}\bigl(b-\norm{\psi_j}^2/(j+1)\bigr)$,
which follows from the second joint estimate, and the third by $\ip{h_j}{\psi_j}\le\norm{h_j}\norm{\psi_j}\le\bigl(\sqrt{1-q^{2j}}+q^j\sqrt b\bigr)\norm{\psi_j}$.
In the sum of the three bounds the coefficient of $b$ is $q^j\bigl[\frac{\Lambda+1}{2(\Lambda-1)}-\frac12\bigr]=\frac{q^j}{\Lambda-1}$ and the coefficient of $\norm{\psi_j}^2$ is $-q^j\bigl[\frac12+\frac{\Lambda+1}{2(j+1)(\Lambda-1)}\bigr]=-q^j\frac{(j+2)\Lambda-j}{2(j+1)(\Lambda-1)}$, so the terms other than $\sqrt{1-q^{2j}}\norm{\psi_j}$ are at most
\[
 q^j\left[
 \frac{b}{\Lambda-1}
 -\frac{(j+2)\Lambda-j}{2(j+1)(\Lambda-1)}\norm{\psi_j}^2
 +\sqrt b\norm{\psi_j}\right]
 \le q^j\left(\frac1{\Lambda-1}+\frac{(j+1)(\Lambda-1)}{2((j+2)\Lambda-j)}\right)b,
\]
where the last inequality completes the square in $\norm{\psi_j}$, and $\norm{\psi_j}\le\sqrt{(j+1)b}$ in the remaining term.
Thus, with $b\le L$,
\begin{equation}\label{si:eq:capmoment}
 \left|\ip{u_0}{T^ju_0}-q^j\right|
 \le q^j\max\left\{
      \frac{j+2}{2}+\sqrt{j+1},
      \frac1{\Lambda-1}+\frac{(j+1)(\Lambda-1)}{2((j+2)\Lambda-j)}
      \right\}L
       +\sqrt{(j+1)(1-q^{2j})L}.
\end{equation}
The Gram-matrix expansion from \cref{main:thm:generalmargin} gives
\begin{equation}\label{si:eq:capexcess}
 \Delta_p\le K_{p,\Lambda}L+G_p\sqrt L.
\end{equation}
Since $\norm{p(W)}<q^d$, we have $E_p\le\Delta_p$ and $G_p\le H_p\sqrt{E_p}$.
Young's inequality gives
$G_p\sqrt L\le\Delta_p/2+H_p^2L/2$, proving
$L\ge\tau_{p,\Lambda}$. Its range follows from $K_{p,\Lambda}\ge K_p$
and \cref{main:thm:generalmargin}. This proves the block estimate;
block iteration and one-step contraction give the root factor.
The proof only used retained residuals from the current block and
caps at least $\Lambda$, proving uniformity over histories and varying caps.

The cap correction vanishes whenever
$(\Lambda+1)/(\Lambda-1)\le(1+\sqrt2)(1+\sqrt3)$.
Indeed, the second entry in the maximum of \eqref{si:eq:capmoment}
increases from its value at $j=1$ by at most $(j-1)/2$, whereas the
first increases by at least $(j-1)/2$. At $j=1$, comparison of the
entries is equivalent to the stated condition. In particular,
$\Lambda\ge7/5$ gives $(\Lambda+1)/(\Lambda-1)\le6<(1+\sqrt2)(1+\sqrt3)$.
\end{proof}

\begin{corollary}[Positive energy excess and the sharper factor]
\label{si:cor:cappedexcesscertificate}
The preceding block bound also holds under the weaker assumption
$\Delta_p=q^{2d}+E_p-\norm{p(W)}^2>0$, if
$\tau_{p,\Lambda}$ is replaced by the following number, in which
$K_{p,\Lambda}$ is defined in \eqref{si:eq:cappedgeneralK}:
\begin{equation}\label{si:eq:cappedsharptau}
 \widehat\tau_{p,\Lambda}
 =\frac{4\Delta_p^2}
 {\left(\sqrt{G_p^2+4K_{p,\Lambda}\Delta_p}+G_p\right)^2}.
\end{equation}
This number belongs to $(0,1)$ and improves the preceding factor when
$\norm{p(W)}<q^d$. For $\Lambda\ge7/5$, it equals $\widehat\tau_p$ computed with
\eqref{si:eq:generalK}.
\end{corollary}
\begin{proof}
The proof through \eqref{si:eq:capexcess} does not use
$\norm{p(W)}<q^d$. Solving that increasing scalar inequality for $\sqrt L$
gives \eqref{si:eq:cappedsharptau}. Its range and comparison follow as in
\cref{main:cor:excesscertificate}, with $K_{p,\Lambda}\ge K_p$.
\end{proof}

The strict inequality $\Lambda>1$ is necessary for a statement covering
every contraction. For $g(x)=qx$ on the real line, with $0<q<1$, a cap
of one permits only nonnegative coefficients. Starting from a nonzero
point, the latest Picard residual is the smallest same-sign residual in
every window, so it is the minimum-norm convex combination. Induction
gives exactly Picard iteration at every memory, although
$p(z)=(z-q)/(1-q)$ is a degree-one certificate with $p(q)=0$.

\subsection{Stability under small errors}

The block argument involves finitely many steps, so it tolerates small
errors in the recurrence. The next lemma quantifies this; the errors it
allows include those caused by the nonlinear remainder and by
approximate least-squares solutions.

\begin{lemma}[A robust block estimate]
\label{si:lem:compactblock}
Fix $J\in\R^{n\times n}$ with $\norm J=q\in(0,1)$, polynomial data
satisfying either \cref{si:thm:cappedgeneral} or
\cref{si:cor:cappedexcesscertificate}, and a cap $\Lambda>1$; set $\delta=(\Lambda-1)/2$.
Let $\chi=(1-\tau)^{1/2}q^{d+1}$ be the corresponding capped block factor, where $\tau$ is $\tau_{p,\Lambda}$ of \cref{si:thm:cappedgeneral} or its analogue in \cref{si:cor:cappedexcesscertificate}. Suppose vectors
$r_0,\ldots,r_{d+1}$ and $y_1,\ldots,y_{d+1}$ satisfy $r_{j+1}=Jy_{j+1}+\xi_j$ for $0\le j\le d$
and
\begin{equation}
 \ip{y_{j+1}}{r_j-y_{j+1}} \ge0,\qquad
 \ip{y_{j+1}}{r_0-y_{j+1}} \ge0, \qquad\label{si:eq:compactgeometry}
 \ip{y_{j+1}}{(1+\delta)r_j-\delta r_0-y_{j+1}}\ge0.    
\end{equation}
There is a nondecreasing function $\mathcal M:[0,1]\to[0,\infty)$,
with $\mathcal M(t)\to0$ as $t\downarrow0$, such that for every $T>0$,
\begin{equation}\label{si:eq:compactblockbound}
 \norm{r_0}\le T,\quad
 \max_{j\le d}\norm{\xi_j}\le\varepsilon T,
 \quad0\le\varepsilon\le1
 \quad\Longrightarrow\quad
 \norm{r_{d+1}}\le\chi\norm{r_0}+\mathcal M(\varepsilon)T.
\end{equation}
The modulus is independent of the vectors and of any earlier history.
There is an unconstrained version with the corresponding ordinary
energy block factor:
replace \eqref{si:eq:compactgeometry} by $\ip{y_{j+1}}{r_j-y_{j+1}}=0$, $\ip{y_{j+1}}{r_0-y_{j+1}}=0$.
\end{lemma}
\begin{proof}
First take all defects to be zero. The proof of
\cref{si:thm:cappedgeneral,si:cor:cappedexcesscertificate} uses only the
recurrence, the three displayed geometric inequalities, and the
operator certificate. It does not require that each $y_{j+1}$ belong
to a specified feasible set. Thus the same proof gives
$\norm{r_{d+1}}\le\chi\norm{r_0}$ for these abstract vectors. If $r_0=0$,
the first inequality forces $y_1=0$ and the recurrence forces $r_1=0$;
induction gives the same conclusion.

By homogeneity it is enough to take $T=1$. The first geometric
inequality gives $\norm{y_{j+1}}\le\norm{r_j}$, so $\norm{r_j}\le q^j+\varepsilon\sum_{i=0}^{j-1}q^i \le1+j$, $\norm{y_{j+1}}\le1+j$.
Define $\mathcal M(\varepsilon)$ as the supremum of
$(\norm{r_{d+1}}-\chi\norm{r_0})_+$ over all such tuples with
$\norm{r_0}\le1$ and defects bounded by $\varepsilon$.
These finite-dimensional feasible sets are closed and bounded, and
contain the zero tuple. Thus the supremum is finite, is attained, and
is nondecreasing in $\varepsilon$.
If $\mathcal M(\varepsilon)$ did not tend to zero, there would be a sequence
$\varepsilon_\nu\downarrow0$ and feasible tuples for which the excess
stays bounded away from zero. The displayed uniform bounds give a
convergent subsequence of all their vectors. In the limit the defects
vanish and the geometric inequalities persist. The zero-defect bound
then contradicts the positive limiting excess. This proves the claim.
The unconstrained case uses the same argument and the affine projection
identities in \cref{main:thm:generalmargin,main:cor:excesscertificate}.
\end{proof}

This compactness proof suffices for the local theorem in a fixed
finite-dimensional space. A direct energy proof gives an explicit modulus of order
$\sqrt\varepsilon$ with constants independent of dimension
(supplementary \cref{si:prop:explicitenergyperturbation}). Neither proof requires continuity
of the coefficient selection.

\section{Nonlinear extensions and error bounds}\label{si:sec:sinonlinear}
The local theorem of the article uses a fixed coefficient bound and a
remainder that vanishes at the fixed point. This section makes its error bound explicit, shows that without control of the minimizing coefficients some exact least-squares selections diverge, and gives the proof of the theorem together with its version for inexact coefficient solves.

\subsection{An explicit bound for the effect of errors}
\label{si:sec:explicitenergyperturbation}

The following direct energy estimate proves \cref{si:lem:compactblock}
with a dimension-independent formula. The constants are conservative;
their purpose is to make the perturbation bound explicit.

\begin{proposition}[An explicit square-root perturbation modulus]
\label{si:prop:explicitenergyperturbation}
Use the hypotheses and notation of \cref{si:lem:compactblock}. Then
\eqref{si:eq:compactblockbound} holds with
\begin{equation}\label{si:eq:explicitenergymodulus}
 \mathcal M(\varepsilon)=(d+1)\varepsilon
               +C_d\sqrt{\varepsilon},
 \qquad 0\le\varepsilon\le1.
\end{equation}
Here the explicit constant is
\begin{equation}\label{si:eq:explicitenergyuniformconstant}
 C_d
 =\frac12\sqrt{\frac{d+2}{d}}\left[
 12d+11+\frac{\Lambda+1}{\Lambda-1}+2\sqrt2\,\sqrt{3d+7}+\sqrt{2d+2}\right].
\end{equation}
In the unconstrained case the term $(\Lambda+1)/(\Lambda-1)$ is omitted.
The modulus is uniform over $q\in(0,1)$, the polynomial coefficients,
the operator, and the ambient dimension. It applies to the
Young-inequality and positive-excess energy factors in the main article,
and its proof uses $p(J)$ only through the vector $p(J)r_0$.
\end{proposition}

\begin{proof}
Put $N=d+1$, and let $\tau$ be the loss constant corresponding to the
chosen factor, so $\vartheta=\sqrt{1-\tau}$ and $\chi=q^N\vartheta$.
For this proof put $U=\sum_{i=0}^d|a_i|q^i$ and
$V^2=\sum_{i=0}^da_i^2q^{2i}$, and use the lag sums $C_j$
and the constants $G_p,E_p,\Delta_p$ from
\cref{main:thm:generalmargin}.
In the unconstrained case, read $K_{p,\Lambda}$ below as $K_p$;
the one-sided choice in \cref{si:rem:onesidedmoments} is also admissible.
For $1\le j\le d$, set
\[
 D_j={}q^j\left(5j+8+\frac{\Lambda+1}{\Lambda-1}+\sqrt{2(3j+7)}
                 +\sqrt{3j+7}+\sqrt{2j+2}\right)
       +\sqrt{(1-q^{2j})(3j+7)}.
\]
Put
\begin{equation}\label{si:eq:explicitenergyconstants}
 C=\sum_{j=1}^d|C_j|D_j+(7d+3)U^2+2V^2,
 \qquad \Omega=\sum_{i=1}^Nq^{-i}.
\end{equation}

Normalize the envelope bound in \eqref{si:eq:compactblockbound} to one,
and write $R=\norm{r_0}\le1$ and $\norm{\xi_j}\le\varepsilon$.
Put $R'=R+\Omega\varepsilon$. If $R'=0$, all vectors vanish and the
claim follows. Otherwise define $h=\Omega\varepsilon/(R')$, $u_j=r_j/(q^jR')$, $v_j=y_{j+1}/(q^jR')$, $f_j=u_j-v_j$, $\zeta_j=\xi_j/(q^{j+1}R')$, $T=J/q$.
The first geometric condition implies
$\norm{r_{j+1}}\le q\norm{r_j}+\varepsilon$. Consequently $\norm{u_j},\norm{v_j}\le1$, $\norm{u_0}=1-h$, $\sum_{j=0}^{d}\norm{\zeta_j}\le h\le1$, $u_{j+1}=Tv_j+\zeta_j$.
Write $\ell=1-\norm{u_N}^2$, and let $D=(I-T^*T)^{1/2}$ be the defect
operator, as in the proof of \cref{main:thm:generalmargin}.
Since both $\norm{u_{j+1}}$ and $\norm{Tv_j}$ are at most one,
their squared norms differ by at most $2\norm{\zeta_j}$.
Using $\ip{v_j}{f_j}\ge0$ and telescoping gives $\mathcal E:=\sum_{j=0}^d\bigl(\norm{f_j}^2+2\ip{v_j}{f_j}+\norm{Dv_j}^2\bigr)\le\ell+2h\le3$.
The same comparison from any index to the endpoint gives
$\norm{v_j}^2\ge1-\ell-2h$.

\emph{Raw powers.}
Let $t_j=T^ju_0$ and $\psi_j=t_j-v_j$. Expanding the recurrence gives $t_j=v_j+\sum_{i=0}^jT^{j-i}f_i -\sum_{i=0}^{j-1}T^{j-1-i}\zeta_i$.
Every noise prefix has norm at most $h$. The telescoping column
estimate in \eqref{si:eq:generalrawenergy} therefore gives
\begin{align}
 \norm{u_0}^2-\norm{t_j}^2
 &\le\left(\sqrt{(j+1)\mathcal E}+\sqrt j\,h\right)^2
 \le(j+1)\ell+(7j+3)h,\notag\\
 1-\norm{t_j}^2&\le(j+1)\ell+(7j+5)h.
 \label{si:eq:explicitenergyraw}
\end{align}
Here we used $\mathcal E\le\ell+2h\le3$, $h^2\le h$, and
$2\sqrt{3j(j+1)}\le4j+1$. Likewise, $\norm{\psi_j}^2 \le(\sqrt{(j+1)\mathcal E}+h)^2 \le(j+1)\ell+(3j+7)h$,
using $2\sqrt{3(j+1)}\le j+4$.

\emph{Moments.}
Fix $1\le j\le d$ and put $b=1-\norm{v_j}^2\le\ell+2h$.
Telescoping the energy only through the projection at index $j$ gives $\sum_{i<j}\bigl(\norm{f_i}^2+2\ip{v_i}{f_i}+\norm{Dv_i}^2\bigr) +\norm{f_j}^2+2\ip{v_j}{f_j} \le\norm{u_0}^2-\norm{v_j}^2+2h=b+h^2$.
Apply the joint column estimate from \cref{main:thm:generalmargin}
to $(Dt_0,\ldots,Dt_{j-1},\psi_j)$.
Each noise column, like each correction column, has norm equal to
that of its input by contraction-defect telescoping; their combined
norm is at most $h$. Therefore
\[
 \begin{split}
 \left(\norm{u_0}^2-\norm{t_j}^2+\norm{\psi_j}^2\right)^{1/2}
 &\le\sqrt{(j+1)(b+h^2)}+h,\\
 \norm{\psi_j}
 &\le\sqrt{(j+1)\bigl(b-2\ip{v_j}{f_j}+h^2\bigr)}+h.
 \end{split}
\]
Since $b\le1$, $h\le1$, and $2\ip{v_j}{f_j}\le b$, these imply
\[
 \begin{split}
 1-\norm{t_j}^2+\norm{\psi_j}^2
 &\le(j+1)b+\bigl(j+4+2\sqrt{2j+2}\bigr)h,\\
 \norm{\psi_j}^2
 &\le(j+1)\bigl(b-2\ip{v_j}{f_j}\bigr)+(3j+7)h.
 \end{split}
\]
For the second inequality we used
$j+2+2\sqrt{2j+2}\le3j+7$.

Set $h_j=u_0-q^jv_j$. The reflected competitors give $0\le\ip{v_j}{h_j}\le q^j(\Lambda+1)\ip{v_j}{f_j}/(\Lambda-1)$, $\norm{h_j}^2\le1-q^{2j}+q^{2j}b$.
In the unconstrained case $\ip{v_j}{h_j}=\ip{v_j}{f_j}=0$,
so every bound below holds with the terms containing $\Lambda$
deleted; in particular only the first entry of the maximum in
\eqref{si:eq:explicitenergymoments} remains.
For $\mu_j=\ip{u_0}{t_j}$, polarization gives $\mu_j-q^j=q^j\bigl(\norm{t_j}^2+\norm{v_j}^2-\norm{\psi_j}^2-2\bigr)/2+\ip{v_j}{h_j}+\ip{h_j}{\psi_j}$.
The first term is nonpositive; as in the proof of \cref{si:thm:cappedgeneral}, it equals $-q^j\bigl(b+(1-\norm{t_j}^2)+\norm{\psi_j}^2\bigr)/2$ and is therefore at most $-q^j\bigl(b+\norm{\psi_j}^2\bigr)/2$. The joint estimate bounds the negative
part as in \eqref{si:eq:capmoment}, with the additional noise term
$q^j\bigl((j+4)/2+\sqrt{2j+2}\bigr)h$; in the term $\ip{h_j}{\psi_j}$ of the negative part, bound $\norm{\psi_j}$ by the raw-power estimate $\norm{\psi_j}^2\le(j+1)\ell+(3j+7)h$ rather than through $b$, and use $\sqrt{b\ell}\le\ell+\sqrt{2h}$.
Explicitly, the first joint estimate gives $-q^j\bigl(b+(1-\norm{t_j}^2)+\norm{\psi_j}^2\bigr)/2\ge-q^j(j+2)b/2-q^j\bigl((j+4)/2+\sqrt{2j+2}\bigr)h$, and
$\ip{h_j}{\psi_j}\ge-\bigl(\sqrt{1-q^{2j}}+q^j\sqrt b\bigr)\bigl(\sqrt{(j+1)\ell}+\sqrt{(3j+7)h}\bigr)$,
where $q^j\sqrt{(j+1)b\ell}\le q^j\sqrt{j+1}\,\bigl(\ell+\sqrt{2h}\bigr)$ and $q^j\sqrt b\sqrt{(3j+7)h}\le q^j\sqrt{(3j+7)h}$ because $b\le1$; with $b\le\ell+2h$ and $h\le\sqrt h$ this gives
\[
 \begin{aligned}
 \mu_j-q^j&\ge-q^j\Bigl(\frac{j+2}{2}+\sqrt{j+1}\Bigr)\ell-\sqrt{(j+1)(1-q^{2j})\ell}\\
 &\quad-\Bigl(q^j\bigl[\tfrac{3j}{2}+4+2\sqrt{2(j+1)}+\sqrt{3j+7}\bigr]+\sqrt{(1-q^{2j})(3j+7)}\Bigr)\sqrt h ,
 \end{aligned}
\]
and $2\sqrt{2(j+1)}\le(2+\sqrt2)\sqrt{j+1}$.
For the positive part, truncate $\norm{\psi_j}$ at
$\sqrt{(j+1)\bigl(b-2\ip{v_j}{f_j}\bigr)}$, that is, put $\psi_j'=\min\bigl\{\norm{\psi_j},\sqrt{(j+1)(b-2\ip{v_j}{f_j})}\bigr\}$.
The truncation changes its value by at most $\sqrt{(3j+7)h}$, since $\norm{\psi_j}\le\psi_j'+\sqrt{(3j+7)h}$ by the second joint estimate.
The scalar function $z\mapsto-z^2/2+\sqrt b\,z$ has derivative at
most one for $z\ge0$, so the completion-of-the-square argument in
\eqref{si:eq:capmoment} bounds the positive contribution by $q^j\bigl(1/(\Lambda-1)+(j+1)(\Lambda-1)/(2((j+2)\Lambda-j))\bigr)b+q^j\sqrt{(3j+7)h}$,
apart from $\sqrt{1-q^{2j}}\norm{\psi_j}$.
In detail, $\ip{v_j}{h_j}\le q^j\frac{\Lambda+1}{\Lambda-1}\ip{v_j}{f_j}\le q^j\frac{\Lambda+1}{2(\Lambda-1)}\bigl(b-\psi_j'^2/(j+1)\bigr)$ by the definition of $\psi_j'$, while $-q^j(b+\norm{\psi_j}^2)/2\le-q^j(b+\psi_j'^2)/2$ and $q^j\sqrt b\norm{\psi_j}\le q^j\sqrt b\,\psi_j'+q^j\sqrt{(3j+7)h}$; the sum of these three bounds is
$q^j\bigl[b/(\Lambda-1)-\frac{(j+2)\Lambda-j}{2(j+1)(\Lambda-1)}\psi_j'^2+\sqrt b\,\psi_j'\bigr]+q^j\sqrt{(3j+7)h}$, exactly as in the proof of \cref{si:thm:cappedgeneral}, and completing the square in $\psi_j'$ gives the stated bound.
For this remaining term, the earlier bound on $\norm{\psi_j}^2$ gives $\sqrt{1-q^{2j}}\norm{\psi_j} \le\sqrt{(j+1)(1-q^{2j})\ell} +\sqrt{(1-q^{2j})(3j+7)h}$.
Finally use $b\le\ell+2h$ and $h\le\sqrt h$.
The negative estimate contributes at most the noise coefficient $q^j\bigl[3j/2+4+(2+\sqrt2)\sqrt{j+1}+\sqrt{3j+7}\bigr]$,
and the positive estimate contributes at most
$q^j\bigl((\Lambda+1)/(\Lambda-1)+\sqrt{3j+7}\bigr)$, since the
coefficient of $b$ above is at most $(\Lambda+1)/(2(\Lambda-1))$ (its second summand is at most $1/2$), so that $b\le\ell+2h$ contributes at most $q^j\frac{\Lambda+1}{\Lambda-1}h$ to the noise.
Both are bounded by the $q^j$ part of $D_j$; for the negative bound,
use $2\sqrt{j+1}\le j+2$. Thus
\begin{equation}\label{si:eq:explicitenergymoments}
 |\mu_j-q^j|
 \le q^j\max\left\{
       \frac{j+2}{2}+\sqrt{j+1},
       \frac1{\Lambda-1}+\frac{(j+1)(\Lambda-1)}{2((j+2)\Lambda-j)}
       \right\}\ell
       +\sqrt{(j+1)(1-q^{2j})\ell}+D_j\sqrt h.
\end{equation}

\emph{Polynomial energy.}
For unconstrained projections, the identity in
\cref{si:rem:onesidedmoments} also holds under the present normalization.
It gives the one-sided bound $\mu_j-q^j\le\sqrt{(j+1)(1-q^{2j})\ell} +\sqrt{(1-q^{2j})(3j+7)h}$.
Consequently the one-sided $K_p$ is admissible in the unconstrained case
with the same noise constants $D_j$.
The operator telescoping identity remains exact: $\ip{t_i}{t_j} =\mu_{j-i}-\sum_{a=0}^{i-1}\ip{Dt_a}{Dt_{a+j-i}}$ for $0\le i\le j\le d$.
Both dissipation sums in Cauchy--Schwarz are bounded by
$(d+1)\ell+(7d+3)h$ from \eqref{si:eq:explicitenergyraw}.
The diagonal moment is $\mu_0=(1-h)^2$, whose loss relative to one
is at most $2h$. Explicitly,
$\norm{p(J)u_0}^2=\sum_{i,j}\hat a_i\hat a_j\ip{t_i}{t_j}\ge\mu_0\sum_i\hat a_i^2+\sum_{h'\ge1}C_{h'}\mu_{h'}-\bigl((d+1)\ell+(7d+3)h\bigr)U^2$,
with $\mu_0\sum_i\hat a_i^2\ge V^2-2hV^2$ and $\sum_{h'}C_{h'}\mu_{h'}\ge\sum_{h'}C_{h'}q^{h'}-\sum_{h'}|C_{h'}|\,|\mu_{h'}-q^{h'}|$; inserting \eqref{si:eq:explicitenergymoments}, \eqref{si:eq:generalgrambaseline}, and $h\le\sqrt h$, the Gram reduction and the same
moment bound as in \eqref{si:eq:capmoment} give
\begin{equation}\label{si:eq:explicitenergyexcess}
 \norm{p(J)}^2\ge\norm{p(J)u_0}^2
 \ge q^{2d}+E_p-K_{p,\Lambda}\ell-G_p\sqrt\ell-C\sqrt h.
\end{equation}
The two additional Gram losses are $(7d+3)U^2h$ and $2V^2h$,
explaining the last terms in \eqref{si:eq:explicitenergyconstants}.
The upper bound uses only $\norm{p(J)u_0}\le\norm{p(J)}\norm{u_0}\le\norm{p(J)}$.

Set $f(x)=K_{p,\Lambda}x+G_p\sqrt x$ for $x\ge0$.
For the positive-excess factor, $f(\tau)=\Delta_p$; for the
Young-inequality factor, $f(\tau)\le\Delta_p$ because
$G_p\le H_p\sqrt{E_p}$ and $E_p\le\Delta_p$.
Since $f(x)-f(y)\ge K_{p,\Lambda}(x-y)$ for $x\ge y\ge0$,
\eqref{si:eq:explicitenergyexcess} implies $\ell\ge\tau-C\sqrt h/K_{p,\Lambda}$: if $\ell<\tau$, then $K_{p,\Lambda}(\tau-\ell)\le f(\tau)-f(\ell)\le\Delta_p-(\Delta_p-C\sqrt h)$.
Concavity of the square root at $1-\tau>0$ then gives $\norm{u_N}=\sqrt{1-\ell}\le\sqrt{1-\tau+C\sqrt h/K_{p,\Lambda}}\le\vartheta+C\sqrt h/(2K_{p,\Lambda}\vartheta)$.
Undoing the normalization yields $\norm{r_N}\le\chi R+\chi\Omega\varepsilon+q^NC\sqrt{\Omega\varepsilon R'}/(2K_{p,\Lambda}\vartheta)$.
To bound the noise coefficient uniformly, Cauchy--Schwarz gives
$q^j\sqrt{3j+7}+\sqrt{(1-q^{2j})(3j+7)}\le\sqrt{2(3j+7)}$.
Since $q^j\le1$ and $j\le d$, every $D_j$ is at most $D_*=5d+8+(\Lambda+1)/(\Lambda-1)+2\sqrt2\,\sqrt{3d+7}+\sqrt{2d+2}>2$.
Also $\sum_{j=1}^d|C_j|\le U^2-V^2$ and
$K_{p,\Lambda}\ge NU^2$. Hence $C\le(7d+3+D_*)U^2$, $C/(2K_{p,\Lambda}\vartheta) \le(7d+3+D_*)/(2\sqrt{N(N-1)})$.
For the second bound, $\tau\le1/N$ follows from
$\Delta_p\le q^{2d}+E_p\le U^2$ and
$\tau\le\Delta_p/K_{p,\Lambda}$.
Finally, $q^N\Omega\le N$, $q^{2N}\Omega(1+\Omega)\le N(N+1)$, $R'\le1+\Omega$.
Substitution proves
\eqref{si:eq:explicitenergymodulus} and~\eqref{si:eq:explicitenergyuniformconstant}.
Rescaling restores an arbitrary envelope bound. Every estimate uses
operator norms and inner products only, so the statement also holds
for bounded operators on a real Hilbert space under the same
certificate and block-geometry hypotheses.
\end{proof}

The square root in \eqref{si:eq:explicitenergymodulus} describes stability
with respect to recurrence defects. It is unrelated to a square-root
speedup in $1-q$ for the unperturbed AA convergence factor.

\subsection{Failure without control of minimizing selections}

\begin{proposition}[Failure without control of minimizing selections]\label{si:prop:unboundedcounterexample}
There is an analytic scalar fixed-point map, contractive in a
neighborhood of its fixed point, whose derivative at that point has a
degree-one polynomial certificate, for which
ordinary exact AA($m$), for every fixed $m\ge2$, admits divergent
minimizing-coefficient selections from arbitrarily close one-point
initializations.
\end{proposition}
\begin{proof}
Fix $0<q<1$, put $h=1-q$, and take $g(x)=qx+x^2$, $r(x)=x(x-h)$, $x_*=0$.
Here $g'(0)=q$, and $p(z)=(z-q)/(1-q)$ has $p(g'(0))=0$.
Moreover $|g'(x)|\le(1+q)/2<1$ for $|x|<(1-q)/4$, so the map is
contractive on that neighborhood.
Start with $x_0=t$, where $0<t<h/4$. The first step is
$x_1=qt+t^2$, so $0<x_1<t$ and $t+x_1<h$. The two residuals are distinct,
and their unique affine minimum-norm combination is zero. The next
iterate is the secant intercept $x_2=tx_1/(t+x_1-h)<0$.
Thus $x_0,x_1,x_2$ are distinct and none is a zero of $r$.
At the next solve, an exact minimizing selection can send the third
iterate to any prescribed real number. For any three distinct retained points $x_i$, invertibility of
its Vandermonde matrix gives coefficients satisfying, for any chosen $M$, $\sum_i\alpha_i=1$, $\sum_i\alpha_ix_i=M$, $\sum_i\alpha_ix_i^2=hM$.
They give $\sum_i\alpha_ir(x_i)=0$, hence are exact global least-squares
minimizers, while $x_{k+1}=\sum_i\alpha_ig(x_i)=(q+h)M=M$.
For a window deeper than two, use its last three points and assign zero
coefficients to the older points. The mixed residual remains zero and
is therefore still globally minimal.
Choose $x_n=(n-1)h$ for every $n\ge3$. All retained triples remain
distinct, the sequence diverges, and it never visits either residual-zero
point $0$ or $h$. The stopping convention therefore does not intervene.
At the first escaping solve,
$\sum_i\alpha_ix_i^2=2h^2$ forces
$\norm{\alpha}_1\ge2h^2/\max_{0\le i\le2}|x_i|^2\to\infty$ as $t\downarrow0$.
\end{proof}
This counterexample rules out a local theorem covering every exact
minimizing selection. It does not settle the behavior of a prescribed coefficient-selection rule, such as the minimum-Euclidean-norm rule.

\subsection{Proof of the local nonlinear theorem and inexact solves}\label{si:sec:nonlinearproof}
\begin{proof}[Proof of \cref{main:thm:nonlinear}]
Differentiability at $x_*$ gives
\(g(x_*+w)=x_*+Jw+\mathcal N(w)\) with
\(\norm{\mathcal N(w)}\le\varpi(\norm w)\norm w\), where $\varpi$ is nondecreasing and
$\varpi(t)\to0$ as $t\downarrow0$.
Choose $r_*>0$ so that this holds for $\norm w\le r_*$ and
$\varpi(r_*)\le(1-q)/2$, and put $a_*=(1-q)/2$.
For every error in this ball,
\begin{equation}\label{si:eq:nonlinearresidualequi}
 a_*\norm w\le\norm{r(x_*+w)}\le(1+q+\varpi(r_*))\norm w,
\end{equation}
because $\norm{(J-I)w}\ge(1-q)\norm w$. In particular the fixed point
is the only residual-zero point in the ball.

\emph{Invariance and preliminary convergence.}
Let $w_j=x_j-x_*$ and let the retained residual envelope be $T_k=\max_{\max(0,k-m)\le j\le k}\norm{r_j}$.
If the current history lies in the ball and $T_k\le T_*$, every retained
error has norm at most $T_k/a_*$, and coefficient feasibility bounds the
next error by
\(\norm{w_{k+1}}\le GT_k/a_*\), where \(G=\Lambda(q+\varpi(r_*))\).
Choose $T_*>0$ with $\max\{1,G\}\,T_*/a_*\le r_*$, so the new iterate also
lies in the ball. For
$\overline{\mathcal N}_k=\sum_j\alpha_j\mathcal N(w_j)$, the identities
$y_{k+1}=(J-I)\sum_j\alpha_jw_j+\overline{\mathcal N}_k$ and
$w_{k+1}=J\sum_j\alpha_jw_j+\overline{\mathcal N}_k$ give
\begin{equation}\label{si:eq:nonlinearperturbation}
 r_{k+1}=Jy_{k+1}+\xi_k,
 \qquad \xi_k=\mathcal N(w_{k+1})-\overline{\mathcal N}_k,
\end{equation}
and therefore $\norm{\xi_k}\le\eta(T_k)T_k$, $\eta(t)=(G+\Lambda)\varpi\bigl(\max\{1,G\}\,t/a_*\bigr)/a_*\to0$.
Reduce $T_*$ so that $\mu=q+\eta(T_*)<1$. Since $r_k$ is feasible,
$\norm{y_{k+1}}\le\norm{r_k}$, and \eqref{si:eq:nonlinearperturbation} gives
$\norm{r_{k+1}}\le q\norm{r_k}+\eta(T_k)T_k\le\mu T_k$.
Starting with
$\norm{x_0-x_*}\le\min\{r_*,T_*/(1+q+\varpi(r_*))\}$, induction proves that all
iterates remain in the ball, $T_{k+1}\le T_k\le T_*$, and
$T_{k+m+1}\le\mu T_k$. Thus $T_k\to0$ and $x_k\to x_*$. Every constant in
this induction is independent of the choice among minimizing coefficients.

\emph{The asymptotic factor.}
It remains to consider a nonterminating run. Throughout a block beginning
at $k$, monotonicity of $T_k$ bounds the defects by $\eta(T_k)T_k$. At every
projection in that block the start $r_k$ is retained, because $d\le m$, so
$r_k$, the current residual, and their reflected combination are feasible.
They give exactly the geometry of \cref{si:lem:compactblock}, which yields,
for all sufficiently large $k$,
\begin{equation}\label{si:eq:nonlineargeneralforced}
 \norm{r_{k+d+1}}\le\chi\norm{r_k}+\mathcal M(\eta(T_k))\,T_k,
 \qquad\mathcal M(\eta(t))\longrightarrow0\quad(t\downarrow0).
\end{equation}
Fix $\lambda\in(\chi^{1/(d+1)},1)$ and choose $\epsilon>0$ so that
$\chi\lambda^{-d-1}+\epsilon\lambda^{-(d+1+m)}<1$.
There is $k_0\ge m$ for which $\mathcal M(\eta(T_k))\le\epsilon$ whenever
$k\ge k_0$. Choose $C_\lambda$ such that $\norm{r_j}\le C_\lambda\lambda^j$
for $k_0-m\le j\le k_0+d$. If this bound holds through index $k+d$, then
$T_k\le C_\lambda\lambda^{k-m}$, and
\eqref{si:eq:nonlineargeneralforced} gives $\norm{r_{k+d+1}}\le C_\lambda\lambda^{k+d+1} \bigl(\chi\lambda^{-d-1}+\epsilon\lambda^{-(d+1+m)}\bigr) \le C_\lambda\lambda^{k+d+1}$.
Induction proves the bound at every later index. Letting
$\lambda\downarrow\chi^{1/(d+1)}$ and using
\eqref{si:eq:nonlinearresidualequi} gives \eqref{main:eq:nonlinearroot}.
A terminating run has factor zero, so the same conclusion holds.
\end{proof}

\begin{corollary}[Inexact coefficient solves]
\label{si:prop:nonlinearinexact}
The conclusion of \cref{main:thm:nonlinear} also holds for feasible inexact
coefficient solves. Let $\widetilde y_{k+1}$ be the selected combination
and $y^*_{k+1}$ the exact capped minimum over the same retained residuals,
and assume that a prescribed function $\varepsilon(t)\ge0$, independent of
the initialization and of the selected solves, tends to zero at zero and
satisfies $0\le\norm*{\widetilde y_{k+1}}^2-\norm{y^*_{k+1}}^2 \le\varepsilon(T_k)^2T_k^2$
at every nonterminating step, with $T_k$ as in the proof of
\cref{main:thm:nonlinear}. The selected coefficients must still have
sum one and $\ell^1$ norm at most $\Lambda$.
\end{corollary}
\begin{proof}
Convex projection gives
$\norm{\widetilde y_{k+1}-y^*_{k+1}}^2
\le\norm{\widetilde y_{k+1}}^2-\norm{y^*_{k+1}}^2\le\varepsilon(T_k)^2T_k^2$.
The bounds on the new iterate and its nonlinear remainder in the
preceding proof use only coefficient feasibility, so the actual update has
the representation $r_{k+1}=Jy^*_{k+1}+\widehat\xi_k$ with
$\norm{\widehat\xi_k}\le(q\varepsilon(T_k)+\eta(T_k))T_k$. After replacing
this coefficient by a nondecreasing upper envelope near zero, the exact
projections $y^*_{k+1}$ obey all the geometric inequalities used in the
invariance and block arguments, which prove the result with this enlarged
defect modulus.
\end{proof}

\begin{remark}[Uncapped coefficients]\label{si:rem:uncapped}
The cap cannot simply be dropped: for a scalar analytic map whose
derivative has a degree-one certificate, some exact minimizing selections of
ordinary AA($m$), $m\ge2$, diverge from arbitrarily close starts
(\cref{si:prop:unboundedcounterexample}). The conclusion persists without a fixed cap if the selected minimizing coefficients of ordinary AA($m$) have a uniform $\ell^1$ bound while the retained history lies near $x_*$: the proof of \cref{main:thm:nonlinear}, with the unconstrained version of \cref{si:lem:compactblock}, then gives the ordinary factor of \cref{main:thm:generalmargin} or \cref{main:cor:excesscertificate}.
\end{remark}

\end{document}